\documentclass[12pt,a4paper]{article}
\usepackage[english]{babel}
\usepackage[utf8]{inputenc}
\usepackage{amsmath}
\usepackage{amsfonts}
\usepackage{amssymb}
\usepackage{graphicx}
\usepackage{listings}
\usepackage{listingsutf8}
\usepackage{tikz}
\usepackage{cite}
\usepackage{algorithm2e}
\usepackage{array}
\usepackage{makecell}
\usepackage{colortbl}
\usepackage{float}
\usepackage{geometry}
\DeclareMathOperator{\trace}{Tr}

\begin{document}

\allowdisplaybreaks

\setcounter{secnumdepth}{3}
\setcounter{tocdepth}{2}

\nocite{*}
\newcounter{currentcounter}

\newcounter{currentcountercontinuity}
\newcounter{eqcountercontinuity}

\newcounter{algocounter}

\newcounter{mytheoremcounter}
\newcounter{myconjecturecounter}
\newcounter{corollarycounter}

\newcounter{hypocounter}

\newcounter{example}
\newcommand\exc{\theexample}

\newenvironment{general}[1]{\refstepcounter{currentcounter} \bigskip \noindent  \textbf{#1 \thecurrentcounter :}\itshape}{\normalfont \medskip}

\newenvironment{generalbis}[2]{\refstepcounter{currentcounter} \bigskip \noindent \textbf{#1 \thecurrentcounter ~#2 :}\itshape}{\normalfont \medskip}

\newenvironment{proof}{\noindent \textbf{Proof :} \newline \noindent \hspace*{0.2cm}}{\hspace*{\fill}$\square$ \bigskip}

\newenvironment{proofbis}[1]{\noindent \textbf{#1:} \newline \noindent \hspace*{0.2cm}}{\hspace*{\fill}$\square$ \bigskip}

\newenvironment{famous}[1]{\medskip \noindent \textbf{#1:} \newline \noindent  \itshape}{\normalfont \medskip}

\newenvironment{mytheorem}{\refstepcounter{mytheoremcounter} \bigskip \noindent \textbf{Theorem \thechapter.\themytheoremcounter :}\itshape}{\normalfont \medskip}

\newenvironment{mycor}{\medskip \noindent \textbf{Corollary \thecurrentcounter :}  \itshape}{\normalfont \medskip}

\newenvironment{myconjecture}{\refstepcounter{myconjecturecounter} \bigskip \noindent \textbf{Conjecture \themyconjecturecounter :}\itshape}{\normalfont \medskip}

\newenvironment{court}[1]{\refstepcounter{currentcounter} \smallskip \noindent \textbf{#1 \thesubsection.\thecurrentcounter :} \newline \noindent \itshape}{\normalfont \smallskip}

\newenvironment{myalgo}[1]{\refstepcounter{algocounter} \medskip \noindent \textbf{Algorithm \thealgocounter :} #1 \newline \noindent}{\smallskip}

\newcolumntype{C}{>{$}c<{$}}
\newcolumntype{L}{>{$}l<{$}}
\newcolumntype{R}{>{$}r<{$}}
\newcolumntype{T}{>{$ \scriptscriptstyle}r<{$}}
\newcolumntype{S}{>{$ \scriptscriptstyle}c<{$}}

\newcommand\toline{\smallskip \newline}
\newcommand\refp[1]{(\ref{#1})}
\renewcommand\mod{\mathrm{~mod~}}
\newcommand\pregcd{\mathrm{gcd}}
\renewcommand\gcd[2]{\pregcd(#1,#2)}
\newcommand\prelcm{\mathrm{lcm}}
\newcommand\lcm[2]{\prelcm(#1,#2)}
\newcommand\naturaliso{\cong}
\newcommand\iso{\simeq}
\newcommand\cross[1]{#1^{\times}}
\newcommand\crosslong[1]{\cross{(#1)}}
\newcommand\poldegree[1]{\mathrm{deg}(#1)}
\newcommand\cardinal[1]{\# \left(#1\right)}
\newcommand\cardinalshort[1]{\##1}
\renewcommand\det{\mathrm{det}}
\newcommand\indicator[1]{[#1]}

\newcommand\sign{\mathrm{sign}}
\newcommand\signature{\mathrm{sgn}}

\newcommand\bb[1]{\mathbb{#1}}
\newcommand\kk{\bb{K}}
\newcommand\ok{\mathcal{O}_{\kk}}
\newcommand\zz{\bb{Z}}
\newcommand\qq{\bb{Q}}
\newcommand\rr{\bb{R}}
\newcommand\cc{\bb{C}}
\newcommand\ff{\bb{F}}
\newcommand\of{\mathcal{O}_{\ff}}
\newcommand\goth[1]{\mathfrak{#1}}
\newcommand\zsz[1]{\zz/#1\zz}
\newcommand\zszx[1]{(\zsz{#1})^{\times}}
\newcommand\dualsimple[1]{#1^{\vee}}
\newcommand\dual[1]{(#1)^{\vee}}

\newcommand\dx[1]{\mathrm{d}#1}
\newcommand\ddx[1]{\frac{\mathrm{d}}{\dx{#1}}}
\newcommand\dkdxk[2]{\frac{\mathrm{d}^{#2}}{\dx{#1}^{#2}}}
\newcommand\partialx[1]{\partial #1}
\newcommand\partialdx[1]{\frac{\partial}{\partialx{#1}}}
\newcommand\partialkdxk[2]{\frac{\partial^{#2}}{\partialx{#1}^{#2}}}
\newcommand\integral[3]{\int_{#1}{#2}\dx{#3}}
\newcommand\integralsimple[1]{\integrale{G}{#1}{\lambda}}
\newcommand\sprod[2]{\langle #1, #2 \rangle}

\newcommand\limi[1]{\lim_{#1 \to + \infty}}

\newcommand\mn[2]{\mathrm{M}_{#1}(#2)}
\newcommand\mnz[1]{\mn{#1}{\zz}}
\newcommand\sln[2]{\mathrm{SL}_{#1}(#2)}
\newcommand\slnz[1]{\sln{#1}{\zz}}
\newcommand\gln[2]{\mathrm{GL}_{#1}(#2)}
\newcommand\glnz[1]{\gln{#1}{\zz}}
\newcommand\resc{\mathrm{resc}}
\newcommand\com{\mathrm{com}}

\newcommand\bars[1]{\underline{#1}}
\newcommand\taubar{\bars{\tau}}
\newcommand\qbar{\bars{q}}
\newcommand\sigmabar{\bars{\sigma}}
\newcommand\rbar{\bars{r}}
\newcommand\nbar{\bars{n}}
\newcommand\mbar{\bars{m}}
\newcommand\mubar{\bars{\mu}}
\newcommand\xbar{\bars{x}}
\newcommand\Xbar{\bars{X}}
\newcommand\abar{\bars{a}}
\newcommand\bbar{\bars{b}}
\newcommand\alphabar{\bars{\alpha}}
\newcommand\betabar{\bars{\beta}}
\newcommand\omegabar{\bars{\omega}}
\newcommand\taubarsj[1]{\taubar^{-}(#1)}
\newcommand\onebar{\bars{1}}

\newcommand\hh{\bb{H}}
\newenvironment{psmallmatrix}
  {\left(\begin{smallmatrix}}
  {\end{smallmatrix}\right)}
  
\newcommand\declarefunction[5]{#1 := \begin{cases}
\hfill #2 \hfill & \to ~ #3 \\
\hfill #4 \hfill & \to ~ #5 \\
\end{cases}}

\newcommand\sume{\sideset{}{_e}\sum}

\newcommand\tendstowhen[1]{\xrightarrow[#1]{}}

\newcommand\pending{$\bigskip \newline \blacktriangle \blacktriangle \blacktriangle \blacktriangle \blacktriangle \blacktriangle \blacktriangle \blacktriangle \bigskip \newline \hfill$}

\newcommand\pendingref{\textbf{[?]}}

\newcommand\opc[1]{\mathcal{O}^{+,\times}_{#1}}
\newcommand\opck{\opc{\kk}}
\newcommand\opcf{\opc{\goth{f}}}
\newcommand\eps{\varepsilon}
\newcommand\units{\ok^{\times}}
\newcommand\unitsmodf{\mathcal{O}_{\goth{f}}^{\times}}
\newcommand\unitsf{\of^{\times}}
\newcommand\norm[1]{\mathcal{N}(#1)}
\newcommand\normlr[1]{\mathcal{N}\left(#1\right)}
\newcommand\zetaf{\zeta_{\goth{f}}}

\newcommand\ta{\tilde{a}}
\newcommand\tabar{\tilde{\abar}}
\newcommand\talpha{\tilde{\alpha}}
\newcommand\talphabar{\tilde{\alphabar}}
\newcommand\ts{\tilde{s}}
\newcommand\ttt{\tilde{t}}
\newcommand\tlambda{\tilde{\lambda}}
\newcommand\tc{\tilde{c}}
\newcommand\tcone{\tc_1}
\newcommand\tctwo{\tc_2}
\newcommand\te{\tilde{e}}
\newcommand\teone{\te_1}
\newcommand\tetwo{\te_2}
\newcommand\tmu{\tilde{\mu}}
\newcommand\tA{\tilde{A}}
\newcommand\tildef{\tilde{f}}
\newcommand\matrixA{\mathcal{A}}
\newcommand\tmatrixA{\widetilde{\mathcal{A}}}
\newcommand\matrixB{\mathcal{B}}
\newcommand\tmatrixB{\widetilde{\mathcal{B}}}
\newcommand\lineX{\mathcal{X}}
\newcommand\tlineX{\widetilde{\mathcal{X}}}
\newcommand\tx{\tilde{x}}
\newcommand\tbeta{\tilde{\beta}}

\newcommand\spvgamma{\Gamma_{\goth{f}, \goth{b}, \goth{a}}(\eps; h)}
\newcommand\spvgammasign{\Gamma_{\goth{f}, \goth{b}, \goth{a}}^{\pm}(\eps; h)}
\newcommand\spvgr{G_{r, \goth{f}, \goth{b}, \goth{a}}^{\pm}(u_1, \dots, u_r; h)}
\newcommand\spvgrrho{G_{r, \goth{f}, \goth{b}, \goth{a}}^{\pm}([\eps_{\rho(1)}|\dots|\eps_{\rho(r)}]; h_\rho)}
\newcommand\spvgrcomplete{I_{r,\goth{f}, \goth{b}, \goth{a}}(\eps_1, \dots, \eps_r ; \bars{h}, \bars{\mu}, \bars{\nu})}
\newcommand\spvgrcompletebis{I_{r,\goth{f}, \goth{b}, \goth{a}}(\eps_1, \dots, \eps_r ; \bars{h}, \bars{\mu}, \bars{\nu}, (k_{\rho}, \eps_{\rho})_{\rho})}

\newcommand\spvgammac{\Gamma_{\goth{f}, \goth{b}, \goth{a}}(\eps; h, \sigma_{\cc})}
\newcommand\spvgrc{G_{r, \goth{f}, \goth{b}, \goth{a}}^{\pm}([\eps_{\rho(1)}|\dots|\eps_{\rho(r)}]; h_\rho, \sigma_{\cc})}
\newcommand\spvgrrhoc{G_{r, \goth{f}, \goth{b}, \goth{a}}^{\pm}([\eps_{\rho(1)}|\dots|\eps_{\rho(r)}]; h_\rho)}
\newcommand\spvgrcompletec{I_{r,\goth{f}, \goth{b}, \goth{a}}(\eps_1, \dots, \eps_r ; \bars{h}, \bars{\mu}, \bars{\nu}, \sigma_{\cc})}

\newcommand\ula{u_{L, \goth{a}}}
\newcommand\ulah{u_{L, \goth{a}, \bars{h}}}
    
\newcommand\hombase{\mathrm{Hom}}
\newcommand\myhom[3]{\hombase_{#1}(#2, #3)}    
\newcommand\homlong[3]{\hombase_{#1}\left(#2, #3\right)}   
\newcommand\homlz{\myhom{\zz}{L}{\zz}}
\newcommand\homlprimez{\myhom{\zz}{L'}{\zz}}
\newcommand\homlc{\myhom{\zz}{L}{\cc}}
\newcommand\homlambdaz{\myhom{\zz}{\Lambda}{\zz}}
\newcommand\homlambdac{\myhom{\zz}{\Lambda}{\cc}}
\newcommand\homvq{\myhom{\qq}{V}{\qq}}
\newcommand\zexc{z}
  
\newcommand\zfone{\mathcal{Z}_{\goth{f}}^1}
\newcommand\zfoneplus{\mathcal{Z}_{\goth{f}}^{1, +}}
\newcommand\zfoneminus{\mathcal{Z}_{\goth{f}}^{1, -}}
\newcommand\zf{\mathcal{Z}_{\goth{f}}}
  
\newcommand\fracpart[1]{\left\{#1\right\}}
\newcommand\entirepart[1]{\left\lfloor#1\right\rfloor}  
\newcommand\shortsetseparator[1][]{#1|}
\newcommand\setseparator[1][]{~\shortsetseparator[#1]~}
\newcommand\coeff{\mathrm{coeff}}
\newcommand\normalbone[1]{\left(\left(#1\right)\right)}
\newcommand\tensor{\otimes}
\newcommand{\at}[2][]{#1|_{#2}}
\newcommand\omitvar[1]{\widehat{#1}}
\newcommand\kronecker{\delta}
\newcommand\plgt{x}
\newcommand\cohomd{\partial}
\newcommand\cohomdx{\partial^{\times}}
\newcommand\compset[1]{#1^c}
\newcommand\signdet{\sign\,\det}

\newcommand\copen{c^{\circ}}
\newcommand\cclosed{c}
\newcommand\cdual{c^{\vee}}
\newcommand\copendual{c^{\vee, \circ}}
\newcommand\dirac{\delta}

\newcommand\myspan{\mathrm{Span}}
\newcommand\prespanconvex[1]{\mathcal{C}(#1)}
\newcommand\spanconvex{\prespanconvex{V}}
\newcommand\prespancones[1]{\mathcal{K}(#1)}
\newcommand\spancones{\prespancones{V}}
\newcommand\spanconesrr{\prespancones{V_{\rr}}}
\newcommand\spanconesk[1]{\mathcal{K}^{#1}(V)}
\newcommand\spanconeskw[1]{\mathcal{K}^{#1}(W)}
\newcommand\spanconeskdeux[1]{\mathcal{K}(V, #1)}
\newcommand\spanconesktrois[2]{\mathcal{K}(#1, #2)}
\newcommand\spanqcones{\mathcal{K}_{\qq}(V)}
\newcommand\spanqconesrr{\mathcal{K}_{\qq}(V_{\rr})}
\newcommand\prespanwedges[1]{\mathcal{L}(#1)}
\newcommand\spanwedges{\prespanwedges{V}}
\newcommand\spanqwedges{\mathcal{L}_{\qq}(V)}
\newcommand\spanqwedgesrr{\mathcal{L}_{\qq}(V_{\rr})}

\newcommand\coefficient[3]{\mathrm{coeff}[#2^{#3}]\left(#1\right)}
\newcommand\badposition{\mathrm{(BP)}}
\newcommand\symbolparagraph{\S}

\newcommand\alphajk[2]{\alpha^{(#1)}_{#2}}
\newcommand\djk[2]{d^{(#1)}_{#2}}
\newcommand\vjkl{v^{(j)}_{k,j'}}
\newcommand\yj[1]{y^{(#1)}}
\newcommand\gammaj[1]{\gamma^{(#1)}}
\newcommand\yjk[2]{y^{(#1)}_{#2}}
\newcommand\ujk[2]{u^{(#1)}_{#2}}
\newcommand\Ujk[2]{U^{(#1)}_{#2}}

\newcommand\disjointunion{\sqcup}
\newcommand\dkxs{d_k(x, s)}

\newcommand\hexp{h_e}

\newcommand\smoothedgr{G_{n-2, a_1, \dots, a_{n-1}}(v)(w, x, L, L')}
\newcommand\smoothedgromitj{G_{n-2, a_1, \dots, \omitvar{a_j}, \dots, a_{n}}(v)(w, x, L, L')}
\newcommand\smoothedbn{B_{n, a_1, \dots, a_n}(v)(w,x, L, L')}
\newcommand\congruencegrouphnn{\Gamma_0(N, n)}
\newcommand\generaldenomclasscohom[1]{\mathcal{D}(#1)}
\newcommand\denomclasscohom{\generaldenomclasscohom{N, n}}
\newcommand\dNn[2]{\generaldenomclasscohom{#1, #2}}

\newcommand\cun[1]{c^1_{#1}}
\newcommand\cdeux[1]{c^2_{#1}}
\newcommand\Cun[1]{C^{1}_{#1}}
\newcommand\Cdeux[1]{C^{2}_{#1}}
\newcommand\fun{f^1}
\newcommand\fdeux{f^2}
\newcommand\mathcalCpm[2]{\mathcal{C}^{#1}_{#2}}
\newcommand\mathcalCplus[1]{\mathcal{C}^{+}_{#1}}
\newcommand\mathcalCmoins[1]{\mathcal{C}^{-}_{#1}}
\newcommand\good{\textit{good }}
\newcommand\wellplaced{well placed }
\newcommand\hj[2]{H^{#1}_{#2}}
\newcommand\hjplus[1]{\hj{+}{#1}}
\newcommand\hjminus[1]{\hj{-}{#1}}

\newcommand\preclassgroup{\mathrm{Cl}}
\newcommand\preclassgroupplus{\preclassgroup^{+}}
\newcommand\classgroup[1]{\preclassgroup(#1)}
\newcommand\classgroupplus[1]{\preclassgroupplus(#1)}
\newcommand\classgroupk{\classgroup{\kk}}
\newcommand\classgroupplusk{\classgroupplus{\kk}}
\newcommand\classgroupf{\classgroupplus{\goth{f}}}
\newcommand\wideclassgroupf{\classgroup{\goth{f}}}
\newcommand\classfieldf{\kk^{+}(\goth{f})}

\newcommand\BL{B_L}
\newcommand\BLambda{B_{\Lambda}}
\newcommand\Bok{\tilde{B}}

\newcommand\lemmasignstranspose{lemma \ref{lemmasignrelations}, (ii) }
\newcommand\lemmasignstoprows{lemma \ref{lemmasignrelations}, (i) }
\newcommand\lemmasignsujk{lemma \ref{lemmasignrelations}, (iii) }

\newcommand\cohendifferent{[\hspace{1sp}\cite{cohen}, Proposition 4.8.19]}
\newcommand\cohenhnf{[\hspace{1sp}\cite{cohen}, Theorem 2.4.3]}
\newcommand\cohensnf{[\hspace{1sp}\cite{cohen}, Theorem 2.4.12]}

\newcommand\matrixgothu{\goth{U}}
\newcommand\tmatrixgothu{\widetilde{\goth{U}}}
\newcommand\matrixm{\mathcal{M}}
\newcommand\tmatrixm{\widetilde{\mathcal{M}}}

\newcommand\coefH{R}
\newcommand\normeps{\norm{\eps}}
\newcommand\chif{\chi_{\goth{f}}}
\newcommand\helper{\goth{H}}
\newcommand\helpers{\underline{\helper}}
\newcommand\helperdown{p_{\helper}}
\newcommand\DaL{\goth{D}(a, L)}
\newcommand\Dtaok{\goth{D}(\ta, \ok)}
\newcommand\tildeD{\widetilde{\goth{D}}}

\newcommand\sigmaplus{\sigma}
\newcommand\conjko{\kappa}

\newcommand\unitsnormf{E_{\goth{f}}}
\newcommand\auxomega{\Omega}

\newcommand{\rvline}{\hspace*{-\arraycolsep}\vline\hspace*{-\arraycolsep}}

\def\Xint#1{\mathchoice
   {\XXint\displaystyle\textstyle{#1}}%
   {\XXint\textstyle\scriptstyle{#1}}%
   {\XXint\scriptstyle\scriptscriptstyle{#1}}%
   {\XXint\scriptscriptstyle\scriptscriptstyle{#1}}%
   \!\int}
\def\XXint#1#2#3{{\setbox0=\hbox{$#1{#2#3}{\int}$}
     \vcenter{\hbox{$#2#3$}}\kern-.5\wd0}}
\def\mulint{\Xint\times}
\newcommand\intmult{\mulint}

\newcommand\Upsilonf{\Upsilon_{\goth{f}}}
\newcommand\adapted{\textit{adapted }}

\newcommand\uzf{U_{\goth{f}}}
\newcommand\uzfplus{\uzf^+}
\newcommand\clzf{\mathrm{Cl}(\goth{f}, \zz)}
\newcommand\clzfplus{\mathrm{Cl}^{+}(\goth{f}, \zz)}
\newcommand\klzfplus{\kk^{+}(\goth{f}, \zz)}

\begin{center}
\end{center}
\begin{center}
\large \noindent Geometric families of multiple elliptic Gamma functions and arithmetic applications, III

\bigskip \medskip \noindent \textit{Higher elliptic units}
\end{center}
\medskip
\begin{center}
Pierre L. L. Morain\footnotemark[1]\footnotetext[1]{Sorbonne Université and Université Paris Cité, CNRS, INRIA, IMJ-PRG, F-75005 Paris, France. This article was written during the PhD of the author which was funded by the École polytechnique, Palaiseau, France.}
\end{center}
\medskip
\begin{center}
\noindent \textbf{Abstract:}
\end{center}
This is the third paper in a series where we study arithmetic applications of the higher elliptic Gamma functions. In the first two articles in this series, we defined geometric families of these functions and proved that they satisfy coboundary and cocycle relations under the action of special linear groups and associated congruence subgroups. The main purpose of the present paper is to present a construction of conjectural higher elliptic units above number fields with exactly one complex place as special values of higher elliptic Gamma functions, upgrading the construction carried out by Bergeron, Charollois and Garc\'ia for complex cubic fields. These higher elliptic units are obtained by evaluating a multiplicative $(n-2)$-cocycle built from higher elliptic Gamma functions against a $(n-2)$-cycle associated to some group of totally positive units of a given number field of degree $n$ with exactly one complex place. We conjecture that these higher elliptic units are algebraic units which belong to prescribed abelian extensions of the base field where they are evaluated and that they satisfy a Kronecker limit formula which relates the logarithm of their modulus to values of derivatives of partial zeta functions at $s = 0$ in the base field. We showcase our conjecture on various examples for number fields of degree 3, 4 and 5.

\normalsize \bigskip \bigskip

\noindent \textbf{Acknowledgments:} This article is based on chapters III and IV of the author's PhD thesis \cite{thesis} at Sorbonne Université. The author would like to thank his advisors Pierre Charollois and Antonin Guilloux at Sorbonne Université for their guidance, their support and the many discussions on this project around higher elliptic units. The author is also grateful to Mateo Crabit Nicolau and Mart\'i Roset Julià for valuable conversations on related topics in number theory. Finally, the author wishes to thank Henri Darmon and Aurel Page for reviewing the author's PhD thesis manuscript, upon which this article expands.  

\bigskip \bigskip

\noindent \textbf{AI disclosure:} All material in the paper was written by the author without AI assistance.

\newpage

\tableofcontents

\clearpage

\section{Introduction}\label{introduction}

In this series of papers we are interested in a hierarchy of functions satisfying modular transformation properties for special linear groups and whose special values are conjectured to yield algebraic units. In the first two articles of this series, we focused on the cohomological properties of these functions. In the present paper, we focus on the construction of conjectural higher elliptic units above number fields with exactly one complex place as special values of these functions. This work is inspired by the construction of elliptic units above imaginary quadratic fields in the context of the theory of \textit{Complex multiplication}. The elliptic units are obtained as special values of the $\theta$ function defined on $\cc \times \hh$ by:
$$\theta(z, \tau) = \prod_{n \geq 0} \left(1- e^{-2i\pi z} e^{2i \pi (n+1) \tau}\right) \left(1-e^{2i\pi z}e^{2i\pi n \tau}\right).$$
If we introduce the (infinite) $q$-Pochhammer symbol defined for $|q|<1$ by the infinite product $(x; q)_{\infty} = \prod_{n \geq 0}(1 - xq^n)$ then $\theta(z, \tau)$ is simply given by:
$$\theta(z, \tau) = (e^{2i\pi(\tau-z)}\,;\, e^{2i\pi\tau})_{\infty}\cdot(e^{2i\pi z}\,;\, e^{2i\pi \tau})_{\infty}.$$
The elliptic units are given by evaluations of the form:
\begin{equation}\label{ellipticunitsintro}
\frac{\theta\left(\frac1q, \tau\right)^N}{\theta\left(\frac{N}{q}, N\cdot \tau\right)}
\end{equation}
for well-chosen $q, N \in \zz_{\geq 2}$ and $\tau \in \kk - \qq$. These evaluations yield $q$-units in abelian extensions of $\kk$ which are $N$-smoothed versions of Ramachandra's elliptic units (see \cite{Ramachandra}). They also satisfy a version of Kronecker's second limit formula (see section \ref{reviewquadratic} for general results on elliptic units).  

In 2023, a construction of conjectural higher elliptic units above complex cubic fields was formulated by Bergeron, Charollois and Garc\'ia \cite{BCG}. These higher elliptic units are obtained as special values of the elliptic Gamma function of Ruijsenaars \cite{Ruijsenaars}.
It is a meromorphic function of three variables defined on $\cc\times\hh^2$ by the infinite product: 
$$\Gamma(z, \tau, \sigma) = \prod_{m,n \,\geq 0}\frac{1-e^{2i\pi((m+1)\tau + (n+1)\sigma - z)}}{1-e^{2i\pi(m\tau+n\sigma+z)}}.$$
This function may be expressed in terms of multiple $q$-Pochhammer symbols, similarly to how the $\theta$ function could be expressed in terms of $q$-Pochhammer symbols. Indeed, if $|q|, |r| < 1$ we set $(x; q, r)_{\infty} = \prod_{m, n \geq 0} (1 - xq^m r^n)$. Using this symbol, one gets the expression:
$$\Gamma(z, \tau, \sigma) = \frac{(e^{2i\pi(\tau+\sigma-z)}\, ;\, e^{2i\pi\tau}\,,\,e^{2i\pi\sigma})_{\infty}}{(e^{2i\pi z}\, ;\, e^{2i\pi\tau}\,,\,e^{2i\pi\sigma})_{\infty}}.$$
The type of evaluation carried out in \cite{BCG} is an $N$-smoothed evaluation similar to \refp{ellipticunitsintro} taking the shape:
\begin{equation}\label{ellipticunitscubicintro}
\prod_{z \in F} \frac{\Gamma\left(z, \tau, \sigma\right)^N}{\Gamma\left(N\cdot z, N\cdot \tau, N\cdot \sigma\right)}
\end{equation}
for well-chosen $\tau, \sigma \in \kk$ and a well-chosen finite set $F \subset \kk$. Similarly to elliptic units for imaginary quadratic fields, these evaluations satisfy a Kronecker limit formula for complex cubic fields (see \cite{BCG}). We shall review this construction in section \ref{reviewcubic}. 

The central part of the author's PhD thesis concerns the generalisation of such constructions of (higher) elliptic units to arbitrary number fields with exactly one complex place (so-called ``Almost Totally Real number fields'', or ATR fields for short) using the hierarchy of higher elliptic Gamma functions we introduced in this series of papers. Recall that the multiple elliptic Gamma functions $G_r(z, \tau_0, \dots, \tau_r)$ introduced by Nishizawa in \cite{Nishizawa} are meromorphic functions on $\cc \times \hh^{r+1}$ defined by the infinite product:
\begin{equation}\label{defgr}
G_r(z, \tau_0, \dots, \tau_r) = \prod_{m_0, \dots, m_r \,\geq 0} \left(1 - e^{2i\pi(-z + \sum_{j = 0}^r (m_j+1)\tau_j)}\right)  \left(1- e^{2i\pi (z + \sum_{j = 0}^r m_j\tau_j)}\right)^{(-1)^r}.
\end{equation}
The identifications $G_0 = \theta$ and $G_1 = \Gamma$ are clear. These functions may be rewritten in terms of multiple $q$-Pochhammer symbols (see [\hspace{1sp}\cite{Nishizawa}, Definition 3.1]). Indeed, for $|q_0|, \dots, |q_r| <1$ set $(x; q_0, \dots, q_r)_{\infty} = \prod_{m_0, \dots, m_r \geq 0} (1 - xq_0^{m_0}\dots q_r^{m_r})$. Then the $G_r$ function is defined by:
$$G_r(z, \tau_0, \dots, \tau_r) = (e^{2i\pi(\tau_0 + \dots + \tau_r -z)}\,;\, e^{2i\pi\tau_0}, \dots, e^{2i\pi\tau_r})_{\infty}\cdot (e^{2i\pi z}\,;\, e^{2i\pi\tau_0}, \dots, e^{2i\pi\tau_r})_{\infty}^{(-1)^r}.$$ 
For an ATR field of degree $n = r+2$ the general shape of higher elliptic units is a product of the form:
\begin{equation}\label{ellipticunitsgeneralintro}
\prod_{j = 1}^m \frac{G_{n-2}(z_j, \tau_{1, j}, \dots, \tau_{n-1,j})^N}{G_{n-2}(Nz_j, N\tau_{1, j}, \dots, N\tau_{n-1,j})}
\end{equation}
where the evaluation parameters $z_j, \tau_{1, j}, \dots, \tau_{n-1, j}$ should be carefully chosen elements in $\kk$. The construction of higher elliptic units consists of the explicit description of evaluation parameters for which evaluations of the form \refp{ellipticunitsgeneralintro} give algebraic numbers satisying a Kronecker limit formula. In \cite{cheu}, we briefly explained this construction in optimal settings. In this paper, we present the general construction carried out in \cite{thesis} under some simplifying hypotheses (see sections \ref{sectiongcheu} and \ref{sectiongeomsetup}). It is our goal to expand the construction to any ATR field setting. The general conjecture on these special values of higher elliptic Gamma functions (which is precisely stated in section \ref{sectionconjecture}, Main Conjecture \ref{conjecture}) takes the following vague form:

\begin{generalbis}{Conjecture}{[Vague form 1]}\label{conjecturevagueformun}
Let $\kk$ be an ATR field of degree $n \geq 3$ which we view as a subset of $\cc$ via one of its two complex embeddings. Fix an integral ideal $\goth{f} \neq \ok$. Fix an integral ideal $\goth{b}$ of $\ok$ representing a class in $\classgroupf$ and an auxiliary integral ideal $\goth{a}$ of $\ok$ of norm $N$ such that $\ok/\goth{a}$ is cyclic. Then there is an explicit integer $m$ and explicit parameters $z_j$, $\tau_{1, j}, \dots, \tau_{n-1,j}$ in $\kk$, all depending on $\goth{f}, \goth{b}, \goth{a}$ such that the evaluation
\begin{equation}\label{uvagueformun}
u = \left(\prod_{j = 1}^{m} \frac{G_{n-2}(z_j, \tau_{1, j}, \dots, \tau_{n-1,j})^N}{G_{n-2}(Nz_j, N\tau_{1, j}, \dots, N\tau_{n-1,j})}\right)^{\denomclasscohom}
\end{equation}
is an algebraic unit in the narrow ray class field $\classfieldf$ (where the integer $\denomclasscohom$ is defined by $\denomclasscohom = \prod_{p \mid N} p^{\lfloor \frac{n}{p-1} \rfloor} $). In addition, the algebraic unit $u$ satisfies a Kronecker limit formula of the form:
\begin{equation}\label{klfvagueformun}
N \cdot \zeta'_{\goth{f}}(\goth{b}, 0) - \zeta'_{\goth{f}}(\goth{a}\goth{b}, 0) = \frac{1}{\denomclasscohom} \log|u|^2.
\end{equation}
\end{generalbis}

We illustrate this vague form of our main conjecture with a quartic example. Consider $\kk = \qq(z)$ where $z$ is the root of the polynomial $x^4 - x^3 - x^2 - 5x + 1$ lying in the upper half-plane. Let $\goth{f}$ be the unique prime ideal of norm $q = 7$ in $\ok$. The narrow ray class group at $\goth{f}$ has order $6$ and is represented by the integral ideals $\ok, 2\ok, 3\ok, \goth{P}_{23}, 2\goth{P}_{23}$ and $3\goth{P}_{23}$ where $\goth{P}_{23}$ is the unique prime ideal of norm $23$ in $\ok$. Let us fix the smoothing ideal $\goth{a}$ to be the unique prime ideal of norm $N = 11$ in $\ok$. Our construction yields the evaluation parameters:
\begin{align*}
&\tau = \frac{7z^{3}-  13z^{2}-  15z-  19098}{5159}, && \tau' = \frac{-5z^{3}+  4z^{2}+  6z-  2607}{77} \\
&\sigma = \frac{9z^{3}+  12z^{2}+  19z+  38971}{5159}, && \sigma' = \frac{-z^{3}+  z^{2}+  z-  477}{77} \\
&\rho = \frac{3z^{3}+  4z^{2}-  16z-  15641}{5159}, && \rho' = \frac{2z^{3}-  2z^{2}-  3z+  827}{77}.
\end{align*}
For $k = 1, 2, 3, 4, 5, 6$, write $I_k$ for the value
$$I_k = \frac{G_2\left(\frac{k}{7}, \tau, \sigma, \rho\right)^{11}}{G_2\left(\frac{11k}{7}, 11\tau, 11\sigma, 11\rho\right)} \times \frac{G_2\left(\frac{3k}{7}, \tau', \sigma', \rho'\right)^{11}}{G_2\left(\frac{33k}{7}, 11\tau', 11\sigma', 11\rho'\right)}.
$$
We may compute these values to high precision (say 1000 digits) and compute in particular:
\begin{center}
\begin{tabular}{| C | C | C |}
\hline
\text{Class} & \text{Definition of } u & \text{Value of } u \\
\hline
\goth{b}_1 = \ok & u = I_1I_6 & u \approx 204411.4289667...+ i\cdot1113079.0424875...\\
\hline
\goth{b}_2 = 2\ok & u = I_2I_5 & u \approx 0.02745587... - i\cdot0.01233068... \\
\hline
\goth{b}_3 = 3\ok & u = I_3I_4 & u \approx (1.6667587... - i\cdot 2.4168791...) \cdot 10^{-5}\\ 
\hline
\end{tabular}
\end{center}
The complex numbers $I_1I_6, I_2I_5, I_3I_4$ together with their inverses seem to be algebraic units as they are found to coincide up to 1000 digits of precision with six of the roots of the reciprocal polynomial:
\begin{align*}
P & = x^{24}+1 - 447564(x^{23}+x) + 1297730000060(x^{22}+x^2) \\
&- 50092339028926216(x^{21} + x^3) + 1489197308166747656846(x^{20} + x^4) \\
&- 98893338471699871421925(x^{19}+x^5) + 2195426185183831598917668(x^{18}+x^6)\\
& - 11183927451691276697810187(x^{17}+x^7) + 29721644381048971258757251(x^{16} + x^8)\\
&- 53392984783876088115872757(x^{15} + x^{9}) + 76166042306689872075863788(x^{14} + x^{10})\\
& - 94138004635629304533913859(x^{13} + x^{11}) + 101459141139453330556583357x^{12}
\end{align*}
which defines an absolute equation of $\classfieldf$ over $\qq$. Alternatively, we may identify their expected minimal relative polynomial over $\kk$ as:
\begin{align*}
\prod_{k = 1, 2, 3} (X - & I_k I_{7-k})(X-I_k^{-1} I_{7-k}^{-1}) \approx \\
&X^6 + 1 +(-66007z^3 + 315205z^2 - 404670z + 66406) (X^5+X) \\ + & (-2980477827z^3 - 882760955z^2 + 19370690773z - 3656696288) (X^4 + X^2) \\ +& (153504704659z^3 - 136634711339z^2 - 544162307740z+ 108155656031) X^3.
\end{align*}
We may check once again that the polynomial in the right-hand side defines a relative equation of $\classfieldf$ over $\kk$. One may then check the expected Kronecker limit formula \refp{klfvagueformun} as up to 1000 digits:
$$11\cdot \zeta'_{\goth{f}}(\ok, 0) - \zeta'_{\goth{f}}(\goth{a}, 0) \approx \log\left|I_1I_6\right|^2 \approx 27.8784505...$$
and similarly for the remaining 5 classes. We will give more details on how this example was produced using our general construction in section \ref{sectionconjecture}.

We highlight that, if proven, our conjectures on higher elliptic units in the flavour of Conjecture \ref{conjecturevagueformun} would constitute progress towards a solution to Hilbert's 12th problem (see \cite{Hilbert12}) for ATR fields, where the higher elliptic Gamma functions are used to describe the abelian extensions of such number fields (see section \ref{hilbertsection}). Our conjectures also have ties with the ``rank one abelian Stark conjectures'' (see \cite{Stark}, \cite{rankonestark}), and one may view our higher elliptic units as smoothed versions of the conjectural Stark units (see section \ref{starksection}). At the time of writing this paper, Hilbert's 12th problem is solved only when the base field is either $\qq$ or an imaginary quadratic field. Very recently, some progress has been made by Radchenko and Wheeler \cite{RW} towards a solution to Hilbert's 12th problem for real quadratic fields. Indeed, they proved that the conjectural Stark units above real quadratic fields are algebraic numbers which can be expressed in terms of evaluations of Faddeev's modular quantum dilogarithm. This modular function may be expressed in terms of $q$-Pochhammer symbols:
$$\Phi_{\gamma, m, n}(z, \tau) = \frac{(e^{2i\pi(m\tau+z)}\,;\, e^{2i\pi\tau})_{\infty}}{(e^{2i\pi\frac{n(a\tau+b)+z}{c\tau+d}}\,;\, e^{2i\pi\frac{a\tau+b}{c\tau+d}})_{\infty}}  $$
for $\gamma  = \begin{psmallmatrix} a & b \\ c & d \end{psmallmatrix} \in \slnz{2}$, $m, n \in \zz$ and $z, \tau \in \cc \times (\cc - \rr)$ where we explicitly use the usual action by Möbius transformation $\gamma \cdot (z, \tau) = (\frac{z}{c\tau +d}, \frac{a\tau +b}{c\tau +d })$. In \cite{RW}, it is shown that evaluations of this function $\Phi$ at well-chosen points in a real quadratic field yield algebraic numbers. The $q$-Pochhammer symbols and their higher dimension analogues we use in this paper enjoy nice transformation properties under actions of corresponding special linear groups. Thus, it seems that they are the correct building blocks one may use to construct analytic functions which produce collections of algebraic numbers when evaluated in some specific way at points in some fixed number field. 

The main goal of this paper is to formulate a precise form of Conjecture \ref{conjecturevagueformun} (see the Main Conjecture \ref{conjecture}) for general ATR fields and to provide a general construction of higher elliptic units. This is done by giving an explicit description of \refp{uvagueformun} as a cohomological evaluation of a multiplicative $(n-2)$-cocycle built from $N$-smoothed $G_{n-2}$ on subgroups of $\Gamma_{0}(N, n)$ against the fundamental $(n-2)$-cycle associated to some group of totally positive units of $\ok$. This is coherent with the fact that the rank of the group of units in $\ok$ is exactly $n-2$ by Dirichlet's unit theorem. For this purpose we shall use the results on the cohomological properties satisfied by the higher elliptic Gamma functions which we proved in the first two paper in this series (see \cite{firstpaper} and \cite{secondpaper}). 

In practice, the fundamental cycles on which we evaluate take the form 
\begin{equation}\label{shapecycleintro}\Upsilon = \sum_{\rho \in \goth{S}_{n-2}} \Upsilon(\rho) = \sum_{\rho \in \goth{S}_{n-2}} \signature(\rho) \left(\eps_{\rho(1)}, \dots, \prod_{i = 1}^{j} \eps_{\rho(i)}, \dots, \prod_{i = 1}^{n-2} \eps_{i}\right)
\end{equation} 
where the components are indexed by the symmetric group on $n-2$ elements and $\eps_1, \dots, \eps_{n-2}$ is a choice of fundamental units for some group of totally positive units in $\ok$. In section \ref{sectionreview}, we describe a procedure which to some component $\Upsilon(\rho)$, some choice of base point $h_{\rho} \in \kk$ and some choice of ``orientation datum'' associates a complex number $I(\rho)$ given by an explicit product of $G_{n-2}$ functions evaluated at points in $\kk$. The complex number $u$ appearing in \refp{uvagueformun} is then given by $u = (\prod_{\rho} I(\rho))^{\denomclasscohom}$ when the base points and the orientation data are \textit{well chosen}. We shall give an explicit description of what it means for these to be well chosen in terms of arithmetic quantities and properties of the cycle $\Upsilon$. This will lead us to define the notion of compatible sets of admissible base points, for which we expect our construction of higher elliptic units to be correct. This notion of compatibility is expressed as the equality of some special ideals attached to the cycle $\Upsilon$ and to the base points $(h_{\rho})_{\rho}$. Namely, let us fix the data $\kk, \goth{f}$, $\goth{b}$, $\goth{a}$ as in Conjecture \ref{conjecturevagueformun} and a fundamental cycle of the form given in \refp{shapecycleintro} associated to some fundamental units $\eps_1, \dots, \eps_{n-2}$. Let $q$ be the unique positive integer satsifying $q\zz = \goth{f} \cap \zz$. Let us set $L = \goth{f}\goth{b}^{-1}$. Let us fix a $\zz$-basis $B$ of $L$ and for each $\rho \in \goth{S}_{n-2}$ let us fix a base point $h_{\rho} \in L$. For each permutation $\rho$ let us consider the linear form $f_{\rho}$ defined on $\kk$ by:
$$f_{\rho}: y \mapsto \det_{B}\left(h_{\rho}, \eps_{\rho(1)}h_{\rho}, \dots, \left(\prod_{i = 1}^{j}\eps_{\rho(i)}\right)h_{\rho}, \dots, \left(\prod_{i = 1}^{n-2} \eps_{\rho(i)}\right)h_{\rho}, y\right)$$
and let us define $a_{\rho}$ to be the unique linear form in the one-dimensional $\qq$-vector space $\qq f_{\rho}$ satisfying $a_{\rho}(L) =\zz$. We may associate to this linear form a generalised different ideal $\goth{D}(a_{\rho}, L)$ defined by:
$$\goth{D}(a_{\rho}, L)^{-1} := \{x \in \kk \setseparator \forall\, y \in L, a_{\rho}(xy) \in \zz\}.$$
We shall say that the base points $h_{\rho}$ for $\rho \in \goth{S}_{n-2}$ form a compatible set of admissible base points (for the data $\goth{f}, \goth{b}, \goth{a}, \eps_1, \dots, \eps_{n-2}$ if:
\begin{enumerate}
\item \text{[Admissibility]} For each $\rho \in \goth{S}_{n-2}$, $h_{\rho}/q \equiv 1 \mod L$ and $h_{\rho}/N$ generates the cyclic group $\goth{a}^{-1}L/L$.
\item \text{[Compatibility]} The generalised different ideal $\goth{D}(a_{\rho}, L)$ is independent of the permutation $\rho$.
\end{enumerate}
Since our construction of higher elliptic units relies on the existence of such compatible sets, we devote most of section \ref{sectiongeomsetup} to the explicit construction of compatible sets for which we have a good control on the size of the product \refp{uvagueformun}. A vague form of the main result of this paper (see Theorem \ref{theoremchoiceh}) is given as follows:

\begin{generalbis}{Theorem}{[Vague form]}
Suppose that $\kk$ is an ATR field of degree $n = r+2 \geq 3$. Fix an integral ideal $\goth{f} \neq \ok$. Fix an integral ideal $\goth{a}$ of norm $N$ coprime to $q$ such that $\ok/\goth{a}$ is cyclic. Fix an integral ideal $\goth{b}$ coprime to $\goth{f}$. Fix a set of fundamental units $\eps_1, \dots, \eps_r$ for some group of totally positive units in $\ok$ and consider the cycle $\Upsilon$ attached to these units. Under mild assumptions on $\goth{f}, \goth{b}, \goth{a}$ and $\Upsilon$, there is a compatible set of admissible base points $(h_{\rho})_{\rho}$ for the data $\goth{f}, \goth{b}, \goth{a}, \Upsilon$. Moreover, the size of the product \refp{uvagueformun} associated to this choice of base points is bounded by quantities depending only on $\goth{f}$ and $\Upsilon$.
\end{generalbis}

Let us now give an outline of this paper. In section \ref{sectionreview} we review the construction of (higher) elliptic units for imaginary quadratic fields and complex cubic fields. In section \ref{sectiongcheu} we describe precisely the construction of higher elliptic units for general ATR fields. The main conjecture is expressed up to choices, the most important of which is the choice of a compatible set of base points. This choice is at the heart of the construction and it is discussed at length in section \ref{sectiongeomsetup} where we make the construction as explicit as possible. Finally, in section \ref{sectioncomputing} we briefly explain how we put our conjecture to the test, and in section \ref{sectionnumerical} we give numerical evidence to support the conjecture for number fields of degree 3, 4 and 5.

\section{Review of the quadratic and cubic cases}\label{sectionreview}

In this section, we briefly review the main results on elliptic units above for imaginary quadratic fields as well as the construction of conjectural elliptic units above complex cubic fields by Bergeron, Charollois and Garc\'ia.

\newcommand\taustar{\tau^{\ast}}
\newcommand\zstar{z^{\ast}}

\subsection{Elliptic units for imaginary quadratic fields}\label{reviewquadratic}

Recall the definition of the $\theta$ function for $(z, \tau) \in \cc \times \hh$ by:
$$\theta(z, \tau) = \prod_{n \geq 0} \left(1- e^{-2i\pi z} e^{2i \pi (n+1) \tau}\right) \left(1-e^{2i\pi z}e^{2i\pi n \tau}\right).$$
This function enjoys modular transformation properties under an action of $\slnz{2}$. Indeed, there is a rational function $P_{2} : \slnz{2} \to \qq(z, \tau)$ such that for all $\gamma = \begin{psmallmatrix} a & b \\ c & d \end{psmallmatrix}$:
\begin{equation}\label{dedekindrademacher}
\theta\left(\frac{z}{c\tau +d}, \frac{a\tau + b}{c\tau + d}\right) = \theta(z, \tau)\cdot e^{2i\pi P_{2, \gamma}(z, \tau)}.
\end{equation}
The rational function $P_{2, \gamma}$ may be expressed in terms of classic Dedekind sums (see \cite{Rademacher}). 

It is classical to perform a smoothing operation on the $\theta$ function to obtain modular units of one of the two following shapes:
$$\frac{\theta(z, \tau)^{N^2}}{\theta(Nz, \tau)}~~\textit{ or } ~~ \frac{\theta(z, \tau)^N}{\theta(Nz, N\tau)}.$$
In the first case, by adding a small exponential prefactor we obtain the basic ingredient for Siegel units:
\begin{equation}\label{prefacteursiegel}
 _N\theta(z, \tau) = e^{2i\pi\left(\frac{N^2-1}{12}\tau + \frac{N - N^2}{2}\left(z-\frac{1}{2}\right)\right)}\frac{\theta(z, \tau)^{N^2}}{\theta(Nz, \tau)}
 \end{equation}
when $z \in \qq$, $\tau$ is an imaginary quadratic number and $N$ is a rational prime which is inert in the field $\qq(\tau)$ (see [\hspace{1sp}\cite{Kato}, Proposition 1.3]). The second expression gives smoothed versions of Ramachandra's elliptic units (see \cite{Ramachandra}) in the setting where $N$ is a rational prime which splits in $\qq(\tau)$. In this work we focus on the second type of smoothing operation. 

We now briefly review the two main results on these smoothed elliptic units, that is a smoothed version of Kronecker's second limit formula they satisfy and the fact that they are algebraic $S$-units (see \cite{Siegel}, \cite{Stark} and \cite{Ramachandra}, \cite{Robert} for versions of these results). We package these results together in the following form:

\begin{generalbis}{Theorem}{(Kronecker, Siegel-Ramachandra)}\label{theoremklf}
Let $\kk$ be an imaginary quadratic field. Fix an integral ideal $\goth{f} \neq \ok$. Set $q\zz = \goth{f} \cap \zz$ for some positive integer $q$. Fix an embedding of the ray class field at $\goth{f}$ in $\cc$. Fix an integral ideal $\goth{b}$ representing a class in the class group at $\goth{f}$. Let $\goth{a}$ be an integral ideal of prime norm $N$ such that $(N, q) = 1$. Let $\mu, \nu$ be elements of $\kk$ such that $\goth{f}\goth{b}^{-1} = \zz \mu  + \zz \nu$ and $\goth{f}(\goth{a}\goth{b})^{-1} = \zz \mu +\zz \nu/N$ with $\mu/\nu = \tau \in \kk \cap \hh$. Assume further that $\nu/q \equiv 1 \mod \goth{f}\goth{b}^{-1}$. Then the complex number
$$u_{\goth{b}} = \frac{\theta\left(\frac{1}{q}, \tau\right)^N}{\theta\left(\frac{N}{q}, N\tau\right)}$$
is a $q$-unit in the ray class field of $\kk$ at $\goth{f}$. The Shimura reciprocity law is explicitly given by $\sigma_{\goth{b}'}(u_{\goth{b}}) = u_{\goth{b}\goth{b}'}$ where $\goth{b}' \mapsto \sigma_{\goth{b}'}$ is Artin's map. In addition, $u_{\goth{b}}$ satisfies a version of Kronecker's second limit formula:
$$N \cdot \zeta'_{\goth{f}}(\goth{b}, 0) - \zeta'_{\goth{f}}(\goth{a}\goth{b}, 0) = -\frac{1}{w_{\goth{f}}}\log\left|\frac{\theta\left(\frac{1}{q}, \tau\right)^N}{\theta\left(\frac{N}{q}, N\tau\right)}\right|^2$$
where $w_{\goth{f}}$ is the number of roots of unity in $\ok$ which are congruent to $1 \mod \goth{f}$. 
\end{generalbis}

We may illustrate this formula on a simple example for the field $\kk = \qq(\sqrt{-19})$. Fix $z = i \cdot \sqrt{19}$ in the upper half-plane. Consider $\goth{f} = \frac{z-1}{2}\ok$ a prime ideal of norm $q = 5$. Fix $\goth{b} = \ok$ and $\goth{a} = \frac{z-7}{2}$ a prime ideal of norm $N = 17$. The elements $\mu = (z+9)/2$ and $\nu =-(5z+35)/2 $ of $\kk$ are such that $\goth{f}\goth{b}^{-1} = \zz \mu + \zz \nu$ and $\goth{f}(\goth{a}\goth{b}^{-1}) = \zz\mu + \zz (\nu/17)$ with $\nu/5 \equiv 1 \mod \goth{f}$. Setting $\tau = \mu/\nu = \frac{z-41}{170}$, we may compute: 
$$u = \frac{\theta\left(\frac{1}{5}, \tau\right)^{17}}{\theta\left(\frac{17}{5}, 17\tau\right)}$$
and check that 
$$u = \frac{1160\sqrt{-19}-2185+\sqrt{-5068642\sqrt{-19} - 20792753}}{2}$$
is a root of the irreducible polynomial:
$$x^4 + 4370x^3 + 30340914x^2 + 6780625x + 390625 $$
which defines an absolute equation of $\classfieldf$ over $\qq$. In addition, we may check the Kronecker limit formula as:
$$17\cdot \zeta'_{\goth{f}}(\ok, 0) - \zeta'_{\goth{f}}(\goth{a}, 0)  =  -\log|u|^2 \approx 4.3524721... $$

In the rest of this work, we shall be interested in conjectural analogues of elliptic units for ATR fields of higher degree and conjectural generalisations of Theorem \ref{theoremklf}.

\subsection{Review of the complex cubic case}\label{reviewcubic}

Our study of higher elliptic units is based on the treatment of the complex cubic case by Bergeron, Charollois and Garc\'ia in \cite{BCG}. In their article, they give a construction for conjectural elliptic units above a complex cubic field $\kk$ by evaluating the elliptic Gamma function at points in $\kk$. They relate their evaluations to the values of derivatives of partial zeta functions in $\kk$ at $s = 0$, thus proving an analogue for complex cubic fields of Kronecker's second limit formula. In this section, we briefly review their construction and their results which motivated our work.

Let us fix $\kk$ a complex cubic field. Let us denote by $\sigma_{\rr} = \sigma_1$ the real embedding of $\kk$ and by $\sigma_{\cc} = \sigma_2 = \overline{\sigma_3}$ one of the complex embeddings of $\kk$. Fix an integral ideal $\goth{f} \neq \ok$ which will be the finite part of the class field modulus. We denote by $\classfieldf$ the narrow ray class field at $\goth{f}$, which we assume to be totally complex. We denote by $q$ the positive integer satisfying $q\zz = \goth{f} \cap \zz$. Fix an integral ideal $\goth{b}$ of $\ok$ representing a class in $\classgroupf$. Fix an integral ideal $\goth{a}$ of $\ok$ such that $\ok/\goth{a}$ is cyclic of order $N = \norm{\goth{a}}$ coprime to $q\cdot \goth{b}$. Such an ideal $\goth{a}$ is called a smoothing ideal and $N$ is called the smoothing index (see Definition \ref{defsmoothingideal}). This is coherent with the smoothing operation presented in \cite{secondpaper}. Denote by $\opcf$ the group of totally positive units in $\ok$ which are congruent to $1 \mod \goth{f}$. This is a free abelian group of rank $1$ by Dirichlet's unit theorem. Therefore, there is a unique element $\eps \in \units$ such that $\opcf = \eps^{\zz}$ and $\sigma_{\rr}(\eps) < 1$.

Let us now consider a primitive vector $h$  in $L = \goth{f} \goth{b}^{-1}$ such that $h/q \equiv 1 \mod L$ and $h/N$ is a generator of the cyclic group $\goth{a}^{-1}L/L$. Such a vector is called an admissible base point. Let us fix a $\zz$-basis $\BL = [e_0, e_1, e_2]$ of $L$ such that $e_0 = h$ and $[e_0/N, e_1, e_2]$ is a $\zz$-basis of $\goth{a}^{-1}L$. We further assume that the basis $\BL$ is positive for the orientation obtained on $\kk$ by pulling back the standard orientation of $\rr \times \cc \simeq \kk \otimes \rr$. Let us denote by $a = a_h$ the integral linear form defined on $L$ by $a(y) = -\sign(\Im(\sigma_{\cc}(\eps)))\det_{\BL}(h, \eps h, y)$ and denote by $b$ the integral linear form on $L$ defined by $b(y) = a(\eps y)$. Bergeron, Charollois and Garc\'ia define an arithmetic evaluation of the elliptic Gamma function as follows:
$$ \Gamma_{\goth{f}, \goth{b}, \goth{a}}(\eps, h) = \frac{\Gamma_{a, b}(0)\left( \sigma_{\cc}\left(\frac{h}{q}\right), \sigma_{\cc}, L \right)^N}{\Gamma_{a, b}(0)\left( \sigma_{\cc}\left(\frac{h}{q}\right), \sigma_{\cc}, \goth{a}^{-1}L \right)}$$
where $\Gamma_{a, b}(v = 0)$ is the geometric elliptic Gamma function introduced by Felder, Henriques, Rossi and Zhu in \cite{FDuke}. These evaluations of the elliptic Gamma function are expected to yield algebraic numbers in the abelian extension $\classfieldf$ of $\kk$, and they should more precisely satisfy the following conjecture:

\begin{generalbis}{Conjecture}{[\hspace{1sp}\cite{BCG}, Conjecture]}\label{conjecturecubic}
Let $\kk$ be a complex cubic field and $\goth{f} \neq \ok$ be an integral ideal in $\ok$. Fix a complex embedding $\sigmaplus$ of $\classfieldf$ which extends the fixed complex embedding $\sigma_{\cc}$ of $\kk$. Fix a class $\goth{c}$ in $\classgroupf$ and an integral ideal $\goth{b}$ in the class $\goth{c}$. Fix a smoothing ideal $\goth{a}$ of norm $N$ such that $N$ is coprime to $6\cdot q\cdot\goth{b}$ where $q\zz = \goth{f} \cap \zz$. Then, the complex number $ \Gamma_{\goth{f}, \goth{b}, \goth{a}}(\eps, h)$ is independent of the choice of admissible base point $h$ as well as from the choice of ideal $\goth{b}$ and it is the image in $\cc$ of an algebraic unit $u_{\goth{f}, \goth{c}, \goth{a}}$ in $\classfieldf$ under the complex embedding $\sigmaplus$ of $\classfieldf$. Moreover:
\begin{itemize}
\item[(i)] Every embedding of $\classfieldf$ extending the real embedding of $\kk$ maps $u_{\goth{f}, \goth{c}, \goth{a}}$ to the unit circle.
\item[(ii)] If $\goth{c} \to \sigma_{\goth{c}}$ is the Artin map, the explicit reciprocity law is given by $\sigma_{\goth{c}'}(u_{\goth{f}, \goth{c}, \goth{a}}) = u_{\goth{f}, \goth{c}\goth{c}', \goth{a}}$ for any class $\goth{c}'$ in the narrow ray class group at $\goth{f}$. 
\end{itemize}
\end{generalbis}

We note that, if proven, this conjecture would constitute progress towards a solution to Hilbert's 12th problem for complex cubic fields, as some (and possibly all) finite abelian extensions of complex cubic fields would be described by evaluations of the elliptic Gamma function. Another remark we can make in view of the results of \cite{secondpaper} on the smoothing operation is that Conjecture \ref{conjecturecubic} should still hold for $N$ not coprime to $6$ if one is willing to replace the complex number $ \Gamma_{\goth{f}, \goth{b}, \goth{a}}(\eps, h)$ by its power $ \Gamma_{\goth{f}, \goth{b}, \goth{a}}(\eps, h)^{24}$ as $\mathcal{D}(3, 2) = 8$ and $\mathcal{D}(3, 3) = 3$. Very recently, the independence statement has been proven by Gray in \cite{TeymourGray} under mild hypotheses on $N$, that is, the complex number $ \Gamma_{\goth{f}, \goth{b}, \goth{a}}(\eps, h)$ is indeed independent of the choice of base point $h \in L$.

Let us briefly present an example supporting Conjecture \ref{conjecturecubic}. A more detailed presentation of this example is given in section \ref{detailledcubic}. Consider the field $\kk = \qq(z)$ where $z = e^{2i\pi/3}13^{1/3}$ is the root of the polynomial $x^3-13$ in the upper half-plane. Let $\goth{f}$ be the unique prime ideal in $\ok$ of norm $q = 3$, $\goth{b} = \ok$ and $\goth{a}$ the unique prime ideal of norm $N = 5$ in $\ok$. The fundamental unit $\eps = 2z^2 - 3z -4$ for $\opcf$ satisfies $\sigma_{\rr}(\eps) < 1$. Then, for the admissible base point $h = -(21z^2 + 42z + 114)$ we compute
$$\Gamma_{\goth{f}, \goth{b}, \goth{a}}(\eps, h) = \frac{\Gamma\left(\frac{1}{3}, -\frac{\eps^{-1} + 5348}{1965}, -\frac{\eps+467}{1965}\right)^{5}}{\Gamma\left(\frac{5}{3}, -\frac{\eps^{-1}+5348}{393}, -\frac{\eps+467}{393}\right)} \approx -0.0660917... + i\cdot0.0932299...$$
with 1000 digits of precision. This special value of the elliptic Gamma function coincides up to this precision with a root of the polynomial 
\begin{align*}
P_{\mathrm{abs}} = & \,x^{18} + 384x^{17} + 2310x^{16} - 10646490x^{15} + 1596241353x^{14} + 18608357181x^{13} \\ 
&+ 156933809421x^{12} + 215098256580x^{11} + 381407365338x^{10} + 338205493469x^9 \\
& +381407365338x^8 + 215098256580x^7 +  156933809421x^6 + 18608357181x^5 \\
&+ 1596241353x^4  - 10646490x^3 + 2310x^2 + 384x + 1
\end{align*}
which defines an absolute equation of $\classfieldf$ over $\qq$. Details on how this example was produced may be found in \cite{thesis}.

To support their conjecture, Bergeron, Charollois and Garc\'ia proved a key unconditionnal result, that is a version of Kronecker's limit formula for complex cubic fields:

\begin{generalbis}{Theorem}{[\hspace{1sp}\cite{BCG}, Theorem 3.2]}\label{KLFcubic}
The modulus of the complex number $ \Gamma_{\goth{f}, \goth{b}, \goth{a}}(\eps, h) $ is independent of the choice of admissible base point $h$ and it satisfies the Kronecker limit formula:
$$N\cdot \zeta'_{\goth{f}}(\goth{b}, 0) - \zeta'_{\goth{f}}(\goth{a}\goth{b}, 0)  = \log\left|\Gamma_{\goth{f}, \goth{b}, \goth{a}}(\eps, h)\right|^2.$$
\end{generalbis}

Conjecture \ref{conjecturecubic} together with Theorem \ref{KLFcubic} imply that these higher elliptic units $u_{\goth{f}, \goth{c}, \goth{a}}$ should be smoothed versions of Stark units (see section \ref{starksection} for a discussion of the relations with the rank one abelian Stark conjectures). In the cubic case, Conjecture \ref{conjecturecubic} has been tested numerically on thousands of examples and in the following sections we generalise this conjecture to ATR fields of degree $n \geq 4$ in some specific cases supported by numerical evidence.

\section{General construction of higher elliptic units}\label{sectiongcheu}

In this section we develop the general construction of higher elliptic units above ATR fields of degree $n \geq 3$ as products of special values of $G_{n-2}$ functions. We first recall some notations from \cite{secondpaper}. Then, we present the basic ingredients for the construction of these conjectural elliptic units and the precise general conjecture on the algebraic nature of such evaluations.

\subsection{Geometric families of $G_{r}$ functions and cocycle interpretations}

In \cite{firstpaper} we introduced a geometric framework allowing us to construct geometric families out of the multiple elliptic Gamma functions, generalising the construction in \cite{FDuke} to the whole hierarchy of $G_r$ functions. This allowed us to clearly rephrase the modular properties satisfied by these functions as multiplicative $(n-2)$-coboundary relations under an action of $\slnz{n}$. Let us fix some notations to state the main results of \cite{firstpaper} and \cite{secondpaper} which are relevant for the construction of higher elliptic units. Indeed, let us consider $V$ a $\qq$-vector space of finite dimension $n = r+ 2\geq 2$. We fix $L' \subset L$ two $\zz$-lattices of rank $n$ in $V$ such that $L/L' \simeq \zsz{N}$ for some integer $N \geq 2$. The integer $N$ will be the smoothing index in our evaluations. We also define $\Lambda =\homlz$ and $\Lambda' = \homlprimez$. In \cite{firstpaper} we defined a higher elliptic Gamma function attached to any set of $r+1$ linear forms in $\Lambda$. It is a function of the form:
$$G_{r, a_0, \dots, a_r} := \begin{cases} V/L \times \cc \times \homlc & \to \cc \\
(v, w, x) & \mapsto G_{r, a_0, \dots, a_r}(v)(w, x, L) \end{cases} $$
and it satisfies modular transformation properties under the action of $\slnz{n}$ involving an associated collection of Bernoulli rational functions. For the purpose of this paper we focus immediately on the $N$-smoothed versions defined in \cite{secondpaper} by:
\begin{equation}\label{defsmoothedgr}
G_{r, a_0, \dots, a_r}(v)(w, x, L, L') = \frac{G_{r, a_0, \dots, a_r}(v)(w, x, L')^N}{G_{r, a_0, \dots, a_r}(v)(w, x, L)}.
\end{equation}
The benefit of the smoothing operation is that it simplifies the form of the modular transformation property under the action of the congruence subgroup
\begin{equation}\label{definitionmodulargroup}
\congruencegrouphnn := \left\{ g \in \slnz{n} \setseparator[\Big] g \equiv \begin{pmatrix} \ast & \ast &\dots & \ast \\ 0 & \ast & \dots & \ast \\ \vdots & \vdots & \vdots & \vdots \\ 0 & \ast & \dots & \ast \end{pmatrix} \mod N \right\}.
\end{equation}
Indeed, for most configurations of $r+2$ linear forms $a_0, \dots, a_{r+1}$ which are primitive on both $L$ and $L'$, we have the modular property:
\begin{equation}\label{cocyclesmoothed}
\left(\prod_{j = 0}^{r+1} G_{r, a_0, \dots, \omitvar{a_j}, \dots, a_{r+1}}(v)(w, x, L, L')^{(-1)^j}\right)^{\denomclasscohom} = 1 \end{equation}
where the notation $\omitvar{a_j}$ means as usual that the element $a_j$ is omitted and, as before, $\denomclasscohom = \prod_{p \mid N} p^{\lfloor \frac{n}{p-1} \rfloor}$ (see [\hspace{1sp}\cite{secondpaper}, Theorem 1]). The equalities given by \refp{cocyclesmoothed} are interpreted as multiplicative $(n-2)$-cocycle relations for the congruence group $\congruencegrouphnn$.

For the purpose of the construction of higher elliptic units, we shall evaluate these multiplicative ``almost'' $(n-2)$-cocycles against fundamental cycles associated to groups of totally positive units in ATR fields of degree $n$.

\subsection{Evaluation of the $G_{r}$ functions against $r$-cycles}\label{sectionevaluation}

In this section we give the basic geometric setup to construct higher elliptic units using higher elliptic Gamma functions. Consider an ATR number field $\kk$ of degree $n \geq 3$. Fix a complex embedding $\sigma_{\cc}$ of $\kk$ and fix an ordering $\sigma_1, \dots, \sigma_{n-2}, \sigma_{n-1} = \sigma_{\cc}, \sigma_n = \overline{\sigma_{\cc}}$ on the embeddings of $\kk$. We fix a global orientation of $\kk$ as follows: for any basis $B = [e_1, \dots, e_n]$ of the $\qq$-vector space $\kk$, the number $i\cdot\det((\sigma_j(e_k))_{1 \leq j,k\leq n})$ is non zero and real, and we say that the basis $B$ is positive if this number is positive. This applies in particular to any $\zz$-basis of a fractional ideal in $\kk$. Consider an integral ideal $\goth{f} \neq \ok$ and set $q$ to be the positive integer satisfying $q\zz = \goth{f} \cap \zz$. Let us now fix an integral ideal $\goth{b}$ coprime to $\goth{f}$ and set $L = \goth{f}\goth{b}^{-1}$. The integral ideal $\goth{b}$ will represent a class in some class group of $\kk$ (see section \ref{sectionconjecture}). Let us now define precisely the notion of smoothing ideal:

\begin{general}{Definition}\label{defsmoothingideal}
An integral ideal $\goth{a}$ is a smoothing ideal for $\goth{f}$ (resp. for $\goth{f}$ and $\goth{b}$) if the quotient $\ok/\goth{a}$ is a cyclic abelian group of order $N$ such that $N$ is coprime to $\goth{f}$ (resp. $N$ is coprime to $\goth{f}$ and to $\goth{b}$). The integer $N$ is then called the smoothing index.
\end{general}

This definition is coherent with the smoothing operation described in \cite{secondpaper} for the pair of lattices $L \subset \goth{a}^{-1}L$. Let us then fix a smoothing ideal $\goth{a}$ for $\goth{f}$ and $\goth{b}$. We consider specific base points in $L$ associated to the data $\goth{f}, \goth{b}, \goth{a}$ which we will use to evaluate geometric $G_{n-2}$ functions.

\begin{general}{Definition}\label{defadmissible}
A base point $h \in L = \goth{f}\goth{b}^{-1}$ is admissible (for the data $\goth{f}, \goth{b}, \goth{a}$) if $h/q \equiv 1 \mod L$ and $h/N$ generates the cyclic group $\goth{a}^{-1}L/L \simeq \zsz{N}$. A base point $h \in L$ is simply called ``weakly admissible'' if $h/q \equiv k\cdot \eps \mod L$ for some $k \in \zszx{q}$ and $\eps \in \units$ and $h/N$ generates the cyclic group $\goth{a}^{-1}L/L \simeq \zsz{N}$
\end{general}

Most of our constructions can be carried out using weakly admissible base points, the strong admissibility being only required to give a normalisation of our evaluations at the end of the computations. In this section we fix a weakly admissible base point $h \in L$ and we set $m$ to be the maximal integer satisfying $h/m \in L$. Let us also fix a positive $\zz$-basis $\BL$ of the $\zz$-module $L$. In what follows, we present an evaluation of some multiplicative $(n-2)$-cocycle $\psi$ built from $N$-smoothed $G_{n-2}$ functions against the fundamental cycle associated to the unit group
\begin{equation}\label{defopcf}
\opcf = \{\eps \in \units, \eps \equiv 1 \mod \goth{f}, \sigma(\eps) > 0, \text{ for all real embedding } \sigma \text{ of } \kk\}.
\end{equation}
By Dirichlet's unit theorem, this group is a free $\zz$-module of rank $r = n-2$. If we fix a set of fundamental units $\eps_1, \dots, \eps_r$ for this group, we may form the fundamental $r$-cycle $\Upsilon_{\goth{f}}$ defined by:
$$\Upsilonf = \sum_{\rho \in \goth{S}_r} \signature(\rho) [\eps_{\rho(1)}|\dots|\eps_{\rho(r)}] $$
where
$$[\eps_1 | \dots | \eps_r] = \left(\eps_1, \dots, \prod_{j = 1}^k\eps_j, \dots, \prod_{j=1}^r\eps_j\right).$$
This is an analogue for ATR fields of Sczech's cycle for totally real number fields (see \cite{Sczech}). In particular, it follows from [\hspace{1sp}\cite{Sczech}, Lemma 5] that the class of $\Upsilonf$ in $H_{r}(\opcf, \zz)$ is independent of the choice of fundamental units for $\opcf$. However, our construction of higher elliptic units depends on this choice. As a matter of fact, our construction can be carried out for any group of totally positive units $U$ satisfying
$$\opcf \subset U \subset \opck $$
where $\opck$ is the group of totally positive units in $\ok$, as if $\Upsilon_U$ if the fundamental cycle associated to $U$ then we expect the equality of cohomological pairings:
\begin{equation}\label{justifU}
\langle \psi, \Upsilonf \rangle = \sum_{\eps \in U/\opcf} \langle \eps \cdot \psi, \Upsilon_U \rangle.
\end{equation}
For the purpose of this paper, we shall focus on the best practical choice of subgroup $U$ we use in the computations. Indeed, on the one hand, the smallest $U$ is (especially if $U = \opcf$), the worse the properties of its fundamental cycle are expected to be (see section \ref{sectiongeomsetup}). On the other hand, we try to minimize the number of terms to compute in \refp{justifU} if we can help it. Thus, we define the specific unit group:
\begin{equation}\label{defuzf}
\uzfplus := \opck \cap (\zz + \goth{f})
\end{equation}
It is clear that $\opcf \subset \uzfplus \subset \opck$. We shall fix a set $\eps_1, \dots, \eps_r$ of fundamental units for this free group $\uzfplus$ and consider the associated Sczech cycle $\Upsilon = \Upsilon(\eps_1, \dots, \eps_r)$ defined by:
\begin{equation}\label{defUpsilon}
\Upsilon = \sum_{\rho \in \goth{S}_r} \signature(\rho) [\eps_{\rho(1)}|\dots|\eps_{\rho(r)}].
\end{equation}

Let us remark that the group $\uzfplus$ acts on $L$ by multiplication, preserving its orientation. Thus, the unit group $\uzfplus$ may be identified with an algebraic torus of rank $r$ in $\slnz{n}$ by identifying any unit $\eps \in \uzfplus$ with the matrix of multiplication by $\eps$ in the basis $B_L$. The group $\uzfplus$ also naturally acts on the dual space $\Lambda = \homlz$ by inverse multiplication: for any $b \in \Lambda$, $\eps \cdot b =b( \eps^{-1} \cdot)$. In particular, for any $b \in \Lambda$ and any $\beta \in L$, we get $(\eps \cdot b)(\eps \cdot \beta) = b(\beta)$.
 
Let us now describe the basic arithmetic evaluation procedure for $G_r$ functions. Let us consider $u_1, \dots, u_r$ to be fundamental units for the group $\uzfplus$ corresponding to one of the components of $\Upsilon$, that is $u_j = \prod_{i = 1}^j \eps_{\rho(i)}$ for some permutation $\rho \in \goth{S}_{n-2}$. Set $u_0 = 1$. Let us assume that the family $u_0, u_1, \dots, u_r$ is a free family of the $\qq$-vector space $\kk \simeq \qq^n$ (this corresponds to [\hspace{1sp}\cite{thesis}, (H1)]). This implies that the linear form $f := y \mapsto \det_{\BL}(u_0h, u_1h, \dots, u_rh, y)$ is non-zero for any $h \in \kk^{\times}$. We denote by $a = a(u_1, \dots, u_r, h)$ the unique primitive element in $\Lambda$ and by $\lambda = \lambda(u_1, \dots, u_r, h)$ the unique positive integer such that $\lambda m^{r+1} a = f$. The basic ingredient for our construction of higher elliptic units consists of evaluations of smoothed $G_{r}$ functions of the form:
\begin{equation}\label{arithmeticgrshape}
G_{r, \goth{f}, \goth{b}, \goth{a}}^{\pm}(u_1, \dots, u_r ; h) = \frac{G_{r, \pm(a, au_1, \dots, au_r)}(0)\left(\sigma_{\cc}\left(\frac{h}{q}\right), \sigma_{\cc}, L\right)^N}{G_{r, \pm(a,  au_1, \dots, au_r)}(0)\left(\sigma_{\cc}\left(\frac{h}{q}\right), \sigma_{\cc}, \goth{a}^{-1}L\right)}
\end{equation}
where the $\pm$ depends on some orientation choice which we discuss later on. Here $au_j$ stands for the function $y \to a(u_jy)$ which is also $u_j^{-1} \cdot a$. In the notation given by \refp{defsmoothedgr}, this evaluation can be written in terms of a smoothed $G_{r}$ function as:
$$G_{r, \goth{f}, \goth{b}, \goth{a}}^{\pm}(u_1, \dots, u_r ; h) = G_{r, \pm(a, au_1, \dots, au_r)}(0)\left(\sigma_{\cc}\left(\frac{h}{q}\right), \sigma_{\cc}, \goth{a}^{-1}L, L\right).$$
This evaluation is well-defined (see [\hspace{1sp}\cite{thesis}, lemma III.7] for a detailed discussion on this matter). Let us note that if we replace $\sigma_{\cc}$ by the other complex embedding of $\kk$ (which is $\overline{\sigma_{\cc}}$) we get essentially the same information as:
$$G_{r, \goth{f}, \goth{b}, \goth{a}}^{\pm}(u_1, \dots, u_r ; h, \overline{\sigma_{\cc}}) = \overline{G_{r, \goth{f}, \goth{b}, \goth{a}}^{\pm}(u_1, \dots, u_r ; h, \sigma_{\cc})}^{(-1)^r}.$$
This justifies our focus on only one of the two complex embeddings of $\kk$. In the special cases $r = 0, \,r = 1$ (that is, the cases where $\kk$ is imaginary quadratic and complex cubic respectively) we will identify as before:
\begin{align*}
\theta_{\goth{f}, \goth{b}, \goth{a}}^{\pm}(h) & = G_{0,\goth{f}, \goth{b}, \goth{a}}^{\pm}(\emptyset; h) \\
\Gamma_{\goth{f}, \goth{b}, \goth{a}}^{\pm}(\eps; h) & = G_{1, \goth{f}, \goth{b}, \goth{a}}^{\pm}(\eps; h).
\end{align*}
The evaluation against the full Sczech cycle associated to $\uzfplus$ then takes the shape:
\begin{equation}\label{fullarithmeticgrshape}
\langle \psi, \Upsilon \rangle = \prod_{\rho \in \goth{S}_{r}} G_{r, \goth{f}, \goth{b}, \goth{a}}^{\mu_{\rho}}([\eps_{\rho(1)} | \dots | \eps_{\rho(r)}] ; h_\rho)^{\nu_{\rho}}
\end{equation}
for some choices of admissible base points $h_{\rho} \in L$ and some choices of orientation signs $\mu_{\rho}, \nu_{\rho} \in \{ \pm 1\}$ for each $\rho \in \goth{S}_{n-2}$. For this evaluation to make complete sense we suppose that each unit system $1, \eps_{\rho(1)}, \dots, \prod_{j = 1}^r \eps_{\rho(j)}$ is a family of $r+1$ linearly independent vectors in $\kk \simeq \qq^{r+2}$. When this is the case, we shall say that $\eps_1, \dots, \eps_r$ are non-degenerate fundamental units for the group $\uzfplus$. Throughout this paper we shall always assume that this condition is satisfied.

 To define the complete evaluation against the Sczech cycle associated to $\opcf$, one must only describe the action of $\uzfplus/\opcf$ on this pairing. To describe this action we introduce for simplicity the set
\begin{equation}\label{defzfoneplus}\zfoneplus := \{ (k, \eps) \in \zszx{q} \times \uzfplus/\opcf \setseparator k\cdot\eps \equiv 1 \mod \goth{f}\} \simeq \uzfplus/\opcf.
\end{equation}
Here it is understood that $(k, \eps) \in \zfoneplus$ always comes with $k \in \zz_{>0}$. We may then show that:
\begin{equation}\label{fullarithmeticgrshapeprod}
\prod_{\eps \in \uzfplus/\opcf} \langle \eps \cdot \psi, \Upsilon \rangle = \prod_{\rho \in \goth{S}_{r}} \prod_{(k, \eps) \in \zfoneplus} G_{r, \goth{f}, \goth{b}, \goth{a}}^{\mu_{\rho}}([\eps_{\rho(1)} | \dots | \eps_{\rho(r)}] ; k\cdot \eps \cdot h_\rho)^{\nu_{\rho}}.
\end{equation}
In particular, the term corresponding to $(k, \eps) \in \zfoneplus$ does not depend on the particular representative we take in $(\zz - q\zz)_{>0} \times \uzfplus$ (see [\hspace{1sp}\cite{thesis}, lemmas III.28 and III.29]). For future reference, we denote this product by $\spvgrcomplete$, that is:
\begin{equation}\label{definvariant}
\spvgrcomplete = \prod_{\rho \in \goth{S}_{r}} \prod_{(k, \eps) \in \zfoneplus} G_{r, \goth{f}, \goth{b}, \goth{a}}^{\mu_{\rho}}([\eps_{\rho(1)} | \dots | \eps_{\rho(r)}] ; k\cdot \eps \cdot h_\rho)^{\nu_{\rho}}
\end{equation}
This expression is expected to be a ``higher elliptic unit'' above the number field $\kk$ when the parameters are well-chosen. We may give another vague form of our general conjecture:

\begin{generalbis}{Conjecture}{[Vague form 2]}\label{conjecturevagueformdeux}
Let $\kk$ be an ATR field of degree $n = r+2 \geq 3$. Suppose that $\goth{f}, \goth{b}, \goth{a}$ are as above and set $N = \norm{\goth{a}}$. Let $\eps_1, \dots, \eps_r$ be a set of non-degenerate fundamental units for $\uzfplus$. Then, for each permutation $\rho \in \goth{S}_r$, there isa n admissible base point $h_\rho$ and there are orientation signs $\mu_{\rho}, \nu_{\rho} \in \{\pm 1\}$ such that the complex number
\begin{equation}\label{uvagueformdeux}u = \spvgrcomplete^{\denomclasscohom}
\end{equation}
satisfies:
\begin{enumerate}
\item An algebraicity statement: the complex number $u$ is an algebraic unit in $\classfieldf$.
\item A Kronecker limit formula of the form:
\begin{equation}\label{klfvagueformdeux}
N\zeta'_{\goth{f}}(\goth{b}, 0) - \zeta'_{\goth{f}}(\goth{a}\goth{b}, 0) = \frac{1}{\denomclasscohom}\log\left|u\right|^2.
\end{equation}
\end{enumerate}
\end{generalbis}

In the next section, we shall make this conjectural statement more precise, by giving some insight on how the base points $h_{\rho}$ should be chosen, especially for $n \geq 4$. Indeed, the numerical computations carried out to formulate our conjecture show that for $n \geq 4$ the base points $h_{\rho}$ must satisfy some compatibility condition for Conjecture \ref{conjecturevagueformdeux} to hold.

\subsection{Conjecture on higher elliptic units}\label{sectionconjecture}

In this section we state the most precise version of our conjecture on higher elliptic units. Let us first recall the compatibility condition we impose on the families of base points $(h_{\rho})_{\rho}$ we use in the evaluation \refp{definvariant}. 

Recall that $\kk$ is an ATR field of degree $n = r+2$. Let us fix $\eps_1, \dots, \eps_r$ a set of non-degenerate fundamental units for $\uzfplus$. Fix a permutation $\rho \in \goth{S}_{r}$ and set $u_{j, \rho} = \prod_{i = 1}^{j} \eps_{\rho(i)}$ for all $1 \leq j \leq r$. Fix an admissible base point $h_{\rho} \in L$ and set as before $m_{\rho}$ to be the maximal positive integer satisfying $h_{\rho}/m_{\rho} \in L$. Recall that the linear $a_{\rho}$ is the unique primitive linear form on $L$ satisfying 
$$\lambda_{\rho} m_{\rho}^{r+1} a_{\rho}(y) = \det_{\BL}(h_{\rho}, u_{1, \rho} h_{\rho}, \dots, u_{r, \rho}h_{\rho}, y)$$
for some positive integer $\lambda_{\rho}$. To this linear form we may associate a generalised different ideal $\goth{D}(a_{\rho}, L)$ (see section \ref{sectiondefdifferent} for a discussion on these) by:
$$\goth{D}(a_{\rho}, L)^{-1} := \{x \in \kk \setseparator \forall\, y \in L, a_{\rho}(xy) \in \zz\}.$$
The compatibility condition for the base points may then be stated in terms of these generalised different ideals as follows:

\begin{general}{Definition}\label{defcompatible}
The family of admissible base points $(h_{\rho})_{\rho}$ is compatible (for the data $\goth{f}, \goth{b}, \goth{a}, \eps_1, \dots, \eps_r$) if and only if the ideal $\goth{D}(a_{\rho}, L)$ is independent of $\rho$.
\end{general}

This notion of compatibility is more general than the strict version presented in \cite{thesis} as we shall show in section \ref{sectiongeomsetup}. The main benefit of this new version is that it doesn't require us to make additional hypotheses to formulate. Our goal is to show that such families always exist (or at least that they exist for any given ATR field for some choice of fundamental units). We shall show that this is the case under some simplifying hypotheses in section \ref{choiceh}.

A first observation we make regarding the compatibility condition is that it depends on the ideals $(h_{\rho}/m_{\rho})\ok$ rather than on the base points $(h_{\rho})_{\rho}$, as explained by the following lemma:

\begin{general}{Lemma}\label{lemmadepidealhrho}
Let $h, h' \in \goth{f}\goth{b}^{-1}$ be two admissible base points. Let $m$ and $m'$ be the maximal integers satisfying $h/m, h'/m' \in \goth{f}\goth{b}^{-1}$. Suppose that the ideals $(h/m)\ok$ and $(h'/m')\ok$ are equal. Consider the primitive linear forms $a$ and $a'$ given as before by:
\begin{align*}
\lambda m^{r+1} a(y)&  = \det_{\BL}(h, u_1 h, \dots, u_rh, y)\\
\lambda' (m')^{r+1} a'(y) &= \det_{\BL}(h', u_1 h', \dots, u_rh', y) 
\end{align*}
where $\lambda, \lambda'$ are positive integers and $u_1, \dots, u_r$ are any totally positive units for which these linear forms are non-zero. Then $\goth{D}(a, L) = \goth{D}(a', L)$.
\end{general}

\begin{proof}
Consider two such admissible base points $h, h' \in \goth{f}\goth{b}^{-1}$. Since $(h/m)\ok = (h'/m')\ok$, there is a unit $\xi \in \units$ such that $h/m = \xi h'/m'$. Then, we have:
$$\det_{\BL}(h', u_1h', \dots, u_rh', y) = \norm{\xi} \left(\frac{m}{m'}\right)^{r+1} \det_{\BL}(h, u_1 h, \dots, u_rh, \xi^{-1}y).$$
Hence $\lambda = \lambda'$ and $a'= \norm{\xi} \cdot a(\xi^{-1}\cdot)$ where $\norm{\xi} = \pm 1$. For any $x, y \in \kk$, the conditions $a(xy) \in \zz$ and $a(\xi^{-1}xy) \in \zz$ are equivalent, therefore, by definition of $\goth{D}(a, L)$ and $\goth{D}(a', L)$ we get $\goth{D}(a, L) = \goth{D}(a', L)$.
\end{proof}

As a consequence of this simple result, there is still some ambiguity as to what we should take $h_{\rho}$ to be. Indeed, let us define:
\begin{equation}\label{defzfone}
\zfone := \{ (k, \eps) \in \zszx{q} \times \uzf/\unitsmodf \setseparator k \cdot \eps \equiv 1 \mod \goth{f}\} \simeq \uzf /\unitsmodf
\end{equation}
where $\unitsmodf$ is the set of units congruent to $1 \mod \goth{f}$ and $\uzf = \units \cap (\zz + \goth{f})$. It follows from lemma \ref{lemmadepidealhrho} that if $(h_{\rho})_{\rho}$ is a compatible set of admissible base points, then so is $(k_{\rho} \cdot \eps_{\rho} \cdot h_{\rho})_{\rho}$ for any choice of elements $(k_{\rho}, \eps_{\rho}) \in \zfone$. In addition, our construction is also given up to the orientations signs $\mu_{\rho}, \nu_{\rho} \in \{\pm 1\}$. Let us say a few words on these ambiguities in the construction. In [\hspace{1sp}\cite{thesis}, Proposition III.30] we proved the following: if $k, k' \in \zz$ are coprime to $q$ and $\eps, \eps' \in \uzf$ are such that $k \equiv k' \mod q$ and $k \cdot \eps \equiv k'\cdot \eps' \mod \goth{f}$ then:
$$G_{r, \goth{f}, \goth{b}, \goth{a}}^{\mu_{\rho}\norm{\eps}}([\eps_{\rho(1)} | \dots | \eps_{\rho(r)}] ; k\cdot \eps \cdot h_\rho)^{\nu_{\rho}\norm{\eps}} = G_{r, \goth{f}, \goth{b}, \goth{a}}^{\mu_{\rho}\norm{\eps'}}([\eps_{\rho(1)} | \dots | \eps_{\rho(r)}] ; k\cdot \eps \cdot h_\rho)^{\nu_{\rho}\norm{\eps'}}.$$
Hence this term depends only on the class of $(k, \eps)$ in $\zfone$. We may then define:
\begin{multline}\label{fullorientation}
\spvgrcompletebis = \\ \prod_{\rho \in \goth{S}_{n-2}}\prod_{(k', \eps') \in \zfoneplus} G_{r, \goth{f}, \goth{b}, \goth{a}}^{\mu_{\rho}\norm{\eps}}([\eps_{\rho(1)} | \dots | \eps_{\rho(r)}] ; kk'\cdot \eps\eps' \cdot h_\rho)^{\nu_{\rho}\norm{\eps}} \end{multline}
which, for the fixed data $r, \goth{f}, \goth{b}, \goth{a}, \eps_1, \dots, \eps_r, (h_{\rho})_{\rho}$ only depends on the finite possible choices for $\mu_{\rho} \in \{\pm 1\}, \nu_{\rho} \in \{\pm 1\}, (k_{\rho}, \eps_{\rho}) \in \zfone/\zfoneplus$. The collection of signs $\mu_{\rho}, \nu_{\rho}$ and elements $(k_{\rho}, \eps_{\rho}) \in \zfone/\zfoneplus$ will henceforth be called an ``orientation datum''. We state our main conjecture up to this choice of orientation, the main idea being that there is at least one orientation datum for which the special value $\spvgrcompletebis$ is an algebraic unit. A useful remark we may use to reduce the number of possible orientation data to check is the following: the element $(q-1, -1)$ belongs to $\zfone$ and acts as:
$$ G_{r, \goth{f}, \goth{b}, \goth{a}}^{\mu_{\rho}}([\eps_{\rho(1)} | \dots | \eps_{\rho(r)}] ; (q-1)\cdot (-1) \cdot h_\rho)^{\nu_{\rho}} = G_{r, \goth{f}, \goth{b}, \goth{a}}^{\mu_{\rho}(-1)^{r+1}}([\eps_{\rho(1)} | \dots | \eps_{\rho(r)}] ; h_\rho)^{\nu_{\rho}}$$
(see [\hspace{1sp}\cite{thesis}, Proposition III.30]). In particular we may restrict in general our search for orientation data from $(k, \eps) \in \zfone/\zfoneplus$ to $(k, \eps) \in \zfone/\langle\zfoneplus, (q-1, -1)\rangle$. 

Let us make a last remark before stating our main conjecture. This conjecture is stated in a way that makes use of an exact sequence in class field theory to reduce the number of computations we need to perform. The exact sequence is given by:
\begin{equation}\label{eqdefclzfplus}
1 \to \opcf \to \uzfplus \to \zszx{q} \to \classgroupf \to \clzfplus \to 1
\end{equation}
where $\clzfplus$ is the quotient of the group of fractional ideals coprime to $\goth{f}$ by the group of principal fractional ideals generated by some totally positive element $\beta$ such that there is an integer $k \in \zszx{q}$ for which $k\cdot\beta \equiv 1 \mod \goth{f}$. General class field theory associates to this class group a class field $\klzfplus$. The description of our conjectural higher elliptic units comes with an explicit formulation for the reciprocity law, that is, for the action of $\mathrm{Gal}(\classfieldf/\kk)$ on these units. Using the exact sequence above, we may show that the action of the subgroup $\mathrm{Gal}(\classfieldf/\klzfplus)$ conjecturally takes a very simple form.

\begin{general}{Main Conjecture}\label{conjecture}
Suppose that $\kk$ is an ATR field of degree $n = r+2 \geq 3$. Fix an integral ideal $\goth{f} \neq \ok$. Set $q\zz = \goth{f} \cap \zz$ with $q >0$. Let $\sigma_{\cc}$ be one of the two complex embeddings of $\kk$ and fix a complex embedding $\sigmaplus$ of $\classfieldf$ which extends $\sigma_{\cc}$. Fix a smoothing ideal $\goth{a}$ of norm $N$ coprime to $q$. Fix a class $\goth{c}$ in the narrow ray class group at $\goth{f}$. Fix an integer $\conjko > 0$ coprime to $q$ and an integral ideal $\goth{b}$ coprime to $\goth{f}$ such that the integral ideal $\conjko\cdot\goth{b}$ belongs to the class $\goth{c}$. Fix a set of non-degenerate fundamental units $\eps_1, \dots, \eps_r$ for $\uzfplus$. Assume that there exists a compatible set of admissible base points $\bars{h}$ for the data $\goth{f}, \goth{b}, \goth{a}, \eps_1, \dots, \eps_r$. Then:
\begin{enumerate}
\item There is an orientation datum $(\mu_{\rho}, \nu_{\rho}, k_{\rho}, \eps_{\rho})_{\rho}$ depending on $\goth{f}, \conjko, \goth{b}, \goth{a}, \eps_1, \dots, \eps_r, (h_{\rho})_{\rho}$ for which the following Kronecker limit formula holds:
\begin{equation}\label{klfconjecture}
N\zeta'_{\goth{f}}(\goth{c}, 0) - \zeta'_{\goth{f}}(\goth{a}\goth{c}, 0) = \log\left|I_{r,\goth{f}, \conjko\cdot\goth{b}, \goth{a}}(\eps_1, \dots, \eps_r ; \bars{h}, \bars{\mu}, \bars{\nu}, (k_{\rho}, \eps_{\rho})_{\rho})\right|^2.
\end{equation}
Any such orientation datum shall be called ``correct'' (for the data $\goth{f}, \conjko, \goth{b}, \goth{a}, \eps_1, \dots, \eps_r, \bars{h}$).
\item For the fixed data $\goth{f}, \conjko, \goth{b}, \goth{a}, \eps_1, \dots, \eps_r, \bars{h}$, the complex number
$$\mathcal{V}_{\goth{f}, \conjko, \goth{b}, \goth{a}}(\eps_1, \dots, \eps_r, \bars{h}) = I_{r,\goth{f}, \conjko\cdot\goth{b}, \goth{a}}(\eps_1, \dots, \eps_r ; \bars{h}, \bars{\mu}, \bars{\nu}, (k_{\rho}, \eps_{\rho})_{\rho})^{\denomclasscohom} $$
is independent of the choice of correct orientation datum.
\item For the fixed data $\goth{f}, \conjko, \goth{b}, \goth{a}, \eps_1, \dots, \eps_r$, the complex number 
$$\mathcal{V}_{\goth{f}, \conjko, \goth{b}, \goth{a}}(\eps_1, \dots, \eps_r) = \mathcal{V}_{\goth{f}, \conjko, \goth{b}, \goth{a}}(\eps_1, \dots, \eps_r, \bars{h})$$
is independent of the choice of compatible set of admissible base points $\bars{h}$.
\item For the fixed data $\goth{f}, \conjko, \goth{b}, \goth{a}$, the complex number 
$$\mathcal{V}_{\goth{f}, \conjko, \goth{b}, \goth{a}}= \mathcal{V}_{\goth{f}, \conjko, \goth{b}, \goth{a}}(\eps_1, \dots, \eps_r)$$
is independent of the choice of non-degenerate fundamental units $\eps_1, \dots, \eps_r$ for $\uzfplus$ (for which there exists at least one compatible set of base points $\bars{h}$).
\item For the fixed data $\goth{f}, \goth{c}, \goth{a}$, the complex number 
$$\mathcal{V}_{\goth{f}, \goth{c}, \goth{a}}= \mathcal{V}_{\goth{f}, \conjko, \goth{b}, \goth{a}}$$
is independent of the choice of representation of the class $\goth{c}$ by the pair $\conjko, \goth{b}$ such that the integral ideal $\conjko \cdot \goth{b}$ belongs to the class $\goth{c}$.
\item The complex number $\mathcal{V}_{\goth{f}, \goth{c}, \goth{a}}$ is the image in $\cc$ of an algebraic unit $u_{\goth{f}, \goth{c}, \goth{a}}$ in the class field $\classfieldf$ under the complex embedding $\sigmaplus$ of $\classfieldf$.
\item Any embedding of $\classfieldf$ above a real embedding of $\kk$ maps $u_{\goth{f}, \goth{c}, \goth{a}}$ to the unit circle.
\item If $\goth{c}'$ is a class in the narrow ray class group at $\goth{f}$ then the explicit reciprocity law is given by
$$\sigma_{\goth{c}'}(u_{\goth{f}, \goth{c}, \goth{a}}) = u_{\goth{f}, \goth{c}\goth{c}', \goth{a}} $$
where $\goth{c} \to \sigma_{\goth{c}}$ is the Artin map.
\item For a fixed integral ideal $\goth{b}$ coprime to $\goth{f}$, the algebraic units $u_{\goth{f}, \conjko \cdot \goth{b}, \goth{a}}$ obtained when varying $\conjko \in \zszx{q}/\zfoneplus$ are conjugated in the field extension $\classfieldf/\klzfplus$.
\end{enumerate}
\end{general}

\noindent \textbf{Remarks : } \begin{enumerate} 
\item In the case $n = 2$, that is for imaginary quadratic fields, it is already known (see section \ref{reviewquadratic}) that the values 
$$\theta_{\goth{f}, \goth{b}, \goth{a}}^{\pm}(h)^{\dNn{N}{2}}$$
are $q$-units in abelian extensions of imaginary quadratic fields as a consequence of the theory of Complex Multiplication and that they satisfy the second limit formula due to Kronecker.
\item In the case $n = 3$, the Kronecker limit formula is already proven in \cite{BCG} and the invariance statements 2-5 hold for the absolute values of these complex numbers as a consequence. The invariance statements 2-4 have been proven very recently in \cite{TeymourGray} under mild hypotheses on $N$ (that is for the complex numbers themselves and not only for their absolute values). The statements regarding the algebraic nature of these complex numbers remain conjectural.
\item Correct orientations in the conjecture may be computed by checking which ones among the $(2\cdot 2 \cdot \cardinalshort{\zfone/\zfoneplus})^{r!}$ possible choices give the desired Kronecker limit formula \refp{klfconjecture}. For the most simple examples, the quotient $\zfone/\langle \zfoneplus, (q-1, -1)\rangle$ is trivial and we only need to compute the orientation signs $\mu_{\rho}, \nu_{\rho}$. In many cases, the orientations may be taken as:
$$\mu_{\rho} = \nu_{\rho} = \signature(\rho) \cdot \sign(\mathrm{Reg}(\eps_1, \dots, \eps_r)) $$
where $\sign(\mathrm{Reg}(\eps_1, \dots, \eps_r))$ is the sign of the regulator of $\eps_1, \dots, \eps_r$ computing with this order on fundamental units, and the fixed order on the embeddings of $\kk$. These orientations are inspired by orientation conditions described in the study of (signed) Shintani fundamental domains in \cite{Colmez}, \cite{Diazydiaz} and \cite{Espinoza}. Exploring possible connections between our geometric constructions and Espinoza's signed fundamental domains for ATR fields might give us an explicit description of correct orientation data in our conjecture.

\item The most interesting case of the conjecture occurs when the partial zeta functions $\zeta_{\goth{f}}(\goth{c}, s)$ vanish at $s = 0$ with order exactly one. It is well-known that these partial zeta functions share the same order of vanishing (see \cite{Neukirch} for instance) and that for an ATR field this order is always $\geq 1$. For this order of vanishing to be exactly one, the number field $\classfieldf$ must be totally complex. This is not sufficient however, as $\goth{f}$ need not be equal to the conductor $\goth{f}'$ of the field extension $\classfieldf/\kk$. In this case, the order of vanishing of the partial zeta functions at $s= 0$ may be greater than $1$ depending on the possible vanishing of some Euler factors associated to prime factors of the integral ideal $\goth{f}/\goth{f}'$. When the partial zeta functions all vanish with order at least $2$ at $s = 0$, the conjecture still makes sense, but the unit $u_{\goth{f}, \goth{c}, \goth{a}}$ is a root of unity which we don't expect to generate the extension $\classfieldf/\kk$.
\item As the value of $\denomclasscohom$ is identically $1$ if $N$ is supported at primes $p > n+1$ we make the difference in the computations between ``small smoothings'' (which are divisible by some small prime $p \leq n+1$) and ``big smoothings'' which are not. The interest of big smoothings is that the minimal polynomial of the value $\mathcal{V}_{\goth{f}, \goth{c}, \goth{a}}$ should have smaller coefficients and it should therefore be more easily recognizable. The tradeoff however is that the computation time is linear in the smoothing index and therefore computations with ``big smoothings'' may take more time.
\end{enumerate}

Let us now illustrate this precise version of the conjecture on the quartic example presented in the introduction. Recall that $\kk = \qq(z)$ where $z$ is the root of the polynomial $x^4 - x^3 - x^2 - 5x + 1$ in the upper half-plane, that $\goth{f}$ is the unique prime ideal of norm $q = 7$ in $\kk$, and that representatives for $\classgroupf$ are given by $\ok, 2\ok, 3\ok, \goth{P}_{23}, 2\goth{P}_{23}, 3\goth{P}_{23}$. Here we have used the fact that $\clzfplus \simeq \zsz{2}$ and 
$$1 \to \zsz{3} \to \classgroupf \to \clzfplus \to 1.$$
Recall that $\goth{a}$ is the unique prime ideal of norm $N = 11$ in $\ok$. Thus $\denomclasscohom = 1$. The unit group $\opck = \uzfplus$ admits the fundamental units:
$$\eps_1 = 3z^3 + 4z^2 + 6z - 1 ,~\eps_2 = -z^3 + z^2 + 3z.$$
The unit group $\opcf$ is then generated by $\eps_1$ and $\eps_2^2$. The groups $\zfone$ and $\zfoneplus$ are given by $\zfone \simeq \zfoneplus \simeq \{(1, 1), (6, \eps_2)\}$. Let us consider the base points $h_1 = -7z^3 + 49z^2 - 133z + 63$ and $h_2 = -105z^3 + 42z^2 + 84z + 483$ for the case $\goth{b} = \ok$. They are compatible as $\goth{D}(a_1, \goth{f}\ok^{-1}) = \goth{D}(a_2, \goth{f}\ok^{-1}) = \goth{f} \cdot \goth{a}$. For $\goth{b} = \goth{P}_{23}$ we fix $h_3 = -147z^3 + 105z^2 - 2100z - 2142$ and $h_4 = 882z^3 - 630z^2 - 3339z - 3087$. They are also compatible as $\goth{D}(a_3, \goth{f}\goth{P}_{23}^{-1}) = \goth{D}(a_4, \goth{f}\goth{P}_{23}^{-1}) = \goth{f} \cdot \goth{a}$. These base points were explicitly chosen as generators of specific ideals following the general construction described in section \ref{sectiongeomsetup}. A correct orientation datum is then given by $(\mu_1 = 1, \nu_1 = 1, k_1 = 1, \eps_1 = 1)$ and $(\mu_2 = -1, \nu_2 = -1, k_2 = 1, \eps_2 = 1)$.  The complex numbers computed in the introduction are:
\begin{align*}
\mathcal{V}_{\goth{f}, \conjko, \ok, \goth{a}} = ~& G_{2, \goth{f}, \ok, \goth{a}}^{+}(\eps_1, \eps_1\eps_2 ; \conjko \cdot h_1) \cdot  G_{2, \goth{f}, \ok, \goth{a}}^{-}(\eps_2, \eps_2\eps_1 ; \conjko \cdot h_2)^{-1} \\
& \times G_{2, \goth{f}, \ok, \goth{a}}^{+}(\eps_1, \eps_1\eps_2 ; 6\conjko \cdot \eps_2 \cdot h_1) \cdot  G_{2, \goth{f}, \ok, \goth{a}}^{-}(\eps_2, \eps_2\eps_1 ; 6\conjko \cdot \eps_2 \cdot h_2)^{-1}
\end{align*}
and 
\begin{align*}
\mathcal{V}_{\goth{f}, \conjko, \goth{P}_{23}, \goth{a}} = ~& G_{2, \goth{f}, \goth{P}_{23}, \goth{a}}^{+}(\eps_1, \eps_1\eps_2 ; \conjko \cdot h_3) \cdot  G_{2, \goth{f}, \goth{P}_{23}, \goth{a}}^{-}(\eps_2, \eps_2\eps_1 ; \conjko \cdot h_4)^{-1}\\
& \times G_{2, \goth{f}, \goth{P}_{23}, \goth{a}}^{+}(\eps_1, \eps_1\eps_2 ; 6\conjko \cdot \eps_2 \cdot h_3) \cdot  G_{2, \goth{f}, \goth{P}_{23}, \goth{a}}^{-}(\eps_2, \eps_2\eps_1 ; 6\conjko \cdot \eps_2 \cdot h_4)^{-1}
\end{align*}
for $\conjko = 1, 2, 3$. They coincide up to 1000 digits of precision with the roots of the polynomial 
\begin{align*}
&X^6 + 1 +(-66007\zexc^3 + 315205\zexc^2 - 404670\zexc + 66406) (X^5+X) \\ + & (-2980477827\zexc^3 - 882760955\zexc^2 + 19370690773\zexc - 3656696288) (X^4 + X^2) \\ +& (153504704659\zexc^3 - 136634711339\zexc^2 - 544162307740\zexc + 108155656031) X^3
\end{align*}
which defines a relative equation of $\classfieldf$ over $\kk$. One may also check the expected Kronecker limit formula \refp{klfconjecture} up to 1000 digits of precision as:
$$11\cdot \zeta'_{\goth{f}}(\conjko \cdot \goth{b}, 0) - \zeta'_{\goth{f}}(\conjko \cdot \goth{a}\goth{b}, 0) \approx \log\left|\mathcal{V}_{\goth{f}, \conjko, \goth{b}, \goth{a}}\right|^2$$
for $\conjko = 1, 2, 3$ and $\goth{b} = \ok, \goth{P}_{23}$.

\subsection{Discussion of the conjecture}\label{discussionconjecture}

In this section we discuss various aspects of the Main Conjecture \ref{conjecture}, including the optimal setting adpoted in \cite{cheu} and a discussion of the link with Hilbert's 12th problem and the rank one abelian Stark conjectures.

\subsubsection{Optimal setting}\label{sectionoptimal}

In this section we describe a special setting that gives the simplest examples of computations of higher elliptic units. This is the setting in which we carried out most of our computations as it simplifies greatly the research of ATR fields together with their class field moduli for which simple enough fundamental units and compatible sets of base points may be found. We shall say that the pair $(\kk, \goth{f})$ is an optimal setting (where $\kk$ is an ATR field and $\goth{f} \neq (1)$ is an integral ideal in $\ok$) if $\ok/\goth{f}$ is cyclic and $\opck = \uzfplus = \opcf$. This implies in particular that $\zfoneplus = \{(1, 1)\}$. Finding an optimal setting gives us the best practical chances for numerical computations. However, this optimal setting condition is very restrictive. Indeed, if $\eps_1, \dots, \eps_r$ is a set of fundamental units for $\opck$, then $\opck = \opcf$ if and only if $\goth{f}$ is a divisor of the integral ideal $J(\kk) = \sum_{j = 1}^r (\eps_j-1)\ok$ (see [\hspace{1sp}\cite{thesis}, Lemma III.39]). In particular, for any given ATR field $\kk$, there may be only finitely many optimal settings $(\kk, \goth{f})$. 

The main benefit of this optimal setting is that most of the geometric setup we describe in section \ref{sectiongeomsetup} becomes simpler than in the general case. This case has been presented as the simplest case of our general conjecture in \cite{cheu}. In the present paper, we shall present examples in non-optimal settings to showcase our conjecture in more generality.

\subsubsection{Hilbert's 12th problem for ATR fields}\label{hilbertsection}

Hilbert's 12th problem asks for the construction of the abelian extensions of general number fields using analytic functions. It is solved for the field of rational numbers using the exponential function and for imaginary quadratic fields using the elliptic units built from the $\theta$ function. If proven, our Main Conjecture \ref{conjecture} would constitute progress towards a solution to Hilbert's 12th problem for ATR fields. Indeed, when the left-hand side of \refp{klfconjecture} does not vanish, the conjectural elliptic unit $u_{\goth{f}, \goth{c}, \goth{a}}$ we compute by evaluating the multiple elliptic Gamma functions is expected to generate a non-trivial extension of $\kk$ and in many cases the whole class field $\classfieldf$. It would then suffice to prove statements on the exact extensions of $\kk$ generated by these units and possibly some of their roots to obtain a general construction of $\classfieldf$ using analytic functions. 

In the more recent years, progress has been made on a $p$-adic version of Hilbert's 12th problem where $p$-adic analytic functions are used to describe abelian extensions of number fields.  We refer to \cite{DarmonDasgupta} and \cite{DarmonPozziVonk} for constructions of elliptic units above real quadratic fields. More generally we refer to \cite{DasguptaKakde} where the $p$-adic version of Hilbert's 12th problem is proven for totally real number fields as a consequence of the proof of the Brumer-Stark conjecture. It is interesting to note that all these settings, including ours, rely on the same general process of evaluating some arithmetic cocycle against some arithmetic cycle to produce meaningful algebraic invariants.

During the writing of this paper, some progress has been made towards a proof of Hilbert's 12th problem (the archimedean case) for real quadratic fields by Radchenko and Wheeler in \cite{RW}, where they prove that the explicit formulas for the conjectural Stark units produces algebraic numbers which can be expressed in terms of quantum dilogarithms. If statements on the exact field extensions generated by these numbers are proven, one may expect a solution to Hilbert's 12th problem for real quadratic fields using these analytic functions.

\subsubsection{Rank one abelian Stark conjectures}\label{starksection}

We now discuss the link between our conjecture and the rank one abelian Stark conjecture for ATR fields (see \cite{Stark}, \cite{rankonestark}). This conjecture is expressed for narrow ray class fields as follows:

\begin{generalbis}{Conjecture}{[Rank one abelian Stark conjecture for ATR fields]}
Let $\kk$ be an ATR field and $\goth{f} \neq \ok$ be an integral ideal. Write $e_{\goth{f}}$ for the number of roots of unity in $\classfieldf$. Fix a complex embedding $\sigmaplus$ of $\classfieldf$ extending a fixed complex embedding $\sigma_{\cc}$ of $\kk$. There is a unit $u_{Stark}$ in $\classfieldf$ such that 
\begin{itemize}
\item For any class $\goth{c}$ in $\classgroupf$, if $\sigma_{\goth{c}} \in \mathrm{Gal}(\kk^{+}(\goth{f})/\kk)$ is the image of $\goth{c}$ under Artin's map then:
$$\zeta'_{\goth{f}}(\goth{c}, 0) = -\frac{1}{e_{\goth{f}}} \log|\sigmaplus\left(u_{Stark}^{\sigma_{\goth{c}}}\right)|^2.$$
\item Every complex embedding of $\,\kk^{+}(\goth{f})$ above a real embedding of $\,\kk$ maps $u_{Stark}$ to the unit circle.
\item The number field $\kk^{+}(\goth{f})(u_{Stark}^{1/e_{\goth{f}}})$ is an abelian extension of $\,\kk$.
\end{itemize}
\end{generalbis}

If we assume both our Main Conjecture \ref{conjecture} and the rank one abelian Stark conjecture for ATR fields, then we obtain for all classes $\goth{c}$ of the narrow ray class group at $\goth{f}$ the equality:
$$\frac{1}{\denomclasscohom}\log| \sigmaplus(u_{\goth{f}, \goth{c}, \goth{a}})|^2 = -\frac{1}{e_{\goth{f}}} \log|\sigmaplus\left(\left(u_{Stark}^{\sigma_{\goth{c}}}\right)^{(N-\sigma_{\goth{a}})}\right)|^2.$$
In particular, up to some roots of unity of well-controled order, we should get the equality in $\classfieldf$:
\begin{equation}\label{liensmoothstark}
u_{\goth{f}, \goth{c}, \goth{a}}^{-e_{\goth{f}}} = u_{Stark}^{\denomclasscohom\sigma_{\goth{c}}(N-\sigma_{\goth{a}})}.
 \end{equation}
Thus our higher elliptic units give a conjectural analytic formula for smoothed Stark units and not only for their modulus at places of $\classfieldf$ above the unique complex place of $\kk$. Other formulas for Stark units above ATR fields were proposed by Ren and Sczech for complex cubic fields in \cite{RenSczech} and more recently by Black (see \cite{Black}) for general ATR number fields using generalised log-gamma functions. Another approach to the construction of Stark units specifically above ATR fields $\kk$ of degree $n = 2n'$ containing a totally real subfield $\bb{F}$ of degree $n'$ using Hilbert modular forms was presented by Charollois and Darmon in \cite{CharolloisDarmon}. They give a conjectural formula for the Stark unit as the evaluation of some multiplicative $n'$-cocycle (their version of the ``Abel-Jacobi map'') against an $n'$-cycle in the group coholomogy of $\mathrm{SL}_2(\mathcal{O}_{\bb{F}})$. It is worth noting that the presence of a real subfield is often an obstacle in our construction as we need to define higher elliptic units with evaluations of $G_r$ functions at real algebraic irrational points. Nevertheless, it would be interesting future work to try to find a relation between any of these constructions of Stark units and ours.

\section{Existence of compatible base points and explicit constructions}\label{sectiongeomsetup}

In this section we take a deep dive in the construction of higher elliptic units above ATR fields and make the construction as explicit as possible (see section \ref{sectiondetailledgeomsetup}). Then, in section \ref{sectiondifferentideals}, we reformulate the compatibility condition for sets of base points and identify the main obstructions in the search of such compatible sets. These obstructions come from a collection of generalised different ideals $\tildeD_{\rho}$ attached to the components $\Upsilon(\rho)$ of the cycle $\Upsilon$. In section \ref{choiceh} we prove Theorem \ref{theoremchoiceh} by showing that under some simplifying hypothesis we may always find a set of compatible base points for which the complexity of the evaluation of higher elliptic units is well-controlled. To prove this statement we introduce a class of ideals called \textit{helper ideals} and we prove that compatible base points $(h_{\rho})_{\rho}$ may be taken such that for any $\rho \in \goth{S}_{n-2}$:
$$h_{\rho} \ok = \frac{qN}{\goth{a}\goth{b}} \cdot \tildeD_{\rho} \cdot (m_{\rho} \cdot \helper) $$
where $m_{\rho}$ is a positive integer coprime to $q$ and $\helper$ is a common helper ideal. We then determine the values of the parameters that govern the complexity of the evaluation of higher elliptic units for this specific choice of base points.

\subsection{Detailed geometric setup}\label{sectiondetailledgeomsetup}

In this section we give a very detailed description of the formalism introduced in section \ref{sectionevaluation} for the construction of higher elliptic units.

\subsubsection{The linear forms $a, au_1, \dots, au_r$ and their positive dual family}\label{sectionbasicgeomsetuph}

In this section we give a detailed geometric setup to make our conjectural elliptic units completely explicit. This also allows us to understand how the base points $h_{\rho}$, $\rho \in \goth{S}_r$ should be chosen for our evaluations. To this end we shall express explicitly the arithmetic $G_r$ function defined by \refp{arithmeticgrshape} and analyse the properties of the underlying geometric $G_r$ function using the explicit definitions given in \cite{firstpaper} and \cite{secondpaper}.
Let $\kk$ be an ATR field of degree $n = r+2$ and let $\goth{f}$, $\goth{b}$, $\goth{a}$ be given as in section \ref{sectionevaluation}. Recall that we have fixed an ordering on the embeddings of $\kk$ and a global orientation of $\kk$ (see section \ref{sectionevaluation}).
Put $L = \goth{f}\goth{b}^{-1}$. Fix a system of non-degenerate fundamental units $(\eps_1, \dots, \eps_r)$ of $\uzfplus$. Let us fix a permutation $\rho \in \goth{S}_r$ and write $(u_1, \dots u_r)= [ \eps_{\rho(1)} | \dots | \eps_{\rho(r)}]$. For convenience, we will always write $u_0 = 1$. Choose $h \in L$ a weakly admissible base point (see Definition \ref{defadmissible}). Let us denote by $m$ the maximal integer satisfying $h/m \in L$. We recall that the linear form $a = a_h$ is the unique primitive linear form on $L$ satisfying:
\begin{equation}\label{defa}
\lambda m^{r+1} a = \det_{\BL}(h, u_1h, \dots, u_rh, \cdot)
\end{equation}
for some positive integer $\lambda$. This linear form does not depend on the choice of positive basis of $L$, therefore we may assume that the chosen basis $\BL$ is given by $\BL = [e_0 = h/m, e_1, \dots, e_{r+1}]$ such that for any $1 \leq j \leq r$, $u_jh = \sum_{k = 0}^j mc_{jk0}e_k$ where the coefficients $c_{jk0}$ are integers and $c_{jj0} > 0$. This is done by computing the Hermite Normal Form of the matrix expressing $(u_0h, \dots, u_rh)$ in any positive basis of $L$ (see for instance \linebreak\cohenhnf). Explicitly, $(h, u_1h, \dots, u_rh) = m \cdot [e_0, \dots, e_{r+1}]\cdot \matrixgothu$ where $\matrixgothu$ is the $(r+2)\times(r+1)$ matrix:
\begin{equation}\label{defmatu} \matrixgothu =  \begin{pmatrix} 1 & c_{100} & \dots & c_{j00} & \dots & c_{r00}  \\ 0 & c_{110} & \dots & c_{j10} & \dots & c_{r10}  \\ \vdots & \vdots & \vdots & \vdots & \vdots & \vdots  \\ 0 & 0 & \dots & c_{jj0} & \dots & c_{rj0}  \\ \vdots & \vdots & \vdots & \vdots & \vdots & \vdots  \\ 0 & 0 & \dots & 0 & \dots & c_{rr0} \\0 & 0 & \dots & 0 & \dots & 0 \end{pmatrix}.\end{equation}
In this basis $\BL$ we identify the units $u_j$ with the matrix $(c_{jkl})_{k,l} \in \slnz{n}$ such that $u_je_l = \sum_{k = 0}^{r+1}c_{jkl}e_k$. Fix the dual basis $\BLambda = (f_0, \dots, f_{r+1})$ of $\Lambda = \homlz$ such that $f_j(e_k) = \delta_{jk}$ where $\delta_{jk}$ is the Kronecker symbol. In this basis the linear form $a = a_h$ defined in section \ref{sectionevaluation} is exactly $f_{r+1}$. For $1 \leq j \leq r$, the composition of $a$ with multiplication by $u_j$ is written $au_j = \sum_{l = 1}^{r+1} c_{j(r+1)l} f_l$. Since the fundamental units $\eps_1, \dots, \eps_r$ are non-degenerate, the family of linear forms $(a, au_1, \dots, au_r)$ is free in $\myhom{\qq}{\kk}{\qq} \simeq \qq^n$ therefore we may denote as usual by $\gamma$ the unique primitive vector in $L$ satisfying $\det_{\BLambda}(a, au_1, \dots, au_r, \cdot) = s\gamma$ for some positive integer $s$. Since for $0 \leq j \leq r$, $au_j(\gamma) = 0$ and $au_j(h) = 0$, we get 
$$\bigcap_{j = 0}^r \ker(au_j) = \qq \gamma = \qq h $$
and therefore there is a sign $\eta \in \{-1, +1\}$ such that $\eta m \gamma = h$. 

Let us now describe how we construct a positive dual family $\alpha_0, \dots, \alpha_r$ to $a, au_1, \dots, au_r$ in $L$ (see [\hspace{1sp}\cite{firstpaper}, Definition 5]). To achieve this, we first identify the family $\abar = (a = au_0, au_1, \dots, au_r)$ with the matrix obtained by concatenation of the coefficients of the linear forms $au_i$ in the basis $\BLambda$:
\begin{equation}\label{defatable}
\abar = \begin{pmatrix} au_0 \\ au_1 \\ \vdots \\ au_i \\ \vdots \\ au_r \end{pmatrix} = \begin{pmatrix} 0 & 0 & \dots & 0 & \dots & 1 \\ 0 & c_{1(r+1)1} & \dots & c_{1(r+1)k} &\dots& c_{1(r+1)(r+1)} \\ \vdots & \vdots & \vdots & \vdots & \vdots & \vdots \\ 0 & c_{i(r+1)1} & \dots & c_{i(r+1)k} & \dots & c_{i(r+1)(r+1)} \\\vdots & \vdots & \vdots & \vdots \\ 0 & c_{r(r+1)1} & \dots & c_{r(r+1)k} & \dots & c_{r(r+1)(r+1)} \end{pmatrix} = \begin{pmatrix} \textbf{0} & \mathcal{A} \end{pmatrix}
\end{equation}
This matrix has a first column filled with zeroes and the submatrix $\mathcal{A}$ is a square matrix of size $r+1$. Since $\det_{\BL}(a, au_1, \dots, au_r, \cdot)$ is a non-zero linear form on $\Lambda$, the matrix $\mathcal{A}$ is invertible. We shall now use the Smith normal form of $\mathcal{A}$ (see for instance \cohensnf). There are invertible matrices $U, V \in \glnz{r+1}$ and a unique diagonal matrix $S(\mathcal{A}) \in M_{r+1}(\zz)$ such that $U\mathcal{A}V = S(\mathcal{A})$ and
\begin{equation}\label{snfA} S(\mathcal{A}) = \begin{pmatrix} A_r & 0 & \dots & 0 \\ 0 & A_{r-1} & \dots & 0 \\ \vdots & \vdots & \vdots & \vdots \\ 0 & \dots & \dots & A_0  \end{pmatrix}\end{equation}
where $A_0, \dots, A_{r} \in \zz$ are the elementary divisors of $\mathcal{A}$ such that $A_i \mid A_{i+1}$ for all $0 \leq i \leq r-1$. Since $A_0$ is the gcd of all coefficients in $\mathcal{A}$ and $\mathcal{A}$ contains a $1$, it follows that $A_0 = 1$. From this writing we may identify
\begin{align}
\lambda(u_1, \dots, u_r; h) & = \lambda = \prod_{j = 1}^r c_{jj0} \label{identifylambda}\\
s(u_1, \dots, u_r; h) & = s = \prod_{j=1}^r A_j \label{identifys}
\end{align}
where $\lambda$ is defined in section \ref{sectionevaluation} and $\det_{\BLambda}(au_0, \dots, au_r, \cdot) = s\gamma$. We also define
\begin{equation}\label{defth}
t(u_1, \dots, u_r; h) = t = A_r.
\end{equation}
which will be an important parameter in our construction. 

Next, to obtain a positive dual family $\alphabar(u_1, \dots, u_r ; h) = \alphabar = (\alpha_0, \dots, \alpha_r)$ to $a, au_1, \dots, au_r$ in $L$ we set $\mathcal{B} = (\prod_{i = 1}^r A_i^{-1}) \mathrm{com}(\mathcal{A})^{T}$ where $\mathrm{com}(\mathcal{A})$ is the comatrix of $\mathcal{A}$. We claim that the columns of the integral matrix $\mathcal{B}$ give a positive dual family $\alpha_0, \dots, \alpha_r$ to $a, au_1, \dots, au_r$ in $L$. Indeed, if we write the coefficients of the matrix $\mathcal{B}$ as $\mathcal{B} = (b_{ij})_{1 \leq i, j \leq r+1}$ then the family of vectors $\alpha_0, \dots, \alpha_r$ defined by:
\begin{equation}\label{defalpha}
\alpha_j = \sum_{i = 1}^{r+1} b_{ij}e_i
\end{equation}
forms a positive dual family to $a, au_1, \dots, au_r$ in $L$ satisfying:
\begin{equation}\label{defalphaprop}
au_j(\alpha_j) = t, ~~~~ au_k(\alpha_j) = 0, ~\forall\, k \neq j.
\end{equation} 
We shall often write this as the multiplication of two matrices $\abar \cdot \alphabar = tI_{r+1}$. We say that this choice of $\alphabar$ is a uniform positive dual family to $\abar$ because the value $au_j(\alpha_j) = t$ is independent of $j$. We argue that this choice is minimal amongst all uniform positive dual families to $\abar$ in $L$, as explained by the following lemma.

\begin{general}{Lemma}\label{tmin}
Let $\mathcal{A}$ be the matrix defined in \refp{defatable}. 
\begin{itemize}
\item[(i)] Let $\mathcal{B}' \in M_{r+1}(\zz)$ be a square matrix of size $r+1$ such that $\mathcal{A} \cdot \mathcal{B}' = d I_{r+1}$ for some integer $d$. Then $t$ divides $d$.
\item[(ii)] Any uniform positive dual family $\alphabar'$ to $\abar$ in $L$ satisfies $\abar \cdot \alphabar' = d I_{r+1}$ where $t$ divide $d$.
\item[(iii)] If $\alphabar, \alphabar'$ are uniform positive dual families to $\abar$ in $L$ satisfying $\abar \cdot \alphabar = d I_{r+1}$ and $ \abar \cdot \alphabar' = d' I_{r+1}$ then for all $0 \leq j \leq r$, $d'\alpha_j - d\alpha'_j \in \pregcd(d, d')\zz\gamma$.
\end{itemize}
\end{general} 

\begin{proof}
(i) Consider the matrices $U,V \in \glnz{r+1}$ such that $U\mathcal{A}V = S(\mathcal{A})$ is the Smith normal form of $\mathcal{A}$. Then
$$ dI_{r+1} = U (dI_{r+1})U^{-1}= U \mathcal{A} \mathcal{B}' U^{-1} =\left(U \mathcal{A} V\right) \cdot \left(V^{-1} \mathcal{B}' U^{-1}\right) = \begin{pmatrix} t & 0 & \dots & 0 \\ 0 & A_{r-1} & \dots & 0 \\ \vdots & \vdots & \vdots & \vdots \\ 0 & \dots & \dots & 1  \end{pmatrix} V^{-1} \mathcal{B}' U^{-1}.$$
Since $V^{-1} \mathcal{B}' U^{-1} \in M_{r+1}(\zz)$, the first row of $S(\mathcal{A})\times V^{-1} \mathcal{B}' U^{-1}$ is divisible by $t$ and therefore $t$ divides $d$. \smallskip

(ii) Let $\alphabar' = (\alpha'_0, \dots, \alpha'_r)$ be a uniform positive dual family to $\abar$ in $L$ such that $\abar \cdot \alphabar' = dI_{r+1}$. Write the coordinates of the $\alpha'_j$ in the basis $B_L$ as $\alpha'_j = \sum_{i = 0}^{r+1} b'_{ij} e_i$. Then the matrix $\mathcal{B}' = (b_{ij}')_{1 \leq i,j \leq r+1}$ satisfies $\mathcal{A} \cdot \mathcal{B}' = d I_{r+1}$ and thus, by (i), $t$ divides $d$. \smallskip

(iii) Suppose that $\abar \cdot \alphabar = dI_{r+1}$ and $\abar \cdot \alphabar' = d' I_{r+1}$. Write as before the coordinates of $\alpha_0, \dots, \alpha_r$ and $\alpha'_0, \dots, \alpha'_r$ as:
$$\alpha_j = \sum_{i = 0}^{r+1} b_{ij} e_i ~~ \text{ and } ~~  \alpha'_j = \sum_{i = 0}^{r+1} b'_{ij} e_i.$$
Then the matrices $\mathcal{B} = (b_{ij})_{1 \leq i,j \leq r+1} \in M_{r+1}(\zz)$ and $\mathcal{B}' = (b'_{ij})_{1 \leq i,j \leq r+1} \in M_{r+1}(\zz)$ satisfy $\mathcal{A}\mathcal{B} =dI_{r+1}$ and $\mathcal{A} \mathcal{B}' = d' I_{r+1}$. Since $\mathcal{A}$ is invertible, it follows that $d'\mathcal{B} = d\mathcal{B}'$. Thus, for any $0 \leq j \leq r$: 
$$d'\alpha_j - d\alpha'_j = (d'b_{0j}-db'_{0j})e_0 =  (d'b_{0j}-db'_{0j})\eta \gamma \in \pregcd(d, d')\zz \gamma$$ 
as claimed.
\end{proof}

From now on we denote by $\alphabar = \alphabar(u_1, \dots, u_r; h)$ the minimal uniform positive dual family $\alphabar$ to $\abar$ given above. Recalling the construction in section \ref{sectionevaluation} and the definition of the geometric $G_{r, \abar}$ function (see [\hspace{1sp}\cite{firstpaper}, Proposition 7]) we may explicitly write the evaluation \refp{arithmeticgrshape} as:
\begin{equation}\label{grjustify} G_{r, \goth{f}, \goth{b}, \goth{a}}^{+}(u_1, \dots, u_r; h) =  \prod_{\delta \in F(\abar, \alphabar, 0)/\zz\gamma} \frac{G_r\left(\frac{h +q\delta}{q\gamma}, \frac{\alpha_0}{\gamma}, \dots, \frac{\alpha_r}{\gamma}\right)^N}{G_r\left(\frac{N (h+q\delta)}{q\gamma}, \frac{N\alpha_0}{\gamma}, \dots, \frac{N\alpha_r}{\gamma}\right)} \end{equation}
where 
$$F(\abar, \alphabar, 0) = \{ \delta \in L \setseparator \forall\, 0 \leq j \leq r, 0 \leq au_j(\delta) < au_j(\alpha_j) \}$$
and $F(\abar, \alphabar, 0)/\zz\gamma$ is a finite set. In \refp{grjustify} we have identified the elements of $\kk$ with their image in $\cc$ under the embedding $\sigma_{\cc}$ to lighten notations. Since $h = \eta m \gamma$ for some sign $\eta \in \{\pm 1\}$ we have a simplified form:
$$ G_{r, \goth{f}, \goth{b}, \goth{a}}^{+}(u_1, \dots, u_r; h) =  \prod_{\delta \in F/\zz\gamma} \frac{G_r\left(\frac{\eta m}{q} + \frac{\delta}{\gamma}, \frac{\alpha_0}{\gamma}, \dots, \frac{\alpha_r}{\gamma}\right)^N}{G_r\left(\frac{N\eta m}{q}+\frac{N\delta}{\gamma}, \frac{N\alpha_0}{\gamma}, \dots, \frac{N\alpha_r}{\gamma}\right)}$$
where we set $F = F(\abar, \alphabar, 0)$ to lighten notation. We may further improve this formulation to take into account the orientation datum $\mu, \nu, (k, \eps) \in \{-1, +1\} \times \{-1, +1\} \times \zfone/\zfoneplus$ (see section \ref{sectionconjecture}). Indeed, a reformulation of [\hspace{1sp}\cite{thesis}, Formula (III.39)] gives:
\begin{multline}\label{formeexpliciteune}\prod_{(k', \eps') \in \zfoneplus} G_{r, \goth{f}, \conjko \cdot \goth{b}, \goth{a}}^{\mu}(u_1, \dots, u_r; kk' \cdot \eps \eps' \cdot h)^{\nu} = \\ \prod_{\delta \in F/\zz\gamma} \prod_{(k', \eps') \in \zfoneplus} \left(\frac{G_r\left(\frac{\mu \cdot \norm{\eps} \cdot \conjko \cdot k \cdot \eta m}{q} + \frac{\delta}{\gamma}, \frac{\alpha_0}{\gamma}, \dots, \frac{\alpha_r}{\gamma}\right)^N}{G_r\left(\frac{N(\mu \cdot \norm{\eps} \cdot \conjko \cdot k \cdot \eta m)}{q}+\frac{N\delta}{\gamma}, \frac{N\alpha_0}{\gamma}, \dots, \frac{N\alpha_r}{\gamma}\right)}\right)^{\nu \mu^r \norm{\eps}^{r+1}}.
\end{multline}
This expression is shown to be well-defined in any case (see [\hspace{1sp}\cite{thesis}, Lemma III.7]). One must only be careful when treating the case where $\kk$ contains a proper subfield in which case an alternative formula for the $G_r$ functions may be needed to make proper sence of \refp{formeexpliciteune} (see for instance example \ref{quarticreal}). 

To end this section on the general explicit setup to compute \refp{arithmeticgrshape}, we say a few words on the size of the set $F/\zz\gamma$. Indeed, we may identify the set $F/\zz\gamma$ with the quotient space $L/M$ where $M = M(u_1, \dots, u_r; h) = \zz \gamma \oplus (\oplus_{j = 0}^r \zz \alpha_j)$. As this sublattice $M$ of $L$ has index $\det(\mathcal{B}) = t^{r+1}/s$ in $L$, it follows that the $N$-smoothed geometric function $G_{r, \abar}$ is a product of $t^{r+1}/s$ $N$-smoothed ordinary elliptic $G_r$ functions and so is the evaluation \refp{arithmeticgrshape}. In the best case scenario, $t^{r+1} = s$ and $G_{r, \abar}$ is a single $N$-smoothed ordinary elliptic $G_r$ function. However, since $A_0 = 1$ in \refp{snfA} this only happens if $t = s = 1$. When choosing the base point $h$ in the construction, we shall try to minimize the size of the set $F/\zz\gamma$, in order to get fast computations, and we shall always try to obtain this special case $t = s = 1$. This will be the subject of sections \ref{tildeforms} to \ref{sectiondeterminationlt}. The choice of base point $h$ shall then be fully explained in Theorem \ref{theoremchoiceh}.

\subsubsection{Dependence on $u_1, \dots, u_r$ and $h$}\label{tildeforms}

In this section we analyse the dependence of the quantities defined in section \ref{sectionbasicgeomsetuph} on the choice of units $u_1, \dots, u_r$ and on the choice of base point $h$. To this end we make the following important remark: the linear form $a = a(u_1, \dots, u_r ; h)$ may be easily expressed in terms of a linear form $\ta = \ta(u_1, \dots, u_r)$ which depends on $u_1, \dots, u_r$ but not on $h$. Indeed, let us fix $\Bok$ a positive $\zz$-basis of $\ok$ and let us consider the linear form:
\begin{equation}\label{eqdeftf}
\tildef := \left(y \to \det_{\Bok}(1, u_1, \dots, u_r, y)\right).
\end{equation}
If $\tlambda$ is the unique positive integer such that $\tildef$ maps $\ok$ to $\tlambda\,\zz$ then we may define the primitive linear form $\ta$ by:
\begin{equation}\label{eqdefta}
\ta := \left(y \to \tlambda^{-1}\det_{\Bok}(1, u_1, \dots, u_r, y)\right) = \tlambda^{-1} \tildef.
\end{equation}
Using this auxiliary function we may write:
\begin{equation}\label{eqlienata}
a = \left(y \to \ta\left(\frac{\tlambda}{\lambda}\frac{\norm{h/m}}{\norm{L}}\frac{m}{h}y \right) \right).
\end{equation}
Our claim is that most of the geometric setup introduced in section \ref{sectionbasicgeomsetuph} for the linear form $a = a(u_1, \dots, u_r, h)$ has a direct counterpart for the linear form $\ta = \ta(u_1, \dots, u_r)$ and that the two are closely related. In what follows, we will write with a $\sim$ the counterpart of these quantities. For the rest of this section, we assume that the basis $\Bok = [\tilde{e}_0 = 1, \tilde{e}_1, \dots, \tilde{e}_{r+1}]$ is such that $u_j = u_j \te_0 = \sum_{k = 0}^j \tc_{jk0}\te_k$ with $\tc_{jj0} > 0$ for $1 \leq j \leq r$. This is once again done by computing the Hermite Normal Form of a matrix representing $1, u_1, \dots, u_r$ in another positive $\zz$-basis of $\ok$. Explicitly, $(1, u_1, \dots, u_r) = [\tilde{e}_0, \dots, \tilde{e}_{r+1}]\cdot\tmatrixgothu$ where $\tmatrixgothu$ is the $(r+2)\times(r+1)$ matrix:
\begin{equation}\label{deftildematu} \tmatrixgothu =  \begin{pmatrix} 1 & \tilde{c}_{100} & \dots & \tilde{c}_{j00} & \dots & \tilde{c}_{r00} \\ 0 & \tilde{c}_{110} & \dots & \tilde{c}_{j10} & \dots & \tilde{c}_{r10} \\ \vdots & \vdots & \vdots & \vdots & \vdots & \vdots \\ 0 & 0 & \dots & \tilde{c}_{jj0} & \dots & \tilde{c}_{rj0} \\ \vdots & \vdots & \vdots & \vdots & \vdots & \vdots \\ 0 & 0 & \dots & 0 & \dots & \tilde{c}_{rr0}  \\0 & 0 & \dots & 0 & \dots & 0 \end{pmatrix}.\end{equation}
Let us fix $B^{\vee} = [\tilde{f}_0, \dots, \tilde{f}_{r+1}]$ the $\zz$-basis of the dual space $\myhom{\zz}{\ok}{\zz}$ such that $\tilde{f}_j(\tilde{e}_k) = \delta_{jk}$. For any $1 \leq j \leq r$, the unit $u_j$ may be identified with the matrix $(\tilde{c}_{jkl})_{k,l} \in \slnz{n}$ representing multiplication by $u_j$ in $\ok$, i.e. such that $u_j \tilde{e}_l = \sum_{k = 0}^{r+1} \tilde{c}_{jkl}\tilde{e}_k$ for any $0 \leq l \leq r+1$. In the basis $B^{\vee}$ of $\myhom{\zz}{\ok}{\zz}$ we may identify $\ta = \tilde{f}_{r+1}$ and $\ta u_j = \sum_{l = 1}^{r+1} \tilde{c}_{j(r+1)l}\tilde{f}_l$. The assumption that $\eps_1, \dots, \eps_r$ are non-degenerate guarantees that the family of linear form $\ta, \ta u_1, \dots, \ta u_r$ is free in $\myhom{\qq}{\kk}{\qq}$ and it is clear that $\cap_{j = 0}^r \ker(\ta u_j) = \qq$ where as before we write $u_0 = 1$. Let us then denote by $\ts$ the unique positive integer and $\tilde{\eta}$ the sign in $\{-1, +1\}$ such that $\det_{B^{\vee}}(\ta, \ta u_1, \dots, \ta u_r, \cdot ) = \tilde{\eta}\ts$. 

We may now describe a positive dual family $\talphabar = (\talpha_0, \dots, \talpha_r)$ to the family $\tabar = (\ta, \ta u_1, \dots, \ta u_r)$ in $\ok$ as follows. We associate to the family $\tabar$ a square matrix $\tilde{\mathcal{A}}$ of size $r+1$ as above for $\abar$ such that:
\begin{equation}\label{deftildeatable}
\tabar = \begin{pmatrix} \ta u_0 \\ \ta u_1 \\ \vdots \\ \ta u_i \\ \vdots \\ \ta u_r \end{pmatrix} = \begin{pmatrix} 0 & 0 & \dots & 0 & \dots & 1 \\ 0 & \tilde{c}_{1(r+1)1} & \dots & \tilde{c}_{1(r+1)k} &\dots& \tilde{c}_{1(r+1)(r+1)} \\ \vdots & \vdots & \vdots & \vdots & \vdots & \vdots \\ 0 & \tilde{c}_{i(r+1)1} & \dots & \tilde{c}_{i(r+1)k} & \dots & \tilde{c}_{i(r+1)(r+1)} \\\vdots & \vdots & \vdots & \vdots \\ 0 & \tilde{c}_{r(r+1)1} & \dots & \tilde{c}_{r(r+1)k} & \dots & \tilde{c}_{r(r+1)(r+1)} \end{pmatrix} = \begin{pmatrix} \textbf{0} & \tilde{\mathcal{A}} \end{pmatrix}
\end{equation}
Let us write as above the elementary divisors of this matrix $\tilde{\mathcal{A}}$ associated to the linear forms $\ta, \ta u_1, \dots, \ta u_r$ as $[\tA_r, \dots, \tA_1, \tA_0 = 1]$ such that $\tA_i | \tA_{i+1}$. Then we get as before:
\begin{align}
\tlambda(u_1, \dots, u_r) & = \tlambda = \prod_{j = 1}^r \tilde{c}_{jj0}, \label{identifytlambda}\\
\ts(u_1, \dots, u_r) & = \ts = \prod_{j=1}^r \tA_j \label{identifyts}
\end{align}
and if we set 
\begin{equation}\label{defttt}
\ttt = \tA_r 
\end{equation}
we may construct as above a uniform positive dual family $\talphabar$ to $\tabar$ in $\ok$ satisfying $\tabar \cdot \talphabar = \ttt I_{r+1}$ by computing the square matrix $\tilde{\mathcal{B}} = (\prod_{i = 1}^r \tA_i^{-1})\mathrm{com}(\tilde{\mathcal{A}})^{T}$ and extracting its columns. Explicitly, if $\tmatrixB = (\tilde{b}_{ij})_{0 \leq i,j \leq r+1}$ we set
\begin{equation}\label{deftalpha}
\talpha_j = \sum_{i = 1}^{r+1} \tilde{b}_{ij} e_i
\end{equation}
for $0 \leq j \leq r$, so that the elements $\talpha_0, \dots, \talpha_r$ form a positive dual family to $\ta, \ta u_1, \dots, \ta u_r$ in $\ok$ satisfying for all $0 \leq j \leq r$:
\begin{equation}\label{deftalphaprop}
\ta u_j(\talpha_j) = \ttt, ~~~~ \ta u_k(\talpha_j) = 0, ~\forall\, k \neq j.
\end{equation}

In the case $r = 1$ (that is $\kk$ is a complex cubic field and $u_1 = \eps$ is the fundamental unit of $\kk$), $\tlambda$ is the maximal integer satisfying $\eps \in \zz + \tlambda \ok$ and the integer $\ttt$ is defined by $\det(1, \eps, \eps^2) = \tlambda^3\ttt$. If we consider the set $E$ of elements $\beta \in \ok$ for which the $\zz$-module $\mathcal{O}(\beta) = \zz + \zz \eps + \zz \beta$ is a ring, and if for each $\beta \in E$ we denote by $\mathrm{ct}(\beta)$ the content of the ring $\mathcal{O}(\beta)$ in the sense of [\hspace{1sp}\cite{Bhargava}, Definition 14], then $\tlambda = \max\{ \mathrm{ct}(\beta) \setseparator \beta \in E\}$ as $\mathcal{O}(\tlambda \te_2) = \zz + \tlambda \ok$. On the other hand, when $\tlambda = 1$, the integer $\ttt$ is the index of the ring $\zz[\eps]$ in $\ok$. Therefore in general we call $\tlambda = \tlambda(u_1, \dots, u_r)$ the \textit{content} of the unit system $u_1, \dots, u_r$ and we call $\ttt = \ttt(u_1, \dots, u_r)$ the \textit{overflow} of this unit system. We insist that these quantities only depend on the unit system $(u_1, \dots, u_r)$. In what follows, we shall relate the quantities $\lambda, \tlambda, t, \ttt, s, \ts$ and more generally the two construction from sections \ref{sectionbasicgeomsetuph} and \ref{tildeforms}.

\subsubsection{Linking the two constructions}\label{sectionlinktilde}

In this section we give a first fundamental link between the two constructions associated to $a = a_h$ and to $\ta$. In the next sections, this link will be explored in more detail using generalised different ideals. Let us start by defining a specific integer $\ell$ which will correspond in some sense to the level of the computation associated to a base point $h$.

\begin{general}{Definition}\label{definitionlevel}
The level $\ell = \ell(u_1, \dots, u_r ; h)$ associated to the unit system $u_1, \dots, u_r$ and to the weakly admissible base point $h$ is defined as the maximal integer satisfying
$$\talpha_0h, \dots, \talpha_rh \in m\ell L + \zz h.$$
\end{general}

We immediately make the following remark: if in the basis $\BL$ the coordiantes of $\talpha_jh$ are given by $\talpha_jh = \sum_{k = 0}^{r+1} \langle \talpha_j h, e_k \rangle e_k$ then $\ell$ is the $\pregcd$ of the coefficients $\langle\talpha_j h, e_k \rangle$ for $0 \leq j \leq r$ and $1 \leq k \leq r+1$. Thus, in particular, if $d$ is any integer satisfying 
$$\talpha_0h, \dots, \talpha_rh \in m.d. L + \zz h$$
then $d$ must divide all the coefficients $\langle\talpha_j h, e_k \rangle$ for $0 \leq j \leq r$ and $1 \leq k \leq r+1$ and $d \mid \ell$.

The level $\ell$ naturally appears when comparing the families $\alpha_0, \dots, \alpha_r$ and $\talpha_0 h, \dots, \talpha_r h$. We shall now make this comparison and deduce a fundamental relation between the quantities $\lambda, \tlambda, t, \ttt$ and $\ell$.

\begin{general}{Proposition}\label{fundam}
Let $\epsilon = \sign(\norm{h}) \in \{-1,+1\}$.
\begin{itemize}
\item[(i)] The family $\epsilon\talphabar h = (\epsilon\talpha_0 h, \dots, \epsilon\talpha_r h)$ is a uniform positive dual family to $\abar$ in $L$ such that $\abar \cdot (\epsilon\talphabar h) = m \ell t I_{r+1}$. 
\item[(ii)] The following fundamental relation holds in $\zz$:
\begin{equation}\label{formulafundam}
\lambda \ell t = \frac{|\norm{h/m}|}{\norm{L}} \tlambda \ttt.
\end{equation}
\item[(iii)] For any $0 \leq j \leq r$ there is some integer $m_j \in \zz$ such that
\begin{equation}\label{eqlienalpha}
\alpha_j = \frac{\talpha_j h + m_j h}{\epsilon m \ell}.
\end{equation}
\end{itemize}
\end{general}

\begin{proof}
We shall prove simultaneously the three points of this proposition. We use formula \refp{eqlienata} which relates the linear forms $a$ and $\ta$ in conjunction with \refp{defa} and \refp{eqdeftf} to evaluate $\lambda m^{r+1} a (u_k \talpha_j h)$ for any $0 \leq j, k \leq r$. Indeed:
$$\lambda m^{r+1} a (u_k \talpha_j h) = \det_{\BL}(h, u_1 h, \dots, u_r h, u_k\talpha_j h).$$
Using standard linear algebra and the definition of the norm of the fractional ideal $L$ we may rewrite this equality using $\det_{\Bok}$  instead of $\det_{\BL}$ as
\begin{align*}
\lambda m^{r+1} a (u_k \talpha_j h) & = \det_{\BL}(\Bok) \det_{\Bok}(h, u_1h, \dots, u_rh, u_k \talpha_jh)\\
\lambda m^{r+1} a (u_k \talpha_j h) &=\frac{1}{\norm{L}}\det_{\Bok}(h, u_1h, \dots, u_rh, u_k \talpha_jh).
\end{align*}
Then, using the definition of the norm $\norm{h}$ of $h$ we get
\begin{align*}
\lambda m^{r+1} a (u_k \talpha_j h) & = \frac{\norm{h}}{\norm{L}} \det_{\Bok}(1, u_1, \dots, u_k \talpha_j)\\
\lambda m^{r+1} a (u_k \talpha_j h) & = \frac{\norm{h}}{\norm{L}} \tlambda \ta (u_k \talpha_j).
\end{align*}
We now use formula \refp{deftalpha} to conclude that for $k \neq j$, $a (u_k \talpha_j h) = 0$ and for $k = j$: 
$$\lambda m^{r+1} a (u_j \talpha_j h) = \frac{\norm{h}}{\norm{L}} \tlambda \ttt$$
which we may rewrite as:
\begin{equation}\label{tahpdf}
a u_j \talpha_j h = \frac{\norm{h}}{\norm{L}} \frac{\tlambda}{\lambda m^{r+1}} \ttt \neq 0.
\end{equation}
The lattice $L$ is a fractional ideal and $\talpha_0, \dots, \talpha_r \in \ok$, therefore $\talpha_0h, \dots, \talpha_r h \in L$ and since $a \in \homlz$, the rational number 
$$ d = \frac{\norm{h}}{\norm{L}} \frac{\tlambda}{\lambda m^{r+1}} \ttt$$
is actually a non-zero integer whose sign is $\epsilon = \epsilon_h = \sign(\norm{h})$. We have proven that $\epsilon \talphabar h$ is a uniform positive dual family to $\abar$ in $L$ such that 
$$\abar \cdot (\epsilon \talphabar h) = \epsilon d I_{r+1} = |d| I_{r+1}.$$
Let us now rewrite this using the level $\ell$ (see Definition \ref{definitionlevel}). By definition, there are integers $m_0, \dots, m_r$ and elements $\beta_0, \dots, \beta_r \in L$ such that for any $0 \leq j \leq r$:
$$\talpha_jh = \epsilon (\ell m \beta_j - m_j h).$$
Then \refp{tahpdf} may be written in terms of $\beta_j$ as:
$$a(u_k \beta_j) = \begin{cases} 0 & \text{ if } j \neq k \\ \frac{|d|}{m\ell} = \frac{|\norm{h/m}|}{\norm{L}} \frac{\tlambda}{\lambda} \frac{\ttt}{\ell} & \text{ if } j = k.\end{cases} $$
The elements $u_k\beta_j$ belong to $L$ therefore $a(u_k\beta_j) \in \zz$ for any $0 \leq j,k \leq r$ and $\beta_0, \dots, \beta_r$ is a uniform positive dual family to $a, au_1, \dots, au_r$ in $L$. It follows from lemma \ref{tmin} that $t$ divides the integer $d/(m\ell)$. We now prove that $t m \ell = |d|$. Indeed, it follows from lemma \ref{tmin}, (iii) that for any $0 \leq j \leq r$:
$$\frac{|d|}{m\ell} \alpha_j - t \beta_j \in t\zz \gamma.$$ 
Replacing $\beta_j$ by its definition we obtain:
$$\frac{d}{m\ell t} \alpha_j -   \frac{(\talpha_j +m_j)h}{m \ell} \in \zz \gamma$$
and thus:
$$ \talpha_jh \in \frac{d}{t} \alpha_j -m_jh +m\ell \zz\gamma$$ 
This shows in particular that 
$$\talpha_0h, \dots, \talpha_rh \in \frac{|d|}{t} L + \zz \gamma $$
and therefore the inequality $|d|/t \leq m \ell$ holds. This may be written as $|d| \leq m\ell t$ and the fact that $m \ell t$ divides $|d|$ with $|d| > 0$ implies that $|d| \geq m \ell t$. Thus the two quantities must be equal. We have thus proven that $\abar \cdot (\epsilon \talphabar h) = |d| I_{r+1} = m \ell t I_{r+1}$ and that the fundamental equality:
$$ \lambda \ell t = \frac{|\norm{h/m}|}{\norm{L}} \tlambda \ttt$$
holds in $\zz$. The identification of $\beta_j$ with $\alpha_j$ gives the desired relation between the families $\alphabar$, $\talphabar$:
$$\alpha_j = \frac{\talpha_j h + m_j h}{\epsilon m \ell}. $$
\end{proof}

This Proposition is fundamental in two different ways for the construction of higher elliptic units. On the one hand, the relation \refp{formulafundam} between the integers $\lambda, \tlambda, \ell, t, \ttt$ is crucial to understand how the choice of base point $h$ should be made to obtain interesting evaluations of the higher elliptic Gamma functions (and evaluations that can be computed in a small amount of time). On the other hand, the relation between the families $\alphabar$ and $\talphabar$ show that the base point $h$ impacts the parameters $m\alpha_0/h, \dots, m\alpha_r/h$ of the $G_r$ functions essentially through the level $\ell$. For instance, in the cubic case, the functions $\Gamma^{\pm}_{\goth{f}, \goth{b}, \goth{a}}$ are essentially evaluated at points
$$\tau = \pm\frac{\eps + n_0}{\ell\tlambda} = \frac{\tau_0}{\ell}, \sigma = \pm \frac{\eps^{-1} + n_1}{\ell\tlambda} = \frac{\sigma_0}{\ell}$$ for some integers $n_0, n_1$, where $\eps$ is a fundamental unit for $\opcf$ (or a fundamental unit for $\uzfplus$ following the remark made in section \ref{sectionevaluation}). This also means that when comparing two different choices of base points $h$ and $h'$, we only need to compare similar functions for different levels $\ell$ and $\ell'$.

\subsection{Generalised different ideals}\label{sectiondifferentideals}

\subsubsection{Definition of generalised different ideals}\label{sectiondefdifferent}

In this section we define generalised different ideals which we shall use to give some insight on the choice of base point $h$. 

\begin{general}{Definition}\label{definitiondifferent}
Let $f : \kk \to \qq$ be a non-zero $\qq$-linear form on $\kk$. Let $I$ be a fractional ideal of $\kk$. The fractional ideal $\goth{D}(f, I)$ defined by:
$$\goth{D}(f, I)^{-1} = \{ x \in \kk \setseparator \forall\, y \in I, f(xy) \in \zz\}$$
is called the different ideal of $f$ on $I$ by analogy with the usual different ideal $\goth{d}$ of $\kk$ given by $\goth{d} = \goth{D}(\trace, \ok)$.
\end{general}

The generalised different ideals satisfy basic properties:

\begin{general}{Lemma}\label{lemmabasicdifferent}
\begin{itemize}
\item[(i)] If $I' = \beta I$ for some $\beta \in \kk^{\times}$ then $\goth{D}(f, I') = \beta \cdot \goth{D}(f, I)$.
\item[(ii)] If $I = I_1 + I_2$ then $\goth{D}(f, I) = \goth{D}(f, I_1) + \goth{D}(f, I_2)$.
\item[(iii)] For any fractional ideal $I$ of $\kk$, $\goth{D}(f, I) = I \cdot \goth{D}(f, \ok)$.
\item[(iv)] If $f' = (y \to f(\beta y))$ for some $\beta \in \kk^{\times}$ then $\goth{D}(f', I) = \beta \cdot \goth{D}(f, I)$.
\item[(v)] For any $f \in \myhom{\qq}{\kk}{\qq} - \{0\}$, there is a unique element $\xi \in \kk^{\times}$ such that $f = (y \to \trace(\xi y))$. In particular:
$$\goth{D}(f, I) = (\xi) \cdot I \cdot \goth{d}.$$
\item[(vi)] For a fixed fractional ideal $I$, all different ideals $\goth{D}(f, I)$ belong to the class of $I \times \goth{d}$ in the class group of $\kk$.
\item[(vii)] If $f(I) \subset \zz$ then $\goth{D}(f, I) \subset \ok$. In addition, if $f(I) = \zz$ then $\goth{D}(f, I)$ is a primitive integral ideal.
\end{itemize}
\end{general}

\begin{proof}
See [\hspace{1sp}\cite{thesis}, Lemma III.12]. The proofs are straightforward and require only elementary facts on fractional ideals.
\end{proof}

In what follows, we shall be interested in the different ideals $\goth{D}(a, L)$ and $\goth{D}(\ta, \ok)$ associated to the linear forms $a$ and $\ta$. In particular, we may show that a simple relation exists between these two fractional ideals. This relation is closely linked to \refp{formulafundam}.

\begin{general}{Lemma}\label{lemmalinktarget}
The different ideals $\goth{D}(a, L)$ and $\goth{D}(\ta, \ok)$ are primitive integral ideals related by: 
$$\lambda \cdot\goth{D}(a, L) = \tlambda\cdot \frac{\norm{h/m}}{\norm{L}}\cdot\frac{m}{h}L\cdot\goth{D}(\ta, \ok).$$ 
\end{general}

\begin{proof}
First, since $a(L) = \zz$ and $\ta(\ok) = \zz$, it follows from lemma \ref{lemmabasicdifferent}, (vii) that these different ideals are primitive integral. Secondly, recall from \refp{eqlienata} that 
$$a = \left(y \to \ta\left(\frac{\tlambda}{\lambda}\cdot\frac{\norm{h/m}}{\norm{L}}\cdot\frac{m}{h}\cdot y \right) \right).$$
Thus, it follows from lemma \ref{lemmabasicdifferent}, (iv) that
$$\goth{D}(a, L) = \frac{\tlambda}{\lambda}\cdot\frac{\norm{h/m}}{\norm{L}}\cdot\frac{m}{h}\cdot \goth{D}(\ta, L) $$
and using lemma \ref{lemmabasicdifferent}, (iii) we get the desired relation:
$$\lambda \cdot\goth{D}(a, L) = \tlambda\cdot \frac{\norm{h/m}}{\norm{L}}\cdot\frac{m}{h}L\cdot\goth{D}(\ta, \ok).$$ 
\end{proof}

It follows from this lemma that the two main parameters in the construction that should be taken into account when choosing the set of compatible base points $(h_{\rho})_{\rho}$ are the collection of contents $\tlambda_{\rho}$ and the collection of generalised different ideals $\goth{D}(\ta_{\rho}, \ok)$. Indeed, if $(h_{\rho})_{\rho}$ is a compatible set of base points for which $\goth{D}(a_{\rho}, L) = I$ is independent from $\rho$ then:
$$\lambda_{\rho} \cdot m_{\rho}^{r+1} \cdot \frac{h_{\rho}}{\norm{h_{\rho}}} = \frac{L}{\norm{L}} \cdot \tlambda_{\rho} \cdot \frac{\goth{D}(\ta_{\rho}, \ok)}{I}$$
where everything on the left-hand side depends on $h_{\rho}$ and nothing on the right-hand side does. If $\tlambda_{\rho}, \goth{D}(\ta_{\rho}, \ok)$ were independent of $\rho \in \goth{S}_{r}$ then clearly one would be able to construct a set of compatible base points by setting $h_{\rho} = h$ for any admissible base point $h$. In general, we need to handle these two parameters to construct a set of compatible base points. For the remainder of this paper, we shall write $\tildeD_{\rho}$ for the ideal $\goth{D}(\ta_{\rho}, \ok)$, which we call the ``overflow ideal'' associated to the unit system $[\eps_{\rho(1)}|\dots|\eps_{\rho(r)}]$. In the next section we will study the impact of these parameters on the construction.

\subsubsection{The generalised different ideals and the matrices $\matrixA, \tmatrixA$}\label{sectiondifferentexplicit}

Before we move on to the construction of sets of compatible base points, we prove an important lemma on the link between the generalised different ideals and the matrices $\matrixA, \tmatrixA$ we introduced in \refp{defatable} and \refp{deftildeatable}.  We wish to represent the ideals $\goth{D}(a, L)$ and $\goth{D}(\ta, \ok)$ in some $\qq$-bases of $\kk$. This is done by generalising \cohendifferent\, as follows:

\begin{general}{Lemma}\label{calculdifferent}
Let $B = [e_0, \dots, e_{r+1}]$ be a $\zz$-basis of a fractional ideal $I$ and $B' = [e'_0, \dots, e'_{r+1}]$ be any $\qq$-basis of $\kk$. Let $f$ be a non-zero $\qq$-linear form on the $\qq$-vector space $\kk$. We may define $\mathcal{M} \in M_{r+2}(\qq)$ to be the matrix whose coefficients are the rational numbers $f(e_i e'_j)$ for $0 \leq i,j \leq r+1$. Then $\goth{D}(f, I)^{-1}$ is explicitly parametrised in the basis $B'$ as the image of $\zz^{r+2}$ under the $\qq$-linear isomorphism $\iota : \qq^{r+2} \to \kk$ defined by:
$$\iota\begin{pmatrix} x_0 \\ x_1 \\ \vdots \\ x_{r+1}  \end{pmatrix} = (e'_0, \dots, e'_{r+1}) \times \mathcal{M}^{-1} \times \begin{pmatrix} x_0 \\ x_1 \\ \vdots \\ x_{r+1}  \end{pmatrix}.$$
\end{general}

\begin{proof}
The proof is essentially the same as the proof given in \cohendifferent\, for the usual different ideal. Let $\mathcal{M}$ be the matrix $(f(e_ie'_j))_{0 \leq i,j \leq r+1}$. Let  us denote by $m_{ij}$ the coefficients of $\mathcal{M}^{-1}$ so that:
$$\sum_{k = 0}^{r+1} f(e_je'_k)m_{kl} = \delta_{jl}$$
where $\delta_{jl}$ is Kronecker's symbol. Let $x = \iota((x_0, \dots, x_{r+1})^{T})$ and $y = \sum_{j = 0}^{r+1} y_j e_j$ where the $x_i$'s are rational numbers and the $y_i$'s are integers. Then:
\begin{align*}
f(xy) & = f\left(\left(\sum_{j = 0}^{r+1} y_je_j\right) \times \left(\sum_{k = 0}^{r+1} e'_k \sum_{l = 0}^{r+1} m_{kl}x_l\right)\right) \\
f(xy) & =\sum_{j = 0}^{r+1} y_j \sum_{l = 0}^{r+1} \left(\sum_{k = 0}^{r+1} f(e_je'_k)  m_{kl}\right) x_l \\
f(xy) & =\sum_{j = 0}^{r+1} y_j x_j
\end{align*}
Therefore, $f(xy) \in \zz$ for all $y \in I$ if and only if $x_j \in \zz$ for all $0 \leq j \leq r+1$ and $\goth{D}(f, I)^{-1} = \iota(\zz^{r+2})$.
\end{proof}

Using lemma \ref{calculdifferent} we may represent the integral ideals $\goth{D}(a, L)$ and $\goth{D}(\ta, \ok)$ in the bases we have fixed previously. Indeed, recall that we have fixed $\BL = [e_0, \dots, e_{r+1}]$ a positive $\zz$-basis of $L$ and $\Bok = [\te_0, \dots, \te_{r+1}]$ a positive $\zz$-basis of $\ok$ in section \ref{sectiondetailledgeomsetup}. Thus, we define the matrix $\matrixm$ (resp. $\tmatrixm$) to be the matrix whose coefficients are $(a(e_i\te_j))_{0 \leq i,j \leq r+1}$ (resp. $(\ta(\te_i\te_j))_{0 \leq i, j \leq r+1}$). Since $a(L) \subset \zz$ and $\ta(\ok) \subset \zz$ these matrices are integral. It follows from lemma \ref{calculdifferent} that $\matrixm^{-1}$ represents the inverse different ideal $\goth{D} = \goth{D}(a, L)^{-1}$ in the basis $\BL$ and that $\tmatrixm^{-1}$ represents $\tildeD = \goth{D}(\ta, \ok)^{-1}$ in the basis $\Bok$. In particular, the norms of these ideals are related to the determinants of these matrices as:
$$|\det(\matrixm)| = \norm{\goth{D}} ~~ \text{ and } ~~ |\det(\tmatrixm)| = \norm{\tildeD}.$$
We shall now express a link between these two matrices and the matrices $\matrixA$ and $\tmatrixA$ defined in \refp{defatable} and \refp{deftildeatable}.

\begin{general}{Lemma}\label{lemmaAUM}
The two following matrix equalities hold:
$$(\mathbf{0}~\matrixA) = \tmatrixgothu^{T} \matrixm^{T}, ~~ \text{ and }~~ (\mathbf{0}~\tmatrixA) = \tmatrixgothu^{T} \tmatrixm^{T} $$
where $^T$ stands for transposition and $\tmatrixgothu$ is the matrix defined in \refp{deftildematu}. As a consequence, $\tlambda \mid s$ and $\tlambda \mid \ts$.
\end{general}

\begin{proof}
Recall that the matrix $(\mathbf{0}~\matrixA)$ is given by its coefficients $a(u_j e_l)$ for $0 \leq j \leq r$ and $0 \leq l \leq r+1$. Replacing $u_j$ by its expression $u_j = \sum_{k = 0}^j \tc_{jk0} \te_k$ gives $a(u_j e_l) = \sum_{k = 0}^j \tc_{jk0}a(\te_ke_l)$. Thus:
$$(\mathbf{0}~\matrixA) =  \begin{pmatrix} 1 & 0 & 0 & \dots & 0 & 0\\ \tilde{c}_{100} & \tilde{c}_{110} & 0 & \dots & 0 & 0 \\ \vdots & \vdots & \vdots & \vdots & \vdots & \vdots \\ \tilde{c}_{r00} & \tilde{c}_{r10} & \tilde{c}_{r20} & \dots & \tilde{c}_{rr0} & 0\end{pmatrix} \begin{pmatrix} a(\te_0 e_0) &\dots & a(\te_0 e_l) & \dots & a(\te_0 e_{r+1}) \\ \vdots & \vdots & \vdots & \vdots & \vdots \\ a(\te_j e_0) &\dots & a(\te_j e_l) & \dots & a(\te_j e_{r+1}) \\ \vdots & \vdots & \vdots & \vdots & \vdots \\
a(\te_{r+1} e_0) &\dots & a(\te_{r+1} e_l) & \dots & a(\te_{r+1} e_{r+1}) \\\end{pmatrix}  $$
which is exactly the statement that $(\mathbf{0}~\matrixA) = \tmatrixgothu^{T} \matrixm^T$. Let us now analyse this relation more precisely, by introducing two auxiliary square matrices $\tmatrixgothu_0$ and $\matrixm_0$ of size $r+1$ such that
$$\tmatrixgothu = \begin{pmatrix} & & \\ & \tmatrixgothu_0 & \\ & & \\ \hline 0 & \dots & 0\end{pmatrix} \text{ and } \matrixm = \begin{pmatrix} 0 & \dots & 0 & \rvline & \pm1 \\ \hline & & & \rvline & \ast \\  & \matrixm_0 & & \rvline & \vdots \\ & & & \rvline & \ast \end{pmatrix}.$$
Then the relation $(\mathbf{0}~\matrixA) = \tmatrixgothu^{T} \matrixm^T$ reduces to $\matrixA = \tmatrixgothu_0^{T}\matrixm_0^T$ and taking determinants gives 
$$s = \det(\matrixA) = \det(\tmatrixgothu_0^T) \times \det(\matrixm_0^T) = \tlambda.\det(\matrixm_0) \in \tlambda. \zz.$$
Thus $\tlambda \mid s$. The proof that $(\mathbf{0}~\tmatrixA) = \tmatrixgothu^T \tmatrixm^T$ and that $\tlambda \mid \ts$ is essentially the same. 
\end{proof}

\subsection{The choice of base points $(h_{\rho})_{\rho \in \goth{S}_r}$}\label{choiceh}

In this section, we explain how the base points $(h_{\rho})_{\rho \in \goth{S}_r}$ should be chosen in our Main Conjecture \ref{conjecture}. To make this choice, we keep in mind two important aspects:
\begin{itemize}
\item Computations show that for fields of degree $n \geq 4$, the base points $(h_{\rho})_{\rho}$ cannot be chosen independently. To ensure that the evaluation of our higher elliptic unit is \textit{correct}, one must make sure that the base points are compatible in the sense of Definition \ref{defcompatible}.
\item To perform efficient computations, and to get the simplest formulas for higher elliptic units, the values of the parameters $t_{\rho}$ should be the smallest possible, and, after that, the value of the levels $\ell_{\rho}$ should also be the smallest possible.
\end{itemize}
With these two points in mind, we shall explain how to choose the set of base points $(h_{\rho})_{\rho}$ under some specific conditions (see \refp{hypothesisprimediff}, \refp{hypothesisconstantclass} below). The general idea is that for any $\rho \in \goth{S}_r$, the admissible base point $h_{\rho}$ should be chosen if possible as a generator of the ideal $\frac{qN}{\goth{a}\goth{b}} \tildeD_{\rho}$. However, there is no guarantee that this ideal is principal, nor that it is generated by an admissible base point. For instance, it is not possible to choose an admissible base point of this form whenever $\tildeD_{\rho}$ and $\goth{f}$ share a common prime factor. For this reason, we shall always assume in the rest of this paper the following hypothesis:
\refstepcounter{hypocounter}
\begin{equation}\label{hypothesisprimediff}\tag{H\thehypocounter}
\text{The overflow ideal } \tildeD_{\rho} \text{ is coprime to } \goth{f}.
\end{equation}
This is essentially a weaker form of the hypothesis [\hspace{1sp}\cite{thesis}, (H4)] and it depends only on the unit system $[\eps_{\rho(1)}|\dots|\eps_{\rho(r)}]$ and on $\goth{f}$. Another clear obstruction arises when it is impossible to find positive integers $m_{\rho}$ for which the class of the ideal $m_{\rho} \cdot \tildeD_{\rho}$ is constant (meaning independent of $\rho$) in the wide class group at $\goth{f}$. Let us give a criterion for this to be the case. We define similarly to $\clzfplus$ the class group $\clzf$ in the following way: let $I(\goth{f})$ be, as usual, the set of fractional ideals in $\kk$ which are coprime to $\goth{f}$. Let $P(\goth{f}, \zz)$ be the subset of $I(\goth{f})$ consisting of principal fractional ideals generated by some $\beta \in \kk^{\times}$ for which there is some integer $k \in \zszx{q}$ satisfying $k\cdot\beta \equiv 1 \mod \goth{f}$. Then we set $\clzf = I(\goth{f})/P(\goth{f}, \zz)$ which is a quotient of the wide ray class group $\wideclassgroupf$. Since the action of integers coprime to $\goth{f}$ on this class group is trivial, it is easily seen that:
\begin{center}
There are integers $m_{\rho} \in \zszx{q}$ such that the class of the ideal $m_{\rho} \cdot \tildeD_{\rho}$ in  $\wideclassgroupf$ \\
is independent of $\rho \in \goth{S}_{n-2}$. \\
$\Updownarrow$ \\
The class of the ideal $\tildeD_{\rho}$ in $\clzf$ is independent of $\rho \in \goth{S}_{n-2}$.
\end{center}
Hence we shall also always assume for the rest of this paper that the latter holds, that is:
\refstepcounter{hypocounter}
\begin{equation}\label{hypothesisconstantclass}\tag{H\thehypocounter}
\text{The class of the ideal } \tildeD_{\rho} \text{ in the class group } \clzf \text{ is independent of } \rho \in \goth{S}_{n-2}.
\end{equation}
Note that lemma \ref{lemmabasicdifferent}, (vi) implies that the class of the ideal $\tildeD_{\rho}$ in the usual class group of $\kk$ is independent of $\rho$, so it is not hard to believe that this hypothesis is satisfied quite often. Especially, when $\ok/\goth{f}$ is cyclic (this is the case for most simple examples) then $\clzf \simeq \classgroupk$ and therefore this hypothesis is satisfied. In particular, the hypothesis \refp{hypothesisconstantclass} is a weaker version of [\hspace{1sp}\cite{thesis}, (H5)]. 

Once this is sorted, we still need to handle the possibility that the ideal $\frac{N}{\goth{a}\goth{b}} \tildeD_{\rho}$ might not belong to the trivial class of the wide class group at $\goth{f}$. To deal with this problem we use ideals of the form $m\cdot\helper$ where $m > 0$ is an integer coprime to $q$ and $\helper$ is an ideal satisfying $\ok/\helper \simeq (\zsz{\helperdown})^{r+1}$ where $\helper \cap \zz = \helperdown \zz$. The ideal $\helper$ is what we call a \textit{helper} ideal (see Definition \ref{defhelper}). We shall show that it is always possible (under hypotheses \refp{hypothesisprimediff} and \refp{hypothesisconstantclass}) to choose the base points as generators of ideals of the form
\begin{equation}\label{formcompatible}\frac{qN}{\goth{a}\goth{b}} \tildeD_{\rho}\cdot m_{\rho}\cdot \helper = h_{\rho} \ok
\end{equation}
for some integers $m_{\rho} > 0$ and some helper ideal $\helper$ independent of $\rho$. The presence of the ideal $m_{\rho}\cdot \helper$ as a divisor of the ideal generated by $h_{\rho}$ only has a little and well-controlled impact on the geometric setup we have presented in the previous sections. The fact that the ideal $\helper$ should not depend on the permutation $\rho$ is crucial to ensure that the base points $(h_{\rho})_{\rho}$ are compatible. We may now state our main result on the existence of compatible sets of base points:

\begin{general}{Theorem}\label{theoremchoiceh}
Suppose that $\kk$ is an ATR field of degree $n = r+2 \geq 3$. Fix an integral ideal $\goth{f} \neq \ok$. Set $q\zz = \goth{f} \cap \zz$ with $q >0$ and $Q > 0$ to be the maximal integer satisfying $\goth{f} \subset Q\ok$. Fix a smoothing ideal $\goth{a}$ of norm $N$ coprime to $q$. Fix an integral ideal $\goth{b}$ coprime to $\goth{f}$. Fix a set of non-degenerate fundamental units $\eps_1, \dots, \eps_r$ for $\uzfplus$. Let $\tildeD_{\rho}$ denote the overflow ideal of the unit system $[\eps_{\rho(1)}|\dots|\eps_{\rho(r)}]$ for $\rho \in \goth{S}_{r}$. Let us assume that for all $\rho$ the overflow ideal $\tildeD_{\rho}$ is coprime to $\goth{f}$ (see \refp{hypothesisprimediff}) and to $N$. Let us also suppose that the collection of overflow ideals $\tildeD_{\rho}$ satisfy \refp{hypothesisconstantclass} (or simply that $\ok/\goth{f}$ is cyclic). Then there is a helper ideal $\helper$ coprime to $q\cdot N\cdot\prod_{\rho}\ttt_{\rho}$ and there are integers $m_{\rho} > 0$ such that the fractional ideals  
\begin{equation}\label{eqthdefh}\frac{qN}{\goth{a}\goth{b}} \tildeD_{\rho}\cdot m_{\rho}\cdot \helper
\end{equation}
are generated by admissible base points $h_{\rho}$ which are compatible in the sense of Definition \ref{defcompatible}. In addition:
$$\begin{cases} \goth{D}(a_{\rho}, L) & = \frac{\goth{f}}{Q}\cdot \frac{\helperdown}{\helper} \cdot\goth{a}, \\ 
\lambda_{\rho} & = \tlambda_{\rho} \cdot\norm{\tildeD_{\rho}} \cdot\frac{Q q^{n-1}}{\norm{\goth{f}}} \cdot N^{n-2} \cdot \helperdown^{n-2}, \\
 \ell_{\rho} t_{\rho} & = \frac{q}{Q} \cdot N \cdot \helperdown \cdot \ttt_{\rho} \\ \ell_{\rho} & \in N \cdot \helperdown \,\zz  \\ \end{cases} $$
where $\helper \cap \zz = \helperdown \zz$. If we assume further that the integral ideal $I = \goth{f}/Q$ is such that $\ok/I$ is cyclic then $\ell_{\rho} \in \frac{q}{Q} \cdot N \cdot \helperdown \, \zz$ and $\ttt_{\rho}$ divides $t_{\rho}$. Lastly, in the special case where we assume that $\ok/\goth{f}$ is cylic, that $\tlambda_{\rho} = 1$ and that $\ttt_{\rho}$ is coprime to $q$, we have $\ell_{\rho} = q \cdot N \cdot \helperdown \cdot \ttt_{\rho}$ and $t_{\rho} = s_{\rho} = 1$.
\end{general}

Let us make a remark on the assumption that $N$ is coprime to $\tildeD_{\rho}$ for each $\rho \in \goth{S}_{n-2}$. This is tied to the condition that the smoothing operation we are performing is \textit{good} in the sense of \cite{secondpaper}. In practice, we start by computing the overflow ideals and only then do we fix a smoothing ideal satisfying this condition. 

Let us now give an outline of this section. In section \ref{sectionexistence}, we prove the first claim of Theorem \ref{theoremchoiceh}, that is, the existence of a compatible set of base points satisfying \refp{eqthdefh}. The exact description of the parameters $\lambda_{\rho}$ and $\goth{D}(a_{\rho}, L)$ will follow from this, as well as from lemma \ref{lemmalinktarget}. Then, in section \ref{sectiondeterminationlt}, we push our analysis further to obtain information on both $\ell_{\rho}$ and $t_{\rho}$.

\subsubsection{Existence of compatible sets of base points}\label{sectionexistence}

In this section we show the existence of compatible sets of base points satisfying \refp{formcompatible}. To this end we introduce the notion of \textit{helper ideals}: these are certain integral ideals $\helper$ which don't influence the value of $t_{\rho}$ and such that the ideal $\frac{qN}{\goth{a}\goth{b}} \tildeD_{\rho} \cdot \helper$ is generated by a weakly admissible base point. These ideals are defined as follows:

\begin{general}{Definition}\label{defhelper}
An integral ideal $\helper$ is a helper ideal if either $\helper = \ok$ or $\helper$ is a primitive ideal coprime to $q\cdot N\cdot \prod_{\rho}\ttt_{\rho}$ such that $\helperdown/\helper$ is a cyclic group of order $\helperdown$, where $\helper \cap \zz = \helperdown\,\zz$. An extended helper ideal is an ideal of the form $m\cdot\helper$ where $m$ is a positive integer coprime to $q$ and $\helper$ is a helper ideal. 
\end{general}

The use of helper ideals gives us a bit more room to find admissible base points $h$, at the cost of an increase of the level $\ell$ by a factor $\helperdown$. Luckily, the use of helper ideals is enough to guarantee that we find a suitable ideal, as we explain in the following lemma:

\begin{general}{Lemma}\label{lemmacebotarevhelper}
Any class in the class group $\clzf$ contains infinitely many helper ideals of the form $p\ok \cdot \goth{p}^{-1}$ where $p$ is a rational prime and $\goth{p}$ is a prime ideal of norm $p$. 
\end{general}

\begin{proof}
Fix a class $\goth{C}$ in the class group $\clzf$. By Cebotarev's density theorem (see for instance [\hspace{1sp}\cite{Cox}, Theorem 8.17]), there are infinitely many prime ideals of prime norm in the class $\goth{C}^{-1}$, that is to say there are infinitely many pairs $(p, \goth{p})$ where $p$ is a rational prime and $\goth{p}$ is a prime ideal of norm $p$ in $\ok$ such that $\goth{p} \in \goth{C}^{-1}$. For any of these pairs, the integral ideal $p\ok \cdot \goth{p}^{-1}$ is a helper ideal in the class $\goth{C}$ of $\clzf$.
\end{proof}

In particular, this lemma implies that any class in $\clzf$ contains a helper ideal $\helper$ which is coprime to $q \cdot N \cdot \prod_{\rho} \ttt_{\rho}$. We are now ready to prove the existence of compatible sets of base points.

\begin{general}{Proposition}\label{propexistencebasepoints}
Under the hypotheses of Theorem \ref{theoremchoiceh} there is a helper ideal $\helper$ coprime to $q\cdot N \cdot \prod_{\rho} \ttt_{\rho}$ and there are integers $m_{\rho} > 0$ such that the fractional ideals  
\begin{equation}\label{eqpropdefh}\frac{qN}{\goth{a}\goth{b}} \tildeD_{\rho}\cdot m_{\rho}\cdot \helper
\end{equation}
are generated by admissible base points $h_{\rho}$ which are compatible in the sense of Definition \ref{defcompatible}. In addition:
$$\begin{cases} \goth{D}(a_{\rho}, L) & = \frac{\goth{f}}{Q}\cdot \frac{\helperdown}{\helper} \cdot\goth{a}, \\ \lambda_{\rho} & = \tlambda_{\rho} \cdot\norm{\tildeD_{\rho}} \cdot\frac{Q q^{n-1}}{\norm{f}} \cdot N^{n-2} \cdot \helperdown^{n-2} \\
\ell_{\rho}t_{\rho} &= \frac{q}{Q} \cdot N \cdot \helperdown \cdot \ttt_{\rho}\end{cases}$$
where $\helperdown\, \zz = \helper \cap \zz$.
\end{general}

\begin{proof}
Let us first prove that these ideals are generated by an admissible base point for each permutation $\rho$. It follows from the hypothesis \refp{hypothesisconstantclass} that the class of the ideal $\tildeD_{\rho}$ in the class group $\clzf$ is independent of $\rho$. Let us then denote by $\goth{C}$ the common class of the fractional ideals $\frac{N}{\goth{a}\goth{b}}\tildeD_{\rho}$ in $\clzf$. By lemma \ref{lemmacebotarevhelper} there exists a helper ideal $\helper$ coprime to $q\cdot N \cdot \prod_{\rho} \ttt_{\rho}$ in the class $\goth{C}^{-1}$. This implies that the fractional ideals $\frac{N}{\goth{a}\goth{b}} \cdot \tildeD_{\rho} \cdot \helper$ all belong to the trivial class in $\clzf$. This means that there are elements $g_{\rho} \in \kk$ coprime to $\goth{f}$ and integers $k_{\rho} \in \zszx{q}$ such that for each $\rho \in \goth{S}_{n-2}$:
$$\begin{cases} g_{\rho} \ok & =\frac{N}{\goth{a}\goth{b}} \cdot \tildeD_{\rho} \cdot \helper \\ g_{\rho} & \equiv k_{\rho} \mod \goth{f}.\end{cases} $$
For each permutation $\rho$, let us fix a positive integer $m_{\rho}$ satisfying $k_{\rho} m_{\rho} \equiv 1 \mod \goth{q}$. Then, the fractional ideals $m_{\rho} \cdot \frac{N}{\goth{a}\goth{b}} \cdot \tildeD_{\rho} \cdot \helper$ all belong to the trivial class of the wide class group at $\goth{f}$ as they are generated by the element $m_{\rho} g_{\rho} \equiv 1 \mod \goth{f}$. The base point $h_{\rho} = q \cdot m_{\rho} \cdot g_{\rho}$ is thus an admissible base point in the sense of Definition \ref{defadmissible} and it satisfies 
$$h_{\rho}\, \ok = \frac{qN}{\goth{a}\goth{b}} \tildeD_{\rho}\cdot m_{\rho}\cdot \helper. $$
Let us now prove that the base points $h_{\rho}$ are compatible. To this end we use lemma \ref{lemmalinktarget} as for each $\rho$:
$$\lambda_{\rho} \cdot\goth{D}(a_{\rho}, L) = \tlambda_{\rho}\cdot \frac{\norm{h_{\rho}/m_{\rho}}}{\norm{L}}\cdot\frac{m_{\rho}}{h_{\rho}}L\cdot\tildeD_{\rho}.$$ 
By definition of the base point $h_{\rho}$ we have:
$$\frac{m_{\rho}}{h_{\rho}}L\cdot\tildeD_{\rho} = \frac{\goth{f}}{q} \cdot \frac{\goth{a}}{N} \cdot \helper^{-1}$$
and
$$\frac{\norm{h_{\rho}/m_{\rho}}}{\norm{L}} = \frac{q^n}{\norm{\goth{f}}} \cdot N^{n-1} \cdot \norm{\tildeD_{\rho}}\cdot  \helperdown^{n-1}.$$
Putting everything together we obtain:
\begin{equation}\label{eqdeterlambda}\lambda_{\rho} \cdot\goth{D}(a_{\rho}, L) = \tlambda_{\rho}\cdot \frac{q^{n-1} \goth{f}}{\norm{\goth{f}}} \cdot N^{n-2} \cdot \helperdown^{n-2} \cdot \norm{\tildeD_{\rho}} \cdot \goth{a} \cdot \frac{\helperdown}{\helper}.
\end{equation}
Now, it follows from \ref{lemmabasicdifferent}, (vii) that $\goth{D}(a_{\rho}, L)$ is a primitive integral ideal, which means that $\goth{D}(a_{\rho}, L) \subset k \ok$ for some integer $k$ implies that $k = \pm 1$. This means that the largest integer dividing the right-hand side of \refp{eqdeterlambda} is precisely $\lambda_{\rho}$. We note that if $Q$ is the largest integer satisfying $\goth{f} \subset Q\ok$ then:
$$\frac{q^{n-1} \goth{f}}{\norm{\goth{f}}} = \frac{(q/Q)^{n-1} \goth{f}/Q}{\norm{\goth{f}/Q}} $$
and $\norm{\goth{f}/Q}$ is a divisor of $(q/Q)^{n-1}$. Now, since the integral ideals $\goth{f}/Q$, $\helperdown/\helper$ and $\goth{a}$ are mutually coprime and primitive integral we obtain the desired equalities:
$$\begin{cases} \goth{D}(a_{\rho}, L) & = \frac{\goth{f}}{Q}\cdot \frac{\helperdown}{\helper} \cdot\goth{a}, \\ \lambda_{\rho} & = \tlambda_{\rho} \cdot\norm{\tildeD_{\rho}} \cdot\frac{Q q^{n-1}}{\norm{f}} \cdot N^{n-2} \cdot \helperdown^{n-2}\end{cases}$$
where $\helperdown \,\zz = \helper \cap \zz$.
In particular, the generalised different ideal $\goth{D}(a_{\rho}, L)$ is independent of $\rho$ and thus the base points $(h_{\rho})_{\rho}$ are compatible. Lastly, plugging the value we have obtained for $\lambda_{\rho}$ in the fundamental relation \refp{formulafundam} gives the equality:
$$\ell_{\rho}\cdot t_{\rho} = \frac{q}{Q} \cdot N \cdot \helperdown \cdot \ttt_{\rho}.$$
This completes the proof. 
\end{proof}

\textbf{Remarks:}
\begin{itemize}
\item This proposition is the main tool we use to produce sets of compatible base points $(h_{\rho})_{\rho}$ algorithmically and such sets are the ones we generally use in our computations to put our conjectures to the test. 
\item When one of the hypotheses \refp{hypothesisprimediff}, \refp{hypothesisconstantclass} are not satisfied, minor tweakings to the form \refp{formcompatible} should also allow us to find simple sets of compatible base points. Some examples to showcase what sort of computations remain possible are presented in [\hspace{1sp}\cite{thesis}, Section IV.2.6].
\item When choosing a helper ideal, one should try to minimise the value of the positive integer $\helperdown$ defined by $\helper \cap \zz = \helperdown \,\zz$. Indeed, as we explain in section \ref{computingarithmetic}, the complexity of the computation of the arithmetic $G_r$ function \refp{arithmeticgrshape} is linear in $\ell_{\rho}$. In section \ref{sectiondeterminationlt} we show that $\helperdown \mid \ell_{\rho}$, thus the complexity of the evaluation of the higher elliptic unit using the helper ideal $\helper$ is linear in $\helperdown$.
\end{itemize}

In the next section we take a deeper dive into the geometric setup associated to a set of compatible base points as given by Proposition \ref{propexistencebasepoints} and obtain the values of the level $\ell_{\rho}$ and of the parameter $t_{\rho}$ attached in some simple cases.

\subsubsection{Determination of $\ell_{\rho}$ and $t_{\rho}$}\label{sectiondeterminationlt}

In this section we work once again with fixed $\rho$ and drop the subscripts to gain in clarity. We shall assume as before that the overflow ideal $\tildeD$ of the unit system $u_1, \dots, u_r = [\eps_{\rho(1)}|\dots|\eps_{\rho(r)}]$ is coprime to $\goth{f}$ and to $N = \norm{\goth{a}}$. Our goal is to show that $\ell$ is divisible by $q \cdot N \cdot \helperdown$ (see lemma \ref{lemmalmin}). It will then follow that $t$ is a divisor of $\ttt/Q$. In the specific case where $\tlambda = 1$ and $\ttt$ is coprime to $q$ we shall show further that $\ell = q \cdot N \cdot \helperdown \cdot \ttt$ and $t = 1$ (see lemma \ref{lemmalt}).

To get to these results, we must start by giving a characterisation of the family $\talpha_0, \dots, \talpha_r$ introduced in \refp{deftalpha} in terms of the different ideal $\tildeD$. 

\begin{general}{Lemma}\label{talphabasediff}
The matrix $\ttt \cdot \tmatrixm^{-1}$ is integral and satisfies:
$$\ttt\cdot\tmatrixm^{-1} = \begin{pmatrix} \tx_0 & \dots & \tx_r & \rvline & \pm \ttt \\ \hline & & & \rvline & 0 \\  & \tmatrixB \cdot \tmatrixgothu_0^{T} & & \rvline & \vdots \\ & & & \rvline & 0 \\ \end{pmatrix}$$
where $\tmatrixgothu_0$ is the square matrix of size $r+1$ formed by the top $r+1$ rows of $\tmatrixgothu$ and $\tx_0, \dots, \tx_r$ are integers. As a consequence, the fractional ideal $\ttt \cdot \tildeD^{-1}$ is integral and it admits the $\zz$-basis $[\tbeta_0, \dots, \tbeta_r, \ttt]$ where for any $0 \leq k \leq r$:
$$\tbeta_k = \tx_k + \sum_{j = k}^r \tc_{jk0}\talpha_j.$$
In addition $\ttt \in \tildeD$ and 
\begin{equation}\label{eqtttn}\ttt^{n-1}  = \tlambda \cdot \norm{\tildeD} \cdot |\det(\tmatrixB)|.
\end{equation}
\end{general}

\begin{proof}
Let us write $\tmatrixm$ as:
$$\tmatrixm = \begin{pmatrix} &\rvline & & & \\ X & \rvline & & M & \\ &\rvline & & & \\ \hline \rule{0pt}{0.9\normalbaselineskip} Z & \rvline & & Y &\end{pmatrix}$$
where $M$ is a square matrix of size $r+1$. By definition of $\tmatrixm$, the coefficient $Z = \ta(\te_{r+1})$ is $\pm 1$. The matrix $\tmatrixgothu$ has a last row filled with zeroes, thus we denote as before by $\tmatrixgothu_0$ the square matrix formed by its top $r+1$ rows. Then, since $\tmatrixm$ is symmetric, we may rewrite the equality $(\mathbf{0}~\tmatrixA) = \tmatrixgothu^{T} \tmatrixm^{T}$ as $\mathbf{0} = \tmatrixgothu_0^T X$ and $\tmatrixA = \tmatrixgothu_0^T M$. The matrix $\tmatrixgothu_0$ is invertible when seen as a matrix with coefficients in $\qq$ as $\det(\tmatrixgothu_0) = \tlambda \neq 0$. Thus we obtain $X = 0$ and 
$$\tmatrixB \cdot \tmatrixgothu_0^T = \ttt \cdot \tmatrixA^{-1} \cdot \tmatrixgothu_0^T = \ttt \cdot M^{-1}. $$
Hence $\ttt \cdot M^{-1}$ is an integral matrix and
$$\ttt \cdot \tmatrixm^{-1} = \ttt \cdot \begin{pmatrix} & \pm YM^{-1} & & \rvline & \pm 1 \\ \hline & & & \rvline & 0 \\  & M^{-1} & & \rvline & \vdots \\ & & & \rvline & 0 \\ \end{pmatrix}=\begin{pmatrix} \tx_0 & \dots & \tx_r & \rvline & \pm \ttt \\ \hline & & & \rvline & 0 \\  & \tmatrixB \cdot \tmatrixgothu_0^{T} & & \rvline & \vdots \\ & & & \rvline & 0 \\ \end{pmatrix}$$
for some integers $\tx_0, \dots, \tx_r$. As a consequence the matrix $\ttt \cdot \tmatrixm^{-1}$ is integral, and since it represents the fractional ideal $\ttt \cdot \tildeD^{-1}$ in the basis $\Bok$ of $\ok$ this ideal is integral. As a consequence the integral ideal $\tildeD$ is a divisor of $\ttt \ok$ and $\ttt \in \tildeD$. In addition, taking determinants gives the equality:
$$\ttt^{n-1} = |\det(\tmatrixgothu_0)| \cdot |\det(\tmatrixm)| \cdot |\det(\tmatrixB)|   = \tlambda \cdot \norm{\tildeD} \cdot |\det(\tmatrixB)|.$$

Let us now denote by $\tbeta_k$ the element of $\ok$ corresponding to the $(k+1)$-th column of $\ttt \cdot \tmatrixm^{-1}$ for $0 \leq k \leq r$. Then by definition of $\tmatrixm$, the family $[\tbeta_0, \dots, \tbeta_r, \ttt]$ is a $\zz$-basis of $\ttt \cdot \tildeD^{-1}$. We now only need to make the link between this family and the family $\talpha_0, \dots, \talpha_r$ explicit. Recall that the columns of the matrix $\tmatrixB$ represent the elements $\talpha_0, \dots, \talpha_r$ and that:
$$\tmatrixgothu_0^T = \begin{pmatrix} \tilde{c}_{000} = 1 & 0 & 0 & \dots & 0 & 0\\ \tilde{c}_{100} & \tilde{c}_{110} & 0 & \dots & 0 & 0 \\ \vdots & \vdots & \vdots & \vdots & \vdots & \vdots \\ \tilde{c}_{r00} & \tilde{c}_{r10} & \tilde{c}_{r20} & \dots & \tilde{c}_{rr0} & 0\end{pmatrix}.$$
Therefore, we get $\tbeta_k = \tx_k + \sum_{j = k}^r \tc_{jk0} \talpha_j$ as claimed.
\end{proof}

The above result is a general version of [\hspace{1sp}\cite{thesis}, Lemma III.18] where we only dealt with the case where $\tlambda = 1$. We may now move on to the determination of a lower bound for the level $\ell$:

\begin{general}{Lemma}\label{lemmalmin}
There are integers $m_0, \dots, m_r$ such that $\talpha_0 + m_0, \dots, \talpha_r + m_r \in \goth{a} \cdot (\helperdown \cdot \helper^{-1})$. As a consequence, the level $\ell$ is divisible by $N \cdot \helperdown$. Moreover, if $q'$ is the smallest divisor of $q$ satisfying $\talpha_0, \dots, \talpha_r \in \zz + \goth{f}/q'$ then $\ell$ is divisible by $(q/q') \cdot N \cdot \helperdown$. In particular, if $I =\goth{f}/Q$ is such that $\ok/I$ is cyclic, $(q/Q) \cdot N \cdot \helperdown$ divides $\ell$.
\end{general}

\begin{proof}
The groups $\ok/\goth{a} \simeq \zsz{N}$ and $\helperdown \cdot \helper^{-1} \simeq \zsz{\helperdown}$ are cyclic, therefore any of the $\talpha_j$ is congruent to an integer modulo $\goth{a}$ and to an integer modulo $\helperdown \cdot \helper^{-1}$. Since $\goth{a}$ and $\helperdown \cdot \helper^{-1}$ are coprime by assumption, it follows from the Chinese Remainder Theorem that there are integers $m_0, \dots, m_r$ such that $\talpha_0 + m_0, \dots, \talpha_r + m_r \in \goth{a} \cdot \helperdown \cdot \helper^{-1}$. This gives for any $0 \leq j \leq r$, $(\talpha_j + m_j) h \in m \cdot N \cdot \helperdown L$ and by definition of $\ell$ this implies that $N \cdot \helperdown \mid \ell$. 

Let us now discuss the  $q$-part of $\ell$. Since $\goth{f} \cap \zz = q\zz$ the fractional ideal $q\ok \cdot \goth{f}^{-1} = q/\goth{f}$ is integral and therefore its inverse $\goth{f}/q$ contains $\ok$. Since $\talpha_0, \dots, \talpha_r$ are elements of $\ok$ we have $\talpha_0, \dots, \talpha_r \in \zz + \goth{f}/q$. Let $J = \{ k \in \zz - \{ 0\} \setseparator \talpha_0, \dots, \talpha_r \in \zz + \goth{f}/k\}$. This set is non-empty as $q \in J$. In addition, if $k \in J$ and $k' \in \zz -\{0\}$ then $kk' \in J$ and if $k, k' \in J$, the gcd $(k, k')$ belongs to $J$. Thus $J = q' (\zz - \{0\})$ for a unique positive divisor $q'$ of $q$. Since $\goth{f}/q', \goth{a}, \helperdown$ are two-by-two coprime, by the Chinese remainder theorem we get that $\talpha_0, \dots, \talpha_r \in \zz + (\goth{f}/q') \cdot \goth{a} \cdot (\helperdown \cdot \helper^{-1})$. It follows that $\talpha_0h, \dots, \talpha_rh \in \zz h + m\cdot (q/q') \cdot N \cdot \helperdown L$ which implies that $(q/q')\cdot N \cdot \helperdown \mid \ell$. Lastly, if $I = \goth{f}/Q$ is such that $\ok/I$ is cyclic, then $Q \in J$ and therefore $q'$ divides $Q$. Thus $q/Q$ divides $q/q'$ and $(q/Q) \cdot N \cdot \helperdown \mid \ell$.
\end{proof}

It follows immediately from lemma \ref{lemmalmin} and from the equality $\ell t = (q/Q) \cdot N \cdot \helperdown \cdot \ttt$ that $t$ divides $\ttt q/Q$. Let us now focus on the specific case where $\ok/\goth{f}$ is cyclic, $\tlambda =1$ and $\ttt$ is coprime to $q$ and let us give a complete determination of $\ell$ and $t$ in this case.

\begin{general}{Lemma}\label{lemmalt}
Suppose that $\ok/\goth{f}$ is cyclic, that $\tlambda = 1$ and that $\ttt$ is coprime to $q$. Then $\ell = q \cdot N \cdot \helperdown \cdot \ttt$ and $t = s = 1$.
\end{general}

\begin{proof}
Since $\tlambda = 1$, we have $\tmatrixgothu_0 = I_{r+1}$ and $\tbeta_k = \tx_k + \talpha_k$ in lemma \ref{talphabasediff}. Thus, the ideal $\ttt \cdot \tildeD^{-1}$ admits the $\zz$-basis $[\talpha_0+\tx_0, \dots, \talpha_r + \tx_r, \ttt\,]$. It then follows that 
$$(\talpha_j + \tx_j)h \in m \ttt L $$
and therefore by definition of $\ell$ we get $\ttt \mid \ell$. From lemma \ref{lemmalmin} we also have $q \cdot N \cdot \helperdown \mid \ell$ as $Q= 1$ when $\ok/\goth{f}$ is cyclic. By assumption, $\ttt$ is coprime to $q \cdot N \cdot \helperdown$, therefore $q \cdot N \cdot \helperdown \cdot \ttt \mid \ell$. Since $ \ell t = q \cdot N \cdot \helperdown \cdot \ttt$ we conclude that 
$$\begin{cases} \ell & = q \cdot N \cdot \helperdown \cdot \ttt \\ t & = 1.\end{cases}$$
Lastly, $s$ and $t$ share the same prime factors, therefore $s = 1$.
\end{proof}

We point out that the results of Proposition \ref{propexistencebasepoints}, lemma \ref{lemmalmin} and lemma \ref{lemmalt} together form a complete proof of Theorem \ref{theoremchoiceh}.

\section{Computing the conjecture}\label{sectioncomputing}

In this section we explain how we compute the term 
$$I_{r,\goth{f}, \conjko\cdot\goth{b}, \goth{a}}(\eps_1, \dots, \eps_r ; \bars{h}, \bars{\mu}, \bars{\nu}, (k_{\rho}, \eps_{\rho})_{\rho}) $$
appearing in the Main Conjecture \ref{conjecture}, which by \refp{formeexpliciteune} is explicitly given in terms of ordinary $G_r$ functions as:
\begin{multline}\label{explicitcomputation} I_{r,\goth{f}, \conjko\cdot\goth{b}, \goth{a}}(\eps_1, \dots, \eps_r ; \bars{h}, \bars{\mu}, \bars{\nu}, (k_{\rho}, \eps_{\rho})_{\rho}) = \\ \prod_{\rho \in \goth{S}_{n-2}}\prod_{\delta \in F_{\rho}/\zz\gamma_{\rho}} \prod_{(k', \eps') \in \zfoneplus} \left(\frac{G_r\left(\frac{\mu_{\rho} \cdot \conjko \cdot k_{\rho}k'_{\rho} \cdot \eta_{\rho} m_{\rho}}{q} + \frac{\delta}{\gamma_{\rho}}, \frac{\alpha_{0, \rho}}{\gamma_{\rho}}, \dots, \frac{\alpha_{r, \rho}}{\gamma_{\rho}}\right)^N}{G_r\left(\frac{N(\mu_{\rho} \cdot \conjko \cdot k_{\rho}k'_{\rho} \cdot \eta_{\rho} m_{\rho})}{q}+\frac{N\delta}{\gamma_{\rho}}, \frac{N\alpha_{0, \rho}}{\gamma_{\rho}}, \dots, \frac{N\alpha_{r, \rho}}{\gamma_{\rho}}\right)}\right)^{\nu_{\rho} \mu_{\rho}^r}.
\end{multline}
The computation of this product and our testing of the conjecture is discussed at length in \cite{thesis} and \cite{cheu}. We synthesize these discussions here, by giving an insight on how we compute $G_r$ functions in general, and more specifically our general strategy to compute \refp{explicitcomputation}. We also explain how we identify correct orientations datum using the expected Kronecker formula and how we proceed to recognise the obtained complex number $I_{r,\goth{f}, \conjko\cdot\goth{b}, \goth{a}}(\eps_1, \dots, \eps_r ; \bars{h}, \bars{\mu}, \bars{\nu}, (k_{\rho}, \eps_{\rho})_{\rho})^{\denomclasscohom}$ as a unit inside $\classfieldf$. 

\subsection{Computing the $G_r$ functions}\label{computinggr}

Let us now explain how we compute the ordinary $G_r$ functions, practically speaking, to verify the conjecture on numerical examples. Indeed, the definition of the $G_r$ functions by a multi-index infinite product given by formula \refp{defgr} is not suited for computations. In the case $n = 2$, Jacobi's triple product formula gives a beautiful expression of the $\theta$ function in terms of an infinite sum with converging rate in $q^{n^2/2}$:
\begin{equation}\label{jacobitriple}
\theta(z, \tau) = \frac{q^{1/24}}{\eta(\tau)} \sum_{n \in \zz} x^n (-1)^n q^{n(n-1)/2}
\end{equation}
where $x = \exp(2i\pi z)$, $q = \exp(2i\pi\tau)$ and $\eta(\tau)$ is Dedekind's $\eta$ function. This makes the computation of $\theta$ very fast when $\tau$ is not too close to the real axis, and standard techniques in the study of modular functions allow to use the modularity of $\theta$ to reduce the computation of $\theta(z, \tau)$ to the computation of $\theta(z', \tau')$ where $\Im(\tau') \geq 1/2$. Unfortunately, there is no clear generalisation of Jacobi's triple product formula for higher degree $G_r$ functions and we need other techniques to compute the $G_r$ functions efficiently. In the case where the parameters $\tau_0, \dots, \tau_r$ lie in the upper half-plane, we shall make use of [\hspace{1sp}\cite{Nishizawa}, Proposition 3.6] which we write as:
\begin{equation}\label{sumcomputations}
G_r(z, \tau_0, \dots, \tau_r) = \exp\left(-\sum_{j \geq 1} \frac{1}{j} \frac{q_0^j\dots q_r^jx^{-j} +(-1)^r x^{j}}{\prod_{k = 0}^r (1 - q_k^j)}\right)
\end{equation}
where $x = \exp(2i\pi z)$ and $q_k = \exp(2i\pi \tau_k)$ for $0 \leq k \leq r$.
This formula is only valid for $\tau_0, \dots, \tau_r \in \hh^{r+1}$ and $0 < \Im(z) < \sum_{k = 0}^r \Im(\tau_k)$. We call this domain the center strip. Using the properties of the $G_r$ functions, we may always reduce the computation of a term $G_r(z, \tau_0, \dots, \tau_r)$ with parameters $\tau_0, \dots, \tau_r \in \cc - \rr$ to the computation of one or more terms of the form $G_{r-k}(z', \tau'_0, \dots, \tau'_{r-k})$ with $k \geq 0$, $\tau'_0, \dots, \tau'_{r-k} \in \hh$ and $z'$ in the center strip (see \cite{thesis}, \cite{cheu}).

Let us say a few words on the number of terms we must compute in \refp{sumcomputations}. A simple computation shows that for any $J \geq 1$:
$$\left|\sum_{j \geq J} \frac{1}{j} \frac{q_0^j\dots q_r^jx^{-j} +(-1)^r x^{j}}{\prod_{k = 0}^r (1 - q_k^j)}\right| \leq \frac{|x|^{J} + |q_0\dots q_rx^{-1}|^J}{|1 - x|\prod_{k = 0}^r |1 - q_k|} $$
when $z$ lies in the center strip. This bound is then smaller than some $\delta > 0$ whenever 
$$J(\log|x| + \log|q_0\dots q_r/x|) -\log|1-x| - \sum_{k = 0}^r \log|1-q_k| \leq \log(\delta).$$
Thus, to compute $G_r(z, \tau_0, \dots, \tau_r)$ in the center strip with precision $\delta > 0$ we need to compute $J = O(\log(\delta)/y)$ terms in the sum \refp{sumcomputations}, where $y = 2\pi.\min(\Im(z), \Im(\tau_0) + \dots + \Im(\tau_r) - \Im(z))$ is the distance from $z$ to the boundary of the center strip. The main difficulty in performing the computations to test the conjecture comes from the fact that in practice, the imaginary parts of all the $\tau_k$'s are very small, making the parameter $y$ even smaller. It follows from Proposition \ref{fundam}, (iii) that for our arithmetic evaluation attached to a unit system $u_1, \dots, u_r$ and to an admissible base point $h$, the parameter $y$ is $O(\ell^{-1})$ and therefore the computation time is essentially linear in the level $\ell$.

We also briefly remark that when attempting to compute $G_r(z_j, \tau_0, \dots, \tau_r)$ for several values of $z_1, \dots, z_J$ with identical imaginary parts, that is, $\Im(z_1) = \Im(z_2) = \dots = \Im(z_J)$, the convergence rate in the computation of \refp{sumcomputations} is the same, and we may gain time by computing the successive terms $j\cdot\prod_{k = 0}^r (1 - q_k^j)$ only once. 

As a last remark, the case where $\kk$ contains a proper subfield $\ff$ (which is then automatically totally real) should be treated with great caution, as some of the parameters $\tau_j$ involved in the computations may lie in $\rr$. In this case, we may use the alternative formula:
\begin{equation}\label{grreal}
G_r(z, \tau_0, \dots, \tau_r) = \begin{cases} \exp\left(\sum_{j \geq 1} \frac{1}{(2i)^rj} \frac{\sin(\pi j(2z - (\tau_0+\dots +\tau_r)))}{\prod_{k = 0}^{r} \sin(\pi j \tau_k)}\right) & \text{ if } r \text{ is odd} \\
\exp\left(\sum_{j \geq 1} \frac{2}{(2i)^{r+1}j} \frac{\cos(\pi j(2z - (\tau_0+\dots +\tau_r)))}{\prod_{k = 0}^{r} \sin(\pi j \tau_k)}\right) & \text{ if } r \text{ is even} \end{cases}
\end{equation}
to compute the $G_r$ function. We refer to [\hspace{1sp}\cite{thesis}, Section I.5] for a discussion on that matter.

\subsection{Computing the elliptic unit}\label{computingarithmetic}

In this section we explain step by step how we compute the term $I_{r,\goth{f}, \conjko\cdot\goth{b}, \goth{a}}(\eps_1, \dots, \eps_r ; \bars{h}, \bars{\mu}, \bars{\nu}, (k_{\rho}, \eps_{\rho})_{\rho})$ appearing in the Main Conjecture \ref{conjecture}, at least under the hypotheses \refp{hypothesisprimediff} and \refp{hypothesisconstantclass}. 

To perform our computations we used the software Pari/GP \cite{parigp} and we made extensive use of the tools it provides for algebraic number theory. Our general strategy to compute higher elliptic units takes the following form:
\begin{enumerate}
\item Fix an ATR field $\kk$ and a class field modulus $\goth{f} \neq \ok$. Set $q\zz = \goth{f} \cap \zz$. Compute the narrow ray class group at $\goth{f}$. Compute the derivatives of partial zeta functions at $s = 0$, namely $\zetaf'(\goth{c}, 0)$ for all classes $\goth{c} \in \classgroupf$.
\item Compute fundamental units $\eps_1, \dots, \eps_r$ for the group $\uzfplus$. Check that they are non-degenerate (see section \ref{sectionevaluation}). Compute the contents $\tlambda_{\rho}$, the overflow ideals $\tildeD_{\rho}$ and the overflows $\ttt_{\rho}$. The computation of the ideal $\tildeD_{\rho}$ may be done using lemma \ref{calculdifferent}. Check that the overflow ideals $\tildeD_{\rho}$ satisfy \refp{hypothesisprimediff} and \refp{hypothesisconstantclass}. If needed, try changing fundamental units to get better properties for the fundamental cycle.
\item Fix a smoothing ideal $\goth{a}$ of norm $N$ such that $N$ is coprime to $q \cdot \prod_{\rho} \ttt_{\rho}$.
\item Compute the classgroupf $\clzfplus$ and fix representative ideals $\goth{b}_1, \dots, \goth{b}_g$ for its classes (which ought to be coprime to $q \cdot N \cdot \prod_{\rho} \ttt_{\rho}$). Fix integers $\conjko_1, \dots, \conjko_{g'}$ such that the class group $\classgroupf$ is represented by the ideals $\conjko_i \cdot \goth{b}_j$ for $1 \leq i \leq g'$ and $1 \leq j \leq g$. This representation of the class group allows us to reduce the number of computations we perform down the line.
\item Compute helper ideals $\helper(\goth{b}_i)$ and base points $(h_{\rho}(\goth{b}_i))_{\rho}$ as given by Theorem \ref{theoremchoiceh} for all $1 \leq i \leq g$.
\item Compute the parameters $\gamma_{\rho}(\goth{b}_i)$, $\alpha_{j, \rho}(\goth{b}_i)$ as given in section \ref{sectionbasicgeomsetuph}.
\item For any $1 \leq i \leq g$, any $\rho \in \goth{S}_{n-2}$ and any $k \in \zszx{q}$, $k > 0$ compute the term
$$\Omega(i, \rho, k) = \prod_{\delta \in F_{\rho}(\goth{b}_i)/\zz\gamma_{\rho}(\goth{b}_i)}\frac{G_r\left(\frac{k}{q} + \frac{\delta}{\gamma_{\rho}(\goth{b}_i)}, \frac{\alpha_{0, \rho}(\goth{b}_i)}{\gamma_{\rho}(\goth{b}_i)}, \dots, \frac{\alpha_{r, \rho}(\goth{b}_i)}{\gamma_{\rho}(\goth{b}_i)}\right)^N}{G_r\left(\frac{Nk}{q}+\frac{N\delta}{\gamma_{\rho}(\goth{b}_i)}, \frac{N\alpha_{0, \rho}(\goth{b}_i)}{\gamma_{\rho}(\goth{b}_i)}, \dots, \frac{N\alpha_{r, \rho}(\goth{b}_i)}{\gamma_{\rho}(\goth{b}_i)}\right)}.$$
To compute this efficiently we make use of the fact that the computation of multiple terms of the form $G_r(z_j, \tau_0, \dots, \tau_r)$ for several values $z_1, \dots, z_J$ with identical imaginary parts may be done at once (see section \ref{computinggr}). This is where we make the most out of our choice of representation of the narrow ray class group at $\goth{f}$ using the auxiliary class group $\clzfplus$ (see step 4).
\item For $1 \leq i \leq g$, $1 \leq j \leq g'$ identify a correct orientation datum by testing which values of $k_{\rho} \in \zszx{q}$ and $\epsilon_{\rho} = \pm 1$ give the desired Kronecker limit formula:
\begin{equation}\label{klfalgo} N\cdot \zetaf'(\conjko_j \cdot \goth{b}_i, 0) - \zetaf'(\conjko_j \cdot \goth{b}_i \cdot \goth{a}, 0) = \log\left|\prod_{(k', \eps')\in \zfoneplus}\prod_{\rho \in \goth{S}_{n-2}} \Omega(i, \rho, k'\cdot k_{\rho})^{\epsilon_{\rho}}\right|^2.
\end{equation}
This is done by testing for relations of the form \refp{klfalgo} using the \textbf{lindep} command in Pari/GP.
\item Once correct orientation data have been obtained for all classes, attempt to compute a minimal polynomial satisfied by the conjectural higher elliptic unit:
$$I_{r,\goth{f}, \conjko\cdot\goth{b}, \goth{a}}(\eps_1, \dots, \eps_r ; \bars{h}, \bars{\mu}, \bars{\nu}, (k_{\rho}, \eps_{\rho})_{\rho})^{\denomclasscohom}$$
 either over $\kk$ or other $\qq$ using versions of \textbf{lindep} or \textbf{algdep} commands in Pari/GP. An in depth-discussion on how we recognise the higher elliptic units as algebraic units may be found in \cite{thesis}.
\end{enumerate}

In practice, it is also possible to compute the higher elliptic unit for only one class to high precision and to try to recognise it as an algebraic number over $\qq$.

\section{Numerical Evidence}\label{sectionnumerical}

In this section, we provide numerical examples to support our conjecture. They may be computed with high precision in a low amount of time. In what follows, we will give computation times for $1000$ digits precision on a personal computer. Computations\footnotemark[1]\footnotetext[1]{A short version of the code for the computation of the ordinary $G_r$ functions is available at https://plmlab.math.cnrs.fr/pmorain/computations\_of\_higher\_elliptic\_units. On this page, one may find other examples of computations of higher elliptic units in various settings for number fields up to degree 6, and more examples might be added in the future.} were carried out using number fields found in the LMFDB database \cite{lmfdb} as well as the computer algebra system PARI/GP \cite{parigp}, making extensive use of algebraic number theory tools it provides.

In what follows, we will define our fields as $\kk = \qq(z)$ where $z$ is the complex root of some polynomial $P = (x-z)(x-\bar{z})\prod_{j = 1}^r (x-z_j) \in \qq[x]$ lying in the upper half-plane (here the $z_j$ are real numbers with $z_1 < \dots < z_r$). Thus, to define our orientations we will fix the ordering on the real embeddings of $\kk$ such that $\sigma_{j}(z) = z_j$ for $1 \leq j \leq r$. Throughout this section, we focus mainly on computations for $\goth{b} = \ok$ because most of the work on the computations of higher elliptic units does not depend on $\goth{b}$. For ease of presentation, we will often write prime ideals above a prime $p$ as $\goth{P}_p, \goth{P}'_p, \goth{P}''_p, \dots$. In this section we provide a few examples for fields of degree $3 \leq n \leq 5$ (more examples may be found in \cite{thesis} and \cite{cheu} for degrees up to $n = 6$). For ease of presentation we shall present examples where $\ok/\goth{f}$ is cyclic (except example \ref{cubicnoncyclic}) and the hypothesis \refp{hypothesisprimediff} is satisfied for the fundamental cycle. We also focus on the case where the contents $\tlambda_{\rho}$ are all equal to $1$ as it simplifies the formulas we write down for the higher elliptic units. Note that we present examples in non-optimal settings, as opposed to the examples provided in \cite{cheu}. Regarding the smoothing operation, we will focus on ``big smoothings'', that is on smoothing ideals of prime norm $N > n+1$ for which $\denomclasscohom = 1$. Lastly, for clarity, we will write absolute polynomials over $\qq$ as polynomials in the variable $x$ and relative polynomials over $\ok$ as polynomials in the variable $X$.

\subsection{Cubic examples}

In this section we present two examples for the complex cubic field $\qq(\sqrt[3]{2})$. The first example is one where $q = 29$ and $\zfoneplus$ contains $7$ elements. This showcases the computation of higher elliptic units for non-optimal settings in complex cubic fields with $\ok/\goth{f}$ cyclic. The second example is one of the simplest examples where $\ok/\goth{f}$ is not cyclic and where the hypothesis \refp{hypothesisprimediff} is satisfied. For this example $q = 4, Q = 2$ and $\zfoneplus$ contains only one element. 

\stepcounter{example}
\subsubsection{Example \exc}\label{detailledcubic}

Let $\zexc$ be be the complex root of the polynomial $x^3-2$ lying in the upper half-plane. The complex cubic field $\kk = \qq(\zexc)$ has usual class number $1$.  We choose the ideal $\goth{f}$ to be the unique prime ideal of norm $q = 29$ in $\ok$. The corresponding narrow ray class group is $\classgroupf \simeq \zsz{4}$. It is isomorphic to $\clzfplus$ and it is represented by the integral ideals $\ok, 2\ok, 4\ok, 8\ok$.  The unit group $\opck = \uzfplus$ is generated by $\eps_1 = \zexc -1$ and $\opcf$ is generated by $\eps_1^7$. We then have
$$\zfoneplus \simeq \{(1, 1), (7, \eps_1), (20, \eps_1^{2}), (24, \eps_1^{3}), (23, \eps_1^4), (16, \eps_1^5), (25, \eps_1^6) \} $$
and
$$\zfoneplus \simeq \zfoneplus \bigcup \{ (q-k, -\eps)\setseparator (k, \eps) \in \zfoneplus\}.$$
We may compute the geometric data attached to the fundamental cycle of $\uzfplus$ as $\tlambda = 1$, $\tildeD = \ok$ and $\ttt = 1$. Let us fix the smoothing ideal $\goth{a}$ as the unique prime ideal of norm $N = 5$ in $\ok$. The ideal $2 \cdot qN/\goth{a}$ is generated by the admissible base point $h = -58\zexc^2 +116\zexc+58$. Here we have $m = 2$. The corresponding level is $\ell= 5\cdot29 = 145$. We may compute for $\conjko = 1, 2, 4, 8$ and $\goth{b} = \ok$:
\begin{align*}
u_1 = \Gamma^{+}_{\goth{f}, 1, \ok, \goth{a}}(\eps_1, h) = & \prod_{k \in \{1, 7, 16, 20, 23, 24, 25\}}\frac{\Gamma\left(\frac{2\cdot k}{29}, \frac{\zexc^2+\zexc -557}{145}, \frac{\zexc-693}{145}\right)^{5}}{\Gamma\left(\frac{10\cdot k}{29}, \frac{\zexc^2+\zexc -557}{29}, \frac{\zexc-693}{29}\right)} \approx -0.1762247...+i\cdot 0.0267833...\\
u_2 = \Gamma^{+}_{\goth{f}, 2, \ok, \goth{a}}(\eps_1, h) = & \prod_{k \in \{2, 3, 11, 14, 17, 19, 21\}}\frac{\Gamma\left(\frac{2\cdot k}{29}, \frac{\zexc^2+\zexc -557}{145}, \frac{\zexc-693}{145}\right)^{5}}{\Gamma\left(\frac{10\cdot k}{29}, \frac{\zexc^2+\zexc -557}{29}, \frac{\zexc-693}{29}\right)} \approx -0.0564059... +i\cdot 0.0569654...\\
u_4 = \Gamma^{+}_{\goth{f}, 4, \ok, \goth{a}}(\eps_1, h) = & \prod_{k \in \{4, 5, 6, 9, 13, 22, 28\}}\frac{\Gamma\left(\frac{2\cdot k}{29}, \frac{\zexc^2+\zexc -557}{145}, \frac{\zexc-693}{145}\right)^{5}}{\Gamma\left(\frac{10\cdot k}{29}, \frac{\zexc^2+\zexc -557}{29}, \frac{\zexc-693}{29}\right)} \approx -5.5464541... - i\cdot0.8429711...\\
u_8 = \Gamma^{+}_{\goth{f}, 8, \ok, \goth{a}}(\eps_1, h) = & \prod_{k \in \{8, 10, 12, 15, 18, 26, 27\}}\frac{\Gamma\left(\frac{2\cdot k}{29}, \frac{\zexc^2+\zexc -557}{145}, \frac{\zexc-693}{145}\right)^{5}}{\Gamma\left(\frac{10\cdot k}{29}, \frac{\zexc^2+\zexc -557}{29}, \frac{\zexc-693}{29}\right)} \approx -8.7768188... -i\cdot 8.8638850...
\end{align*}
These four values coincide up to 1000 digits with four roots of the reciprocal polynomial 
$$P_{\mathrm{abs}} =  x^{12} + 27x^{11} + 336x^{10} + 1688x^9 + 1635x^8 - 5394x^7 + 4039x^6 -...$$ 
which defines an absolute equation of $\classfieldf$ over $\qq$.
Alternatively, the four computed values $u_1, u_2, u_4, u_8$ coincide up to 1000 digits of precision with the roots of the polynomial 
$$P_{\mathrm{rel}} = X^4 + 1 + (-7\zexc^2 + 9)(X^3+X) + (-30\zexc^2 + 15\zexc + 31)X^2 \in \ok[X]$$
which defines a relative equation of $\classfieldf$ over $\kk$. The computation time for $1000$ digits is 40 seconds. We may also check formula \refp{klfconjecture} as
$$ 5\cdot \zeta'_{\goth{f}}(\ok, 0) - \zeta'_{\goth{f}}(\goth{a}, 0) = \log\left|u_1\right|^2 \approx -1.7245769...$$
and similarly for the other classes. In addition, we can easily check that up to high precision the roots of $P_{\mathrm{abs}}$ are 
$$u_{\conjko}, \overline{u_{\conjko}} \text{ for } \conjko = 1, 2, 4, 8 $$
together with 4 roots $z_1, z_1^{-1}, z_2, z_2^{-1}$ which lie on the unit circle. This is coherent with the seventh point in the Main Conjecture \ref{conjecture}.

\stepcounter{example}
\subsubsection{Example \exc}\label{cubicnoncyclic}

Keep the same number field as in the previous example, that is $\zexc = e^{2i\pi/3} \cdot \sqrt[3]{2}$ and $\kk = \qq(\zexc)$. We choose the ideal $\goth{f} = \goth{P}_2^4$ where $\goth{P}_2$ is the unique prime ideal of norm $2$ in $\ok$. In particular $\goth{P}_2^3 = 2\ok$ and therefore $q = 4, Q = 2$. This is a case where $\ok/\goth{f} \simeq (\zsz{2})^{4}$ is not cyclic. The corresponding narrow ray class group is $\classgroupf \simeq \zsz{2}$. It is isomorphic to $\clzfplus$ and it is represented by the integral ideals $\ok, 3\ok$. The unit group $\opck$ is generated by $\eps_1 = (\zexc -1)$ and $\uzfplus = \opcf$ is generated by $\eps_1^4$. Since $\uzfplus = \opcf$ we have $\zfoneplus \simeq \{(1,1)\}$ and $\zfone \simeq \{(1,1), (3, -1)\}$. We may compute the geometric data attached to the fundamental cycle of $\uzfplus$ as $\tlambda = 2$ and $\ttt = 110$. The overflow ideal $\tildeD$ is the product $\goth{P}_5\goth{P}_{11}$ where $\goth{P}_5$ and $\goth{P}_{11}$ are the unique prime ideals of norms $5$ and $11$ respectively. In particular we have that $\tildeD$ is coprime to $\goth{f}$. Let us fix the smoothing ideal $\goth{a}$ as the unique prime ideal of norm $N = 17$ in $\ok$. The ideal $qN/\goth{a}$ is generated by the admissible base point $h = -96\zexc^2 - 88\zexc - 92$. The corresponding level is $\ell= 2\cdot5\cdot11\cdot17 = 1870$ and we have $t = 2$. Let us set $\tau = (23\zexc^2 + 29\zexc - 3438)/1870$ and $\sigma = (3\zexc^2 - \zexc - 898)/1870$. We may compute for $\conjko = 1, 3$ and $\goth{b} = \ok$:
$$
u_{\conjko} = \frac{\Gamma\left(\frac{\conjko}{4}, \tau, \sigma\right)^{17}}{\Gamma\left(\frac{17\conjko}{4}, 17\tau, 17\sigma\right)} \cdot \frac{\Gamma\left(\frac{\conjko}{4}-\frac{\tau+\sigma}{2}, \tau, \sigma\right)^{17}}{\Gamma\left(\frac{17\conjko}{4}-17\frac{\tau + \sigma}{2}, 17\tau, 17\sigma\right)} 
$$
The two values $u_1 \approx -959.8089713...+i\cdot92.0181162...$ and $u_3 \approx (-103.2385024...-i\cdot9.8976075...)\cdot10^{-5}$ coincide up to 1000 digits with two roots of the reciprocal polynomial 
$$P_{\mathrm{abs}} = x^6 + 1920x^5 + 930435x^4 + 357120x^3 + 930435x^2 + 1920x + 1$$ 
which defines an absolute equation of $\classfieldf$ over $\qq$. We may also check formula \refp{klfconjecture} as
$$ 17\cdot \zeta'_{\goth{f}}(\ok, 0) - \zeta'_{\goth{f}}(\goth{a}, 0) = \log\left|u_1\right|^2 \approx 13.7426178...$$
and similarly for the other class.

\subsection{Quartic examples}\label{quarticexamples}

We now present three examples for quartic fields. We fix the usual ordering $\goth{S}_2 = \{\mathrm{Id}, (21)\}$ of the symmetric group on two elements to simplify notations. We start with an example where $q = 5$ and $\zfoneplus \simeq \zsz{2}$. This illustrates one of the simplest situations outside of optimal settings for quartic fields. Then we present an example with $q = 67$ and $\zfoneplus \simeq \zsz{22}$. This showcases the benefits of using the unit group $\uzfplus$ instead of $\opcf$ to handle large values of $q$ in non-optimal settings. The third example is one of the simplest cases where the quartic field $\kk$ contains a real quadratic field and some of the parameters $\talpha_j$ belong to a real subfield of $\kk$. We use this example to illustrate the modifications in the computations that we make to compute higher elliptic units when $\kk$ admits a proper totally real subfield.
 
\stepcounter{example}
\subsubsection{Example \exc}\label{kappa2}

We first discuss a quartic example with $q = 5$ and $\zfoneplus \simeq \zsz{2}$. Let $\zexc$ be the complex root of the polynomial $x^4 - 3x + 1$ lying in the upper half-plane. Then $\kk = \qq(\zexc)$ has class number $1$. We choose $\goth{f} = \goth{P}_5$ the unique integral ideal of norm $q = 5$ in $\ok$. The corresponding narrow ray class group is $\classgroupf \simeq \zsz{2}$. It is isomorphic to $\clzfplus$ and it is represented by the ideals $\ok, 2\ok$. Let us fix a set of fundamental units
$$\eps_1 = \zexc^3 - 2\zexc^2 + \zexc,~\eps_2 =\zexc^3 + 2\zexc^2 + 3\zexc + 1$$
for the group $\opck = \uzfplus$. They satisfy $\eps_1 \equiv 1 \mod \goth{f}$ and $\eps_2 \equiv -1 \mod \goth{f}$, thus $\opcf$ is generated by $\eps_1$ and $\eps_2^2$. We may identify $\zfoneplus \simeq \{(1,1), (4, \eps_2)\}$ while 
$\zfone \simeq \{(1, 1), (2, -\xi), (3, \xi), (4, \eps_2)\}$
where $\xi = \zexc^3+\zexc^2+x-2$ has negative norm.

The contents associated to the fundamental units $\eps_1, \eps_2$ are $\tlambda_1 = \tlambda_2 = 1$ and the overflows are $\ttt_1 = 23$, $\ttt_2 = 11$. The associated overflow ideal $\tildeD_1$ is the unique integral ideal of norm $23$ in $\ok$ and $\tildeD_2$ is the unique ideal of norm $11$ in $\ok$. We may choose $\goth{b} = \ok$ and $\goth{a}$ the unique integral ideal of norm $N = 7$ in $\kk$. Both ideals $2qN\goth{D}_1/\goth{a}$ and $qN\goth{D}_2/\goth{a}$ are principal with admissible generators $h_1 = 10\zexc^3 + 110\zexc^2 + 160\zexc + 190$ and $h_2 = 5\zexc^3 + 20\zexc^2 + 10\zexc + 25$ which are compatible. The corresponding levels of the computations are $\ell_1 = 5\cdot7\cdot23 = 805$ and $\ell_2 = 5\cdot7\cdot11 = 385$. Let us write the parameters
\begin{align*}
&\tau = 11\zexc^3 + \zexc^2 - 2\zexc + 14898, && \tau' = - 4z^{3}+ 3z^{2}+ 6z+ 79  \\
&\sigma = -\zexc^3 + 2\zexc^2 - 4\zexc - 1737, && \sigma' = - z^{3}- 2z^{2}- 4z+ 567  \\
&\rho = -3\zexc^3 + 6\zexc^2 + 11\zexc + 7577, && \rho' = - z^{3}- 2z^{2}+ 7z+ 138
\end{align*}

Then the two products 
\begin{align*}
v_1 = &\prod_{k = 1, 4}\frac{G_2\left(\frac{-2k}{5}, \frac{\tau}{805}, \frac{\sigma}{805}, \frac{\rho}{805}\right)^{-7}}{G_2\left(\frac{-14k}{5}, \frac{\tau}{35}, \frac{\sigma}{35}, \frac{\rho}{35}\right)^{-1}}\times\frac{G_2\left(\frac{k}{5}, \frac{\tau'}{385}, \frac{\sigma'}{385}, \frac{\rho'}{385}\right)^{7}}{G_2\left(\frac{7k}{5}, \frac{\tau'}{35}, \frac{\sigma'}{35}, \frac{\rho'}{35}\right)}\\
v_2 = &\prod_{k = 2, 3}\frac{G_2\left(\frac{-2k}{5}, \frac{\tau}{805}, \frac{\sigma}{805}, \frac{\rho}{805}\right)^{7}}{G_2\left(\frac{-14k}{5}, \frac{\tau}{35}, \frac{\sigma}{35}, \frac{\rho}{35}\right)^{1}}\times\frac{G_2\left(\frac{k}{5}, \frac{\tau'}{385}, \frac{\sigma'}{385}, \frac{\rho'}{385}\right)^{-7}}{G_2\left(\frac{7k}{5}, \frac{\tau'}{35}, \frac{\sigma'}{35}, \frac{\rho'}{35}\right)^{-1}}
\end{align*}
coincide up to 1000 digits (that is $v_1 \approx v_2$) and are close to the root $\approx 4.3174875...-i\cdot 7.2167313...$ of the polynomial 
$$P_{\mathrm{abs}} = x^8 - 8x^7 + 66x^6 + 39x^5 + 49x^4 + 39x^3 + 66x^2 - 8x + 1$$
 which defines an absolute equation of $\classfieldf$ over $\qq$. The computation time for $1000$ digits is 46 seconds. These two products are obtained using two different possible choices for the orientation datum in the Main Conjecture \ref{conjecture}. Namely, for the computation of $v_1$ we used the orientation datum:
$$\bars{\mu} = 
(-1, 1), ~\bars{\nu} = (-1, 1), ~\bars{k} = (1, 1), ~\bars{\eps} = (1, 1) $$
whereas for the computation of $v_1$ we used the orientation datum:
$$\bars{\mu} = 
(-1, 1), ~\bars{\nu} = (-1, 1), ~\bars{k} = (3, 3), ~\bars{\eps} = (\xi, \xi).$$
This supports the idea that the computation of higher elliptic units is indeed independent of the choice of correct orientation datum. We then check formula \refp{klfconjecture} as:
$$7\cdot\zeta'_{\goth{f}}(\ok, 0) - \zeta'_{\goth{f}}(\goth{a}, 0) \approx \log\left|v_1\right|^2 \approx \log\left|v_2\right|^2 \approx 4.2587554...$$

\stepcounter{example}
\subsubsection{Example \exc}

Let us now present a quartic example with $q = 67$ and $\zfoneplus \simeq \zsz{22}$. Let $\zexc$ be the complex root of the polynomial $x^4 - 4x^3 - 6x^2 + 2x - 2$ lying in the upper half-plane. Then $\kk = \qq(\zexc)$ has class number $1$. We choose $\goth{f} = 67\ok+((x^3 - 4x^2 - 7x + 73)/3)\ok$ one of the two integral ideals of norm $q = 67$ in $\ok$. The corresponding narrow ray class group is $\classgroupf \simeq \zsz{6}$ while $\clzf \simeq \zsz{2}$. The class group $\classgroupf$ is then represented by $\ok, 2\ok, 4\ok, \goth{P}_2, 2\goth{P}_2, 4\goth{P}_2$ where $\goth{P}_2$ is the unique prime ideal of norm $2$ in $\ok$ and $\goth{P}_2^4 = 2\ok$. Let us fix a set of fundamental units
$$\eps_1 =\frac{\zexc^2 - 5\zexc + 1}{3},~\eps_2 = \frac{\zexc^2 - 2\zexc + 1}{3} $$
for the group $\opck = \uzfplus$. They satisfy $\eps_1 \equiv 14 \mod \goth{f}$ and $\eps_2 \equiv 42 \mod \goth{f}$, thus $\opcf$ is generated by $\eps_1^{11}$ and $\eps_1\eps_2^8$. We may identify $\zfone \simeq \zfoneplus \simeq \{(8^k \mod 67, \eps_2^k) \setseparator 0 \leq k \leq 21\}$.

The contents associated to the fundamental units $\eps_1, \eps_2$ are $\tlambda_1 = \tlambda_2 = 1$ and the overflows are $\ttt_1 = 1$, $\ttt_2 = 1$. The associated overflow ideals are therefore trivial. We shall focus on the case $\goth{b} = \ok$, $\conjko = 1, 2, 4$. We choose $\goth{a}$ to be the unique integral ideal of norm $N = 5$ in $\kk$. This gives in particular $\denomclasscohom = 5$. The integral ideal $4qN/\goth{a}$ is generated by the admissible base point $h_1= h_2 = (2412\zexc^3 - 10720\zexc^2 - 9112\zexc+ 4556)/3$. Here $m_1 = m_2 = 4$. The set of base points $(h_1, h_2)$ is compatible. The corresponding levels of the computations are $\ell_1 = \ell_2 = 5\cdot67 = 335$. Let us write the parameters
\begin{align*}
&\tau = \frac{2\zexc^3 - 7\zexc^2 - 16\zexc+ 423}{3}, && \tau' = \frac{-\zexc^2 + 5\zexc + 3056}{3}  \\
&\sigma = \frac{2\zexc^3 - 9\zexc^2 - 9\zexc + 5212}{3}, && \sigma' = \frac{-\zexc^3 + 6\zexc^2 - 3\zexc - 7190}{3}  \\
&\rho = \frac{ -\zexc^3 + 4\zexc^2 + 7\zexc - 1078}{3}, && \rho' = \frac{2\zexc^3 - 11\zexc^2 + 4\zexc+ 12647}{3}.
\end{align*}

Then for $\conjko = 1, 2, 4$ we compute:
$$I_{\conjko} =\prod_{k \in R}\frac{G_2\left(\frac{4k\cdot \kappa}{67}, \frac{\tau}{335}, \frac{\sigma}{335}, \frac{\rho}{335}\right)^{5}}{G_2\left(\frac{20k\cdot \kappa}{67}, \frac{\tau}{67}, \frac{\sigma}{67}, \frac{\rho}{67}\right)}\times\frac{G_2\left(\frac{-4k\cdot \kappa}{67}, \frac{\tau'}{335}, \frac{\sigma'}{335}, \frac{\rho'}{335}\right)^{-5}}{G_2\left(\frac{-20k\cdot \kappa}{67}, \frac{\tau'}{67}, \frac{\sigma'}{67}, \frac{\rho'}{67}\right)^{-1}}$$
where the product ranges over 
$$R := \{1, 3, 5, 8, 9, 14, 15, 22, 24, 25, 27, 40, 42, 43, 45, 52, 53, 58, 59, 62, 64, 66\}.$$
We have $I_1 \approx -35.6914222... - i\cdot9.6891118...$, $I_2 \approx 0.0020486... + 0.0019857...$ and $I_4 \approx 8.7046584... - 3.7478074...$. The values $\mathcal{V}_{\goth{f}, \conjko, \ok, \goth{a}} = I_{\conjko}^{5}$ and their inverses for $\conjko = 1, 2, 4$ coincide up to 1000 digits with six of the roots of a reciprocal polynomial:
$$P_{\mathrm{abs}} = x^{24} + 8040631607760x^{23} + 27979777748387337111319362x^{22} +... $$
which defines an absolute equation of $\classfieldf$ over $\qq$. The computation time for $1000$ digits is 26 minutes. 
 We then check formula \refp{klfconjecture} as:
$$5\cdot\zeta'_{\goth{f}}(\ok, 0) - \zeta'_{\goth{f}}(\goth{a}, 0) \approx \log\left|\mathcal{V}_{\goth{f}, 1, \ok, \goth{a}}\right|^2 \approx 7.2209270...$$

\stepcounter{example}
\subsubsection{Example \exc}\label{quarticreal}

This example showcases some subtle changes in the computations when the field we consider contains a real subfield. Let $\zexc$ be the complex root of the polynomial $x^4-12$ lying in the upper half-plane. Then $\kk = \qq(\zexc)$ has class number $1$ and \textit{contains the real quadratic field} $\qq(\sqrt{3})$. We choose the ideal $\goth{f}$ as the unique prime ideal of norm $q = 2$ in $\ok$. The corresponding narrow ray class group is $\classgroupf \simeq \zsz{2}$ while $\clzfplus$ is trivial. Let us choose the fundamental units:
$$ \eps_1 = \frac{\zexc^2 + 2\zexc + 2}{4},~~\eps_2 = \frac{-\zexc^2 + 4}{2} = 2+ \sqrt{3} $$
for $\opck = \uzfplus = \opcf$. Note that $\eps_2$ belongs to the subfield $\qq(\sqrt{3})$. We may compute the contents $\tlambda_1 = \tlambda_2 = 1$ as well as the overflows $\ttt_1 = \ttt_2 = 1$. The unit groups $\zfone$ and $\zfoneplus$ are both trivial. Let us fix $\goth{b} = (1)$ and $\goth{a} = 23\ok + \beta\ok$ an integral ideal of norm $N = 23$ in $\ok$ where $\beta = (\zexc^2 - 2\zexc - 26)/4$. The ideal $qN/\goth{a}$ is generated by the admissible base points $h_1 = h_2 = (-3\zexc^3 - 9\zexc^2 - 4\zexc + 34)/2$. These are compatible in the sense of Definition \ref{defcompatible}.
The parameters we obtain are:
\begin{align*}
&\tau = 177 + \sqrt{3}, && \tau' = \frac{z^2 - 2z- 26}{4} \\
&\sigma = \frac{z^3 + 3z^2 + 4z - 1354}{4}, && \sigma' = \frac{-z^3 - z^2 + 4z - 114}{4} \\
&\rho = \frac{-z^2 - 2z + 498}{4}, && \rho' = 177 + \sqrt{3}
\end{align*}
In particular $\tau, \rho' \in \rr$, $\Im(\tau'), \Im(\rho) < 0$ and $\Im(\sigma), \Im(\sigma') > 0$. In this special setting, we must use formula \refp{grreal} from section \ref{computinggr} instead of the usual formula \refp{sumcomputations} to perform the computations of the product:
$$\frac{G_{2, \goth{f}, \goth{b}, \goth{a}}^{+}(\eps_1, \eps_1\eps_2; h_2)}{G_{2, \goth{f}, \goth{b}, \goth{a}}^{-}(\eps_2, \eps_1\eps_2; h_1)} = \frac{G_2\left(\frac{1}{2}, \frac{\tau}{46}, \frac{\sigma}{46}, \frac{\rho}{46}\right)^{23}}{G_2\left(\frac{23}{2}, \frac{\tau}{2}, \frac{\sigma}{2}, \frac{\rho}{2}\right)} \times \frac{G_2\left(\frac{-1}{2}, \frac{\tau'}{46}, \frac{\sigma'}{46}, \frac{\rho'}{46}\right)^{-23}}{G_2\left(\frac{-23}{2}, \frac{\tau'}{2}, \frac{\sigma'}{2}, \frac{\rho'}{2}\right)}.
$$
The complex number we obtain coincides up to 1000 digits of precision with the root $\approx 13.9102308... -i\cdot24.0932265...$ of the polynomial 
$$P_{\mathrm{abs}} = x^8 - 28x^7 + 778x^6 - 112x^5 - 749x^4 - 112x^3 + 778x^2 - 28x + 1$$
 which defines an absolute equation of $\classfieldf$ over $\qq$. Note that in [\hspace{1sp}\cite{thesis}, section IV.2.6.3] we treated this example with a different choice of fundamental units for which the associated parameters $\tau, \sigma, \rho, \tau', \sigma', \rho'$ all lie outside of $\rr$. This means that the presence of a subfield does not automatically mean we have to use \refp{grreal} in the computations, but care is needed when treating these examples.

It is interesting to note that previous work on the determination of Stark units in the case of a quartic field $\kk$ with exactly one complex place containing a real quadratic subfield $\ff$ was carried out by Charollois and Darmon \cite{CharolloisDarmon} using multiplicative $1$-cocycles for the Hilbert modular group $\mathrm{SL}_2(\mathcal{O}_{\ff})$. In our construction, the presence of a real subfield is often times an obstacle to the computation of these same units using $2$-cocycles for $\slnz{4}$, showing that the two approaches complement eachother.

\stepcounter{example}
\subsection{A quintic example}\label{quinticexample}

In this section we present a quintic example with $q = 47$ and $\zfoneplus$ contains $23$ elements. This is one of the simplest example of a computation of higher elliptic unit above a quintic ATR field in a non-optimal setting. Let $\zexc$ be the complex root of the polynomial $x^5 - 2x^4 - 2x^3 + 5x^2 + 2x - 1$ lying in the upper half-plane. The number field $\kk = \qq(\zexc)$ is ATR and has narrow class number $1$. We choose $\goth{f}$ to be the unique prime ideal of norm $q = 47$ in $\ok$. The corresponding narrow ray class group is $\classgroupf \simeq \clzfplus \simeq \zsz{2}$. It is represented by the ideals $\ok$ and $5\ok$. The unit group $\opck = \uzfplus$ admits the fundamental units:
$$\eps_1 = z^2, ~~ \eps_2 = z^2 +2z+1, ~~ \eps_3 = -2z^4 + 6z^3 - z^2 - 11z + 6.$$
The remainders modulo $\goth{f}$ are $\eps_1 \equiv 8 \mod \goth{f}$, $\eps_2 \equiv 37 \mod \goth{f}$ and $\eps_3 \equiv 12 \mod \goth{f}$. We compute $\zfoneplus \simeq \{ (6^k \mod 47, \eps_1^{k}) \setseparator 0 \leq k \leq 22 \}$ and 
$$\zfone \simeq \zfoneplus \bigcup \{(q-k, -\eps) \setseparator (k, \eps) \in \zfoneplus\}.$$
Let us fix the ordering $\{ \mathrm{Id}, (32), (21), (231), (312), (31)\}$ of $\goth{S}_3$. Computing contents and overflows gives $\tlambda_1 = \dots = \tlambda_6 = 1$ and:
\begin{align*}
\tildeD_1 &= (-2z^4 + 5z^3 + z^2 - 10z - 2)\ok && \ttt_1 = \ts_1 = 883\\
\tildeD_2 &= (z^4 - 2z^3 - z^2 + z + 3)\ok && \ttt_2 = \ts_2 = 521\\
\tildeD_3 &= (-2z^4 + 6z^3 - 3z^2 - 5z + 6)\ok && \ttt_3 = \ts_3 = 7 \cdot 1051\\
\tildeD_4 &= (-z^4 + 3z^3 - 2z^2 - 2z + 4)\ok && \ttt_4 = \ts_4 = 7 \cdot 97 \\
\tildeD_5 &= (-z^4 + 3z^3 - z^2 - 5z + 3)\ok && \ttt_5 = \ts_5 = 7 \\
\tildeD_6 &= (z^3 - 4z^2 + z + 15)\ok && \ttt_6 = \ts_6 = 11 \cdot 31 \cdot 233.
\end{align*}
We may choose $\goth{a}$ to be the unique prime ideal of norm $N = 17$ in $ \ok$. The ideals $qN\tildeD_j/\goth{a}$ for $1 \leq j \leq 6$ admit the following admissible generators
\begin{align*}
h_1 & = -517z^4 + 1410z^3 + 1316z^2 - 4559z + 1410, \\
h_2 & = -1363z^4 + 1974z^3 + 4559z^2 - 5264z - 6016, \\
h_3 & = 423z^4 - 1081z^3 + 1974z^2 - 846z - 2679, \\
h_4 & = 1833z^4 - 3619z^3 - 1833z^2 + 6721z - 1222, \\
h_5 & = 3243z^4 - 5358z^3 - 8836z^2 + 14288z + 11421, \\
h_6 & = 2115z^4 - 7003z^3 + 1880z^2 + 9353z - 611.
\end{align*}
The levels $\ell_j = 47 \cdot 17 \cdot \ttt_j$ in the computation are then given by:
$$\ell_1 = 705517, ~\ell_2 = 416279, ~\ell_3 = 5878243, ~\ell_4 = 542521, ~\ell_5 = 5593, ~\ell_6 = 63482947.$$
Let us set for $j = 1, \dots, 6$ the signs $\epsilon_1 = \epsilon_4 = \epsilon_6 = -1$ and $\epsilon_2 = \epsilon_3 = \epsilon_5 = 1$. We compute the complex numbers:
$$ v_j = \prod_{k \in R} \frac{G_3\left(\frac{k\cdot \epsilon_j}{47}, \frac{\tau_j}{\ell_j}, \frac{\sigma_j}{\ell_j}, \frac{\rho_j}{\ell_j}, \frac{\varpi_j}{\ell_j}\right)^{17}}{G_3\left(\frac{17k\cdot\epsilon_j}{47}, \frac{17\tau_j}{\ell_j}, \frac{17\sigma_j}{\ell_j}, \frac{17\rho_j}{\ell_j}, \frac{17\varpi_j}{\ell_j}\right)}
$$
where 
$$R = \{1, 2, 3, 4, 6, 7, 8, 9, 12, 14, 16, 17, 18, 21, 24, 25, 27, 28, 32, 34, 36, 37, 42\} \simeq \zfoneplus$$
and the parameters $\tau_j, \sigma_j, \rho_j, \varpi_j$ are given in Table \ref{tablequintic} below. The higher elliptic units we compute are:
$$\mathcal{V}_{\goth{f}, 1, \ok, \goth{a}} = \mathcal{V}_{\goth{f}, 5, \ok, \goth{a}}^{-1} = \frac{v_1v_2v_4v_5}{v_3v_6} \approx (-1.0581949...-i\cdot 9.6144804...)\cdot 10^{17}.$$
Both complex numbers $\mathcal{V}_{\goth{f}, 1, \ok, \goth{a}}$ and $\mathcal{V}_{\goth{f}, 5, \ok, \goth{a}}$ coincide up to 1000 digits with two of the roots of the reciprocal polynomial 
\begin{align*}
P_{\mathrm{abs}} & = ~x^{10} + 211638998633678735x^9 \\
& + 935580113532604138369179331868118360x^8 \\
& + 63162478400025151879460822014371111x^7 \\
& + 342193760469590788211407283208696885x^6 \\
& + 1253281261201662931852705831834834479x^5 + ...
\end{align*}
which defines an absolute equation of $ \classfieldf$ over $\qq$. The computation for $100$ digits already takes around 30 minutes and the computation time for $1000$ is about 1 day. We may check the expected Kronecker limit formula \refp{klfconjecture} as:
$$17 \cdot \zetaf'(\ok, 0) - \zetaf'(\goth{a}, 0) \approx \log|\mathcal{V}_{\goth{f}, 1, \ok, \goth{a}}|^2 \approx 82.8264748... $$
\bigskip

\begin{table}[H]
\caption{Parameters for the quintic example}
\begin{center}
\noindent \begin{tabular}{| C C R C R C R C R C R |}
\hline
\tau_{1} =  &  & 103z^{4} & - & 139z^{3} & - & 365z^{2} & + & 269z & - & 436584\\
\sigma_{1} =  & - & 7z^{4} & - & 102z^{3} & + & 342z^{2} & - & 44z & - & 2256682\\
\rho_{1} =  &  & 36z^{4} & + & 20z^{3} & - & 119z^{2} & - & 26z & + & 536379\\
\varpi_{1} =  & - & 83z^{4} & + & 52z^{3} & + & 397z^{2} & + & 109z & + & 1389464\\
\hline
\tau_{2} =  & - & 201z^{4} & + & 388z^{3} & + & 372z^{2} & - & 878z & + & 1200066\\
\sigma_{2} =  &  & 14z^{4} & + & 30z^{3} & - & 127z^{2} & - & 84z & + & 188803\\
\rho_{2} =  &  & 329z^{4} & - & 858z^{3} & - & 119z^{2} & + & 1673z & - & 2110276\\
\varpi_{2} =  & - & 44z^{4} & + & 129z^{3} & + & 27z^{2} & - & 257z & + & 143983\\
\hline
\tau_{3} =  & - & 108z^{4} & + & 525z^{3} & - & 673z^{2} & - & 45z & - & 10099409\\
\sigma_{3} =  & - & 1153z^{4} & + & 1722z^{3} & + & 3510z^{2} & - & 4772z & - & 7003961\\
\rho_{3} =  &  & 182z^{4} & - & 476z^{3} & - & 637z^{2} & + & 1302z & + & 225932\\
\varpi_{3} =  & - & 383z^{4} & + & 840z^{3} & + & 815z^{2} & - & 3225z & + & 25426283\\
\hline
\tau_{4} =  &  & 27z^{4} & - & 112z^{3} & + & 86z^{2} & + & 76z & - & 10031781\\
\sigma_{4} =  &  & 66z^{4} & + & 28z^{3} & - & 167z^{2} & - & 116z & + & 175233784\\
\rho_{4} =  &  & 28z^{4} & - & 91z^{3} & - & 112z^{2} & + & 280z & - & 56622104\\
\varpi_{4} =  &  & 3z^{4} & + & 63z^{3} & + & 85z^{2} & - & 67z & - & 119853183\\
\hline
\tau_{5} =  &  & 8z^{4} & - & 34z^{3} & + & 29z^{2} & + & 50z & - & 16389\\
\sigma_{5} =  & - & 16z^{4} & + & 33z^{3} & + & 26z^{2} & - & 72z & - & 11413\\
\rho_{5} =  & - & z^{4} & + & 27z^{3} & - & 57z^{2} & - & 29z & + & 27841\\
\varpi_{5} =  &  & z^{4} & - & 6z^{3} & + & 8z^{2} & + & 8z & - & 4174\\
\hline
\tau_{6} =  & - & 3856z^{4} & + & 14280z^{3} & - & 14963z^{2} & - & 222z & + & 774171489\\
\sigma_{6} =  & - & 6287z^{4} & + & 16854z^{3} & + & 1669z^{2} & - & 34113z & + & 229255403\\
\rho_{6} =  &  & 1475z^{4} & - & 7935z^{3} & + & 10401z^{2} & + & 6967z & - & 209402683\\
\varpi_{6} =  & - & 705z^{4} & + & 2446z^{3} & - & 3086z^{2} & - & 906z & + & 39086064\\
\hline
\end{tabular} 
\label{tablequintic}
\end{center}
\end{table}

\bibliographystyle{alpha}
\bibliography{bibliographie}{}

\end{document}